\documentclass[12pt, article]{amsart}
\usepackage[english]{babel}
\usepackage{amsmath,amsthm,amssymb,amsfonts, enumitem}
\usepackage{extarrows}
\usepackage{fullpage}
\usepackage{graphicx}
\usepackage[colorlinks=true, allcolors=blue,backref=page]{hyperref}
\usepackage{soul}
\usepackage{tikz-cd}
\usepackage{tcolorbox}
\usepackage[all]{xy}
\theoremstyle{plain}
  
  \newtheorem{theorem}{Theorem}[section] % Numbering based on section
  \newtheorem{proposition}[theorem]{Proposition}   % Shares numbering with theorems
  \newtheorem{lemma}[theorem]{Lemma} % Shares numbering with theorems
  \newtheorem{corollary}[theorem]{Corollary} % Shares numbering with theorems

\theoremstyle{definition} 
  \newtheorem{notation}[theorem]{Notation}

  \newtheorem{example}[theorem]{Example}
  \newtheorem{definition}[theorem]{Definition}
  \newtheorem{remark}[theorem]{Remark}

  \newtheorem{construction}[theorem]{Construction}

\newcommand{\bb}[1]{\mathbb{#1}}
\newcommand{\cal}[1]{\mathcal{#1}}
\newcommand{\CC}{\mathbb{C}}
\newcommand{\RR}{\mathbb{R}}
\newcommand{\ZZ}{\mathbb{Z}}

\newcommand{\PP}{\mathbb{P}}

\newcommand{\RO}{\RR^E_{\geq 0}}

\newcommand{\cone}{\operatorname{cone}}
\newcommand{\conv}{\operatorname{Conv}}
\newcommand{\FM}{\operatorname{FM}}

\newcommand{\pe}{\underline{\Pi}_E}
\newcommand{\ste}{\Pi_E}
\newcommand{\Bl}{\operatorname{Bl}}
\newcommand{\codim}{\operatorname{codim}}
\newcommand{\sslash}{\mathbin{/\mkern-6mu/}}
\newcommand{\quotient}{ / \! \! /}

\setitemize{itemsep=1em,topsep = 1em}
\setenumerate{itemsep=1em,topsep = 1em}

\begin{document}

\title{
Polymatroids and moduli of points in flags
}

\author{Patricio Gallardo
}
\address{
{\small Department of Mathematics,
University of California, Riverside,
900 University Ave.
Riverside, CA 92521
Skye Hall}
}
\email{pgallard@ucr.edu}

\author{Javier Gonz\'alez Anaya}
\address{
{\small Department of Mathematics, 
Harvey Mudd College, 
Claremont, CA 91711}
}
\email{javiergo@hmc.edu}

\author{Jose Luis Gonzalez}
\address{
{\small Department of Mathematics,
University of California, Riverside,
900 University Ave.
Riverside, CA 92521
Skye Hall}
}
\email{jose.gonzalez@ucr.edu}

\begin{abstract} 
We introduce and study different compactifications of the moduli space of $n$ distinct weighted labeled points in a flag of affine spaces. We construct these spaces via the weighted and generalized Fulton-MacPherson compactifications of Routis and Kim-Sato. For certain weights, our compactifications are toric and isomorphic to the polypermutohedral and polystellahedral varieties, which arise in the work of Crowley-Huh-Larson-Simpson-Wang and Eur-Larson on polymatroids, a generalization of matroids. Moreover, we show that these toric compactifications have a fibration structure, with fibers isomorphic to the Losev-Manin space, and are related to each other via a geometric quotient.
\end{abstract}

\maketitle

\section{Introduction}
Our work presents novel morphisms between two types of seemingly unrelated spaces, one originating in moduli theory, and the other in combinatorics, specifically in the theory of polymatroids. 
From the moduli perspective, we study generalizations of the Hassett spaces of weighted pointed stable rational curves. On the combinatorial side, we encounter as special cases the polypermutohedral and polystellahedral varieties of \cite{crowley2022bergman,eur2024intersection}, which arise from polytopes that generalize the permutohedron and stellahedron.

We begin by describing the moduli-theoretic side.\! Recall that the Fulton-MacPherson~(FM) compactification $X[n]$, introduced in 1994, is a smooth normal crossings compactification of the configuration space of $n$ distinct labeled points in a smooth variety $X$ \cite{fulton1994compactification}. It has since been widely used in moduli theory; see, for example, \cite{chen2009pointed,gallardo2017wonderful,pandharipande1995geometric,kim2014compactification}. 
In 2014, Routis \cite{routis2014weighted} generalized $X[n]$ to a weighted version $X_{\mathbf{w}}[n]$, allowing certain coincidences between points prescribed by a weight vector $\mathbf{w} \in \cal{D}^{\FM}_n$ (Definition~\ref{def:weights}). The moduli spaces constructed via the FM compactifications in the previous examples, as well as those we construct here, inherit a ``universal" flat family from the FM space $X[n]^{+} \rightarrow X[n]$. The fibers of the family are either the configuration of distinct points on $X$, or degenerations of $X$ marked by $n$ points; see Section \ref{sec:FMpreliminary}. 

Given an $n$-tuple of positive integers $\mathbf{a}=(a_1,\dots,a_n)$ with $a_{i+1} \leq a_{i}$, we use the weighted Fulton-MacPherson space to construct compactifications 
$\mathbb{C}_\mathbf{w}^{[\mathbf{a}]}$ of the open locus $\mathcal{U}^{\mathbf{a}} \subsetneq \mathbb{C}^{\mathbf{a}}:=\prod\bb{C}^{a_i}$ parametrizing $n$ distinct weighted labeled points $p_i$ within a flag 
$ 
\CC^{a_n} \subseteq \CC^{a_{n-1}} \subseteq \cdots \subseteq \CC^{a_1},
$ 
where $p_i \in \mathbb{C}^{a_i}$ for all $i$. 
We then leverage $\bb{C}^{[\mathbf{a}]}_\mathbf{w}$ to define a compactification $T^{\mathbf{a}}_{\mathbf{w}}$ of the moduli space of $n$ weighted points in the flag, up to translation along the $\mathbb{C}^{a_n}$-coordinates and global scaling by a $\mathbb{G}_m$ factor; see Section \ref{subsec: taw}. 
We say the above compactifications are \emph{geometric} because they inherit a ``universal" family from the FM space. 
When $\mathbf{a}=(d,\dots,d)$, the space $T_\mathbf{w}^\mathbf{a}$ recovers the moduli spaces of weighted pointed trees of projective spaces from \cite{chen2009pointed,gallardo2017wonderful} as special cases. In particular, when $d=1$, $T_\mathbf{w}^\mathbf{a}$ specializes to the moduli space of weighted pointed stable rational curves \cite{hassett2003moduli}.

On the combinatorial side, we have the polypermutohedral varieties, introduced in 2022 in \cite{crowley2022bergman}. These varieties arise from studying the intersection theory of polymatroids, which are combinatorial abstractions of subspace arrangements, in the same way that matroids are combinatorial abstractions of hyperplane arrangements \cite{eur2024intersection}. The inner normal fan of the polypermutohedron coincides with the Bergman fan of its associated polymatroid, and these varieties play a crucial role in establishing the K\"ahler package for the Chow ring of polymatroids. While their construction is not initially related to moduli theory, they can be seen as a generalization of the permutohedral variety.

In our first result, we relate the above algebro-geometric and combinatorial constructions. We generalize the well-known fact that the permutohedral variety is isomorphic to the toric compactification of a moduli space of weighted points known as the Losev-Manin compactification \cite{losev2000new}. Its proof is given in Section \ref{sec:Polypermuta_moduli}.  

\begin{theorem} \label{GFM: thm:weightedFM_forFlags}
    Consider an $n$-tuple $\mathbf{a} = (a_1, \dots, a_n) \in \mathbb{Z}_{>0}^n$ with $a_{i+1} \leq a_{i}$, and a weight vector $\mathbf{w} \in \cal{D}^{\operatorname{T}}_n$ (Definition~\ref{pp: twa weights}). There exists a smooth, normal crossings, geometric compactification $T^{\mathbf{a}}_{\mathbf{w}}$ of the moduli space of $n$ distinct labeled points in the flag 
    $\CC^{a_n} \subseteq \CC^{a_{n-1}} \subseteq \cdots \subseteq \CC^{a_1}$, up to scaling and translations preserving the flag, such that:
    \begin{itemize}
        \item[(i). ] If $\mathbf{w} = (\varepsilon_1, \ldots, \varepsilon_{n-1}, 1)$ with $\sum_{i=1}^{n-1}\varepsilon_i\leq 1$, then $T^{\mathbf{a}}_{LM}:=T_{\mathbf{w}}^{\mathbf{a}}$ is isomorphic to the polypermutohedral toric variety associated to $\mathbf{a}$. We refer to this as the Losev-Manin compactification of this moduli problem.
        \item[(ii). ] Let $T^\mathbf{a}$ denote the compactification corresponding to the weight $\mathbf{w}=(1,\dots,1)$. Just as in Kapranov's construction of $\overline{M}_{0,n}$, there is a sequence of smooth blow-ups
        \[
        T^{\mathbf{a}} \longrightarrow 
        T^{\mathbf{a}}_{LM} 
        \longrightarrow 
        \mathbb{P}^{a_1 + \ldots + a_{(n-1)}-1}.
        \]
        \item[(iii). ] $T^{\mathbf{a}}_{LM}$ is a nontrivial, locally trivial fibration over $\prod_{i=1}^{n-1} \mathbb{P}^{a_i-1}$, with fiber isomorphic to the standard $(n-2)$-dimensional Losev-Manin moduli compactification $\overline{M}_{0,n+1}^{LM}$.
    \end{itemize}
\end{theorem}

Our next result is motivated by the following moduli considerations: If $T^{\mathbf{a}}_{LM}$ is a toric compactification of the moduli space of $n$ points in the flag $\mathbb{C}^{a_n} \subseteq \mathbb{C}^{a_{n-1}} \subseteq \cdots \subseteq \mathbb{C}^{a_1}$, up to translation and scaling, then we expect that there exists another space parametrizing configurations of points without considering equivalence classes. One would anticipate that this configuration space recovers $T^{\mathbf{a}}_{LM}$ as a quotient, since this behavior is observed for the moduli space of points in the line \cite[Theorem~3.4]{hu2000mori}. 
Moreover, one would expect such space to be modular and admit a birational map to the product of projective spaces $\prod_{i=1}^n \mathbb{P}^{a_i}$, as this latter space is a natural compactification of $n$ points in the flag $\mathbb{C}^{a_n} \subseteq \cdots \subseteq \mathbb{C}^{a_1}$.

The following theorem shows that all our expectations are realized for arbitrary dimensions and numbers of points, and its proof is given in Section~\ref{sec:FM_polystella}. On the geometric side, we construct the moduli spaces $\bb{P}^{[\mathbf{a}]}_H$ via the generalized Fulton-MacPherson compactification, introduced in \cite{kim2009generalization} and denoted $X^{[n]}_D$. In this generalization, one fixes a nonsingular proper subvariety $D \subseteq X$ and parametrizes configurations of points away from $D$. As in the classical case, $X^{[n]}_D$ carries a ``universal" flat family parametrizing degenerations of $X$ marked by $n$ points, which are always away from $D$ but not necessarily distinct. On the combinatorial side, we encounter the polystellahedral variety $PS_{\mathbf{a}}$ of an $n$-tuple $\mathbf{a}$; this variety plays a central role in constructing the augmented Chow rings of polymatroids \cite{eur2024intersection}. As in our first theorem, the geometric and combinatorial constructions are independent of each other, and a priori, there was no reason to expect they were the same.

\begin{theorem}
\label{thm:mainPolyStella}
Consider an $n$-tuple $\mathbf{a} = (a_1, \ldots, a_n) \in \mathbb{Z}_{>0}^n$ with $a_{i+1} \leq a_{i}$. There exists a smooth, normal crossings, geometric compactification $\bb{P}^{[\mathbf{a}]}_H$ of the configuration space of $n$ not necessarily distinct points in the flag $\CC^{a_n} \subseteq \CC^{a_{n-1}} \subseteq \cdots \subseteq \CC^{a_1}$ such that:
    \begin{itemize}
    \item[(i). ] The variety $\bb{P}^{[\mathbf{a}]}_H$ is constructed as an iterated blow-up of $\prod_{i=1}^n \mathbb{P}^{a_i}$ along torus invariant subvarieties. In particular, it is a toric variety itself. 

    \item[(ii). ] $\bb{P}^{[\mathbf{a}]}_H$ is isomorphic to the polystellahedral variety associated to $\mathbf{a}$.

    \item[(iii). ] There exist an open $\left(\bb{P}^{[\mathbf{a}]}_H\right)^{\circ}\subsetneq \bb{P}^{[\mathbf{a}]}_H$ and a geometric quotient such that 
        \[
            T^{\mathbf{a}}_{LM}  \cong
             \left(\bb{P}^{[\mathbf{a}]}_H\right)^{\circ}\quotient \mathbb{G}_m.
        \]
\end{itemize}
\end{theorem}

Along the way, we study pullbacks and refinements of combinatorial building sets (Sections \ref{sec:FibrationPP} and \ref{PP: subsec refinement}), produce associated toric fibrations (Section~\ref{sec:FibrationPP}), give combinatorial interpretations to some toric fans (Section~\ref{sec:PolyGeometry}), and obtain a quotient presentation of the exceptional divisor of a toric blow-up (Section~\ref{subsection.quotients}). These can all be of independent interest.  

We conclude with some remarks about our results. Theorems~\ref{GFM: thm:weightedFM_forFlags} and \ref{thm:mainPolyStella}, provide evidence that our moduli spaces are combinatorially rich objects that generalize both the permutohedral and stellahedral varieties. Moreover, by \cite[Proposition~2.2]{eur2024intersection}, given any $n$-tuple $\mathbf{a}$, there is a map $\bb{P}^{[\mathbf{1}^m]}_H\to\bb{P}_H^{[\mathbf{a}^-]}$, where $\mathbf{a}^- = (a_1,\dots,a_{n-1})$,
$m =  \sum_{i=1}^{n-1} a_i$, and $\mathbf{1}^m=(1,\dots,1)\in\bb{Z}^m_{>0}$. Analogously, we prove in Proposition~\ref{PP: polypermuto refinement} that there is a map $\overline{M}_{0,m+2}^{LM}\to T_{LM}^\mathbf{a}$.
Currently, there is no moduli interpretation for these maps. Altogether, we can summarize the morphisms described thus far in the following diagram:
\[
    \begin{tikzcd}
        \bb{P}^{[\mathbf{1}^m]}_H\ar[d]\ar[r,dashed,shift left=.75ex,"/\!\!/\mathbb{G}_m", bend left=4]
        & \overline{M}_{0,m+2}^{LM}\ar[d]\ar[l,hook',shift left=.75ex, bend left=0] & \overline{M}_{0,m+2} \ar[l]\ar[d,dashed] \\
        \bb{P}^{[\mathbf{a}^-]}_H \ar[d]\ar[r,dashed,shift left=.75ex,"/ \! \! / \mathbb{G}_m", bend left=4] & T^{\mathbf{a}}_{LM} \ar[d]\ar[l,hook', shift left=.75ex, bend left=0] & T^{\mathbf{a}} \ar[l] \\
        \prod_{i=1}^{n-1} \mathbb{P}^{a_i} & \mathbb{P}^{a_1 + \cdots + a_{(n-1)} - 1} &.
    \end{tikzcd}
\]
Here, besides the two vertical arrows described above, the remaining vertical solid arrows, horizontal solid arrows in the second row, and dashed arrows are obtained from Theorems~\ref{GFM: thm:weightedFM_forFlags} and \ref{thm:mainPolyStella}. Lastly, the horizontal inclusion maps indicate the fact that $T_{LM}^\mathbf{a}$ embedds as a torus-invariant divisor inside $\bb{P}_H^{[\mathbf{a}^-]}$, as already presented in the discussion before Theorem 1.5 in \cite{eur2023stellahedral}. 

Many fundamental questions about the birational geometry of $\overline{M}_{0,n}$ naturally extend to the setting of $T^{\mathbf{a}}$. For example, the study of its birational geometry and when it is a Mori dream space---a topic that has drawn substantial attention in the algebraic geometry community. While the classification of values of $n$ for which $\overline{M}_{0,n}$ is a Mori dream space is now largely resolved (see \cite{castravet2015m,gonzalez2016some,hausen2018blowing}), it would be interesting to answer the analogous question of determining for which $n$-tuples $\mathbf{a}$ the space $T^{\mathbf{a}}$ is a Mori dream space; see for example \cite[Theorem~1.3]{ggg-higher-lm}.

\subsection*{Acknowledgements} 
The authors are grateful for the supportive environment provided by their home institutions, the University of California, Riverside (PG and JLG), and Harvey Mudd College (JGA). Jos\'e Gonz\'alez was supported by a grant from the Simons Foundation (Award Number 710443). 
Patricio Gallardo is partially supported by the National Science Foundation under Grant No. DMS-2316749.
The authors are grateful to Matt Larson for pointing out that the moduli spaces studied in \cite{ggg-higher-lm} are polypermutohedral varieties, and other helpful conversations. We also thank the Western Algebraic Geometry Symposium (WAGS) for providing the venue where this connection was originally pointed out to us.

\setcounter{tocdepth}{2}
\tableofcontents

\section{Preliminaries}

In this section we introduce the notation and preliminary results needed for our work.
Throughout this article we work over the field of complex numbers $\mathbb{C}$.

\subsection{Notation} Define $[n]:=\{1,\dots,n\}$ for all $n\in\bb{Z}_{>0}$. Given a finite set $E$, consider the lattice $\bb{Z}^E$ and the vector space $\bb{R}^E=\bb{Z}^E\otimes_\bb{Z}\bb{R}$, both with distinguished basis $\{e_i\,\vert\,i\in E\}$. If $I\subseteq E$, define
\[
    e_I := \sum_{i\in I}e_i\in\bb{Z}^E\subseteq\bb{R}^E.
\]
Given a point $\mathbf{x}=\sum_{i\in E}x_ie_i\in\bb{R}^E$ and a subset $S\subseteq E$, define $\mathbf{x}_S = \sum_{i\in S}x_i$.

We denote the equivalence class of $e_I\in\bb{R}^E$ in $\bb{R}^E/\bb{R}e_E$ by $\overline{e}_I$. The $E$-simplex is defined as 
\[
    \Delta_{E} := \conv\{\overline{e}_i\,\vert\, i\in E\}\subseteq\bb{R}^{E}/\bb{R}e_E.
\]
In particular, the $n$-dimensional simplex is $\Delta_n:=\Delta_{[n+1]}$. Its inner normal fan in $\bb{R}^E/\bb{R}e_E$, denoted $\Sigma_n$, is the collection of all cones of the form
\[
    \cone(\overline{e}_{i_1},\dots,\overline{e}_{i_k}\,\vert\, I=\{i_1,\dots,i_k\}\subsetneq [n+1]),
\]
as $I$ ranges over all proper subsets of $[n+1]$. The associated projective space to $E$, denoted as $\bb{P}^E$, is the toric variety constructed from the inner normal fan of $\Delta_E$ in $\bb{R}^E/\bb{R}e_E$, which we denote as $\Sigma_E$. 
Notice that if $|E|=n$, then $\bb{P}^E$ has dimension $n-1$.
A vector $\mathbf{x}=\sum_{i\in E} x_ie_i$ has homogeneous coordinates $[x_i\,\vert\, i\in E]\in\bb{P}^E$.

\begin{definition}
    Given a fan $\Sigma$ supported in $\bb{R}^E/\bb{R}e_E$, we denote its corresponding toric variety by $X(\Sigma)$. 
\end{definition}

Here we always consider the star subdivision of a fan $\Sigma$ along a cone $\tau\in\Sigma$ to be its barycentric subdivision along $\tau$. In other words, the subdivision of $\Sigma$ by the vector $v_\tau=v_1+\cdots+v_r$, where the $v_i$ are the primitive lattice vectors generating the rays of $\tau$. 

\subsection{Polymatroids, polypermutohedra and polystellahedra}
\label{sec:preliminaryPoly}
Next, we discuss the necessary combinatorial background on polymatroids, which are generalizations of matroids. Just as in the theory of matroids, polymatroids admit several equivalent definitions. Here we will present their definition in terms of rank functions, and then prove the equivalence of this presentation with that one using independence and base polytopes. From the latter perspective, polymatroids are precisely the integral generalized permutohedra from \cite{postnikov2009permutohedra}.
Polymatroids were originally studied by Edmonds in connection with combinatorial optimization problems \cite{edmonds1970}.
In recent years, there has been increased interest in extending central results from matroid theory to polymatroid theory. For example, many authors have focused on generalizing key concepts and properties, such as defining the Bergman fan of a polymatroid \cite{crowley2022bergman}, their augmented Chow rings \cite{eur2024intersection}, and even their Hodge theory \cite{pagaria2023hodge}. Additionally, polymatroids have found applications in the realm of commutative algebra and algebraic geometry; see for example \cite{castillo2020multidegrees}.

\begin{definition}  
A function $f:2^E\to\bb{R}$ is called a \emph{set function}. A set function $f$ is:
    \begin{itemize}
        \item \emph{Normalized}, if $f(\emptyset)=0$.
        \item \emph{Monotone} or \emph{nondecreasing}, if $f(S)\leq f(T)$ whenever $S\subseteq T\subseteq E$.
        \item \emph{Submodular}, if $f(I\cup J) + f(I\cap J)\leq f(I) + f(J)$.
    \end{itemize}
\end{definition}

\begin{definition}[{cf. \cite[Definition~3.8]{ferroni2022matroids}}]
    Let $f:2^E\to\RR$ be a normalized submodular set function.
    \begin{enumerate}
        \item\label{def: generalized perm} We denote the \emph{generalized permutohedron} associated to $f$ as
        \[
            GP(f) = \left\{\mathbf{x}\in\RR^E\,\vert\,\mathbf{x}_E=f(E)\text{ and }\mathbf{x}_S\leq f(S)\text{ for all }S\subseteq E\right\}.
        \]
        \item If $f$ is further assumed to be monotone, then we refer to its corresponding generalized permutohedron as the \emph{base polytope} corresponding to $f$. We denote it by $B(f)$.
        \item Suppose that $f$ is monotone, then we define the \emph{independence polytope} corresponding to $f$ as
        \[
            I(f) = \left\{\mathbf{x}\in\RO\,\vert\,\mathbf{x}_S\leq f(S)\text{ for all }S\subseteq E\right\}.
        \]
    \end{enumerate}
    In particular, note that $B(f)$ is the face of $I(f)$ whose supporting hyperplane has normal vector $e_E=(1,\dots,1)\in\bb{R}^E$.
\end{definition}
\begin{definition}[{cf. \cite[Definition~1.1]{eur2024intersection}}] \label{bg: def cages}
    A \emph{cage} over $[n]$ is an $n$-tuple $\mathbf{a}=(a_1,\dots,a_n)$ of nonnegative integers. If $a_1=\cdots=a_n$ we say the cage is \emph{constant}. A \emph{caging} of a finite set $E$ is a surjective map of finite sets $\pi:A\to E$. We assume without loss of generality that $E= [n] = \{ 1, \ldots, n \}$, in which case the cage corresponding to this caging is $\mathbf{a}=(a_1,\dots,a_{n})$, where $a_i=|\pi^{-1}(i)|$. 

     For simplicity, in the sequel we shall consider only cages $\mathbf{a}$ over $[n]$ such that $a_1\geq\cdots\geq a_n$. This constitutes no loss of generality and simplifies the presentation significantly.
\end{definition}
\begin{definition}
    A \emph{polymatroid} $P=(E,f)$ is a pair of a finite set $E$ and a normalized, submodular, monotone set function $f:2^E\to\bb{R}$. The function $f$ is often called the \emph{rank function} of the polymatroid.
    \begin{itemize}
        \item The \emph{independence polytope} of $P$ is $I(P):=I(f)\subseteq\bb{R}^E_{\geq 0}$.
        \item The \emph{base polytope} of $P$ is $B(P):=B(f)\subseteq\bb{R}^E_{\geq 0}$. It is the face of $I(P)$ defined by the vector $e_E=(1,\dots,1)\in\bb{R}^E_{\geq 0}$. 
    \end{itemize}
    We say that the polymatroid has cage $\mathbf{a}=(a_1,\dots,a_n)$ if $f(i)\leq a_i$ for all $i\in[n]$.
    In particular, a \emph{matroid} is a polymatroid with cage $\mathbf{a}=(1,\dots,1)$.
\end{definition}

A polymatroid $P$ over a set $E$ is equivalently defined by its rank function, its independence polytope or its base polytope. In order to see this let us recall the following characterization result. For vectors $\mathbf{u},\mathbf{v}\in\RO$, write $\mathbf{u} \geq \mathbf{v}$ if $\mathbf{u}-\mathbf{v}\in\RO$.

\begin{proposition}[{\cite[Theorem~44.5]{schrijver2003combinatorial}}]\label{prop: characterization of polymatroids}
     A nonempty polytope $Q\subseteq\RO$ is the independence polytope of a polymatroid if and only if:
    \begin{enumerate}
        \item Given $\mathbf{v}\in\RO$, if $\mathbf{u}\geq \mathbf{v}$ for some $\mathbf{u}\in Q$, then $\mathbf{v}\in Q$.
        \item Given $\mathbf{v}\in\RO$, the sum of coordinates $\mathbf{u}_E$ of every maximal $\mathbf{u}\leq\mathbf{v}\in Q$ is the same.
    \end{enumerate}
    In particular, if $Q$ is the independence polytope of a polymatroid, then the rank function $f_Q$ of the polymatroid is
    \[
        f_Q(S) = \max\{\mathbf{x}_S\,\vert\,\mathbf{x}\in Q\},\quad\text{for all } S\subseteq E.
    \]
\end{proposition}

A polytope $R\subseteq\RO$ is the base polytope of a polymatroid if the set
\[
    Q=\{\mathbf{u}\in\RO\,\vert\,\text{ there exists }\mathbf{v}\in R\text{ such that }\mathbf{v}-\mathbf{u}\in\RO\}
\]
is the independence polytope of a polymatroid. When this is the case, the second item in Proposition~\ref{prop: characterization of polymatroids} guarantees $R$ is the face of $Q$ whose supporting hyperplane has normal vector $e_E=(1,\dots,1)$.

\begin{example}[Permutohedra and stellahedra as polymatroids] \label{bg: example perm and stella as polymatroid polytopes}
The \emph{permutohedron} in $\bb{R}^{n}$ is the polytope
\[
    \pe=\conv\{w\cdot(1,2,\dots,n)\,\vert\,w\text{ is a permutation of }[n]\}\subseteq\bb{R}^n.
\]
It is the base polytope of the polymatroid $(E,f)$ with $E=[n]$ and 
\[
    f(S)=\max\{\mathbf{x}_S\,\vert\,\mathbf{x}\in\pe\}=n + (n-1) + \cdots +(n-|S|+1).
\]  
 
The \emph{permutohedral variety} $\underline{X}_E$ is the smooth projective toric variety whose fan is the inner normal fan of the image of $\pe$ in $\bb{R}^E/\bb{R}e_E$.

On the other hand, the \emph{stellahedron} $\ste$ is the base polytope corresponding to this polymatroid. Alternatively, as a consequence of Proposition~\ref{prop: characterization of polymatroids}, we have that 
        \[
            \ste = \{\mathbf{u}\in\RO\,\vert\,\text{ there exists }\mathbf{v}\in\pe\text{ such that }\mathbf{v}-\mathbf{u}\in\RO\}.
        \]
 The inner normal fan of the stellahedron is called the \emph{stellahedral fan}. It is an unimodular simplicial fan with respect to the lattice $\ZZ^E\subseteq\RR^E$ \cite{eur2023stellahedral}. The \emph{stellahedral variety} $X_E$ is the projective toric variety whose fan is the stellahedral fan.
\end{example}

\begin{definition}[{\cite[Definition~2.6]{eur2024intersection}}]
    Let $\pi:A\to E$ be a caging with cage $\mathbf{a}$, and let $P=(E,f)$ be a polymatroid. The \emph{expansion} of $P$ with respect to $\pi$ is the polymatroid $\pi^*(P)$ on $A$ with rank function $f\circ\pi$. Equivalently, if $p_{\pi}:\RR^{A}\to\RR^E$ is the linear map induced by $\pi$, we have that
    \[
        I(\pi^*(P)) = p_{\pi}^{-1}(I(P))\cap\RR^{A}_{\geq 0}.
    \]
\end{definition}

We can now define the main polytopes that we study in this paper. 

\begin{definition}
\label{defn:Polypermutohedral}
    Let $\pi:A\to E$ be a caging with cage $\mathbf{a}$. Consider the polymatroid $P=(E,f)$ from Example~\ref{bg: example perm and stella as polymatroid polytopes} and its expansion $\pi^*(P)$. Then,
    \begin{itemize}
        \item The \emph{polypermutohedron with cage $\mathbf{a}$} is the base polytope of $\pi^*(P)$.
        \item The \emph{polystellahedron with cage $\mathbf{a}$} is the independence polytope of $\pi^*(P)$.
    \end{itemize}
    Both of these polytopes live in $\bb{R}^A$. However, we often consider the polypermutohedron as a polytope in $\bb{R}^A/\bb{R}e_A$. Let us now define the toric varieties corresponding to these polytopes.
    \begin{itemize}
         \item The \emph{polypermutohedral variety with cage $\mathbf{a}$} is the toric variety defined by the inner normal fan of the polypermutohedron with cage $\mathbf{a}$ in $\bb{R}^A/\bb{R}e_A$. We denote it by $P_\mathbf{a}$.
         
        \item The \emph{polystellahedral variety with cage $\mathbf{a}$} is the toric variety defined by the inner normal fan of the polystellahedron with cage $\mathbf{a}$ in $\bb{R}^A$. We denote it by $PS_\mathbf{a}$.
    \end{itemize}
    In particular, the polypermutohedron and polystellahedron with cage $\mathbf{a}=(1,\dots,1)$ coincide with the permutohedron and stellahedron, respectively.
\end{definition}

As already noted in \cite[Section 6]{eur2024intersection}, for any cage $\mathbf{a}$ the polypermutohedral variety $P_{\mathbf{a}}$ admits an embedding $P_{\mathbf{a}} \hookrightarrow PS_{\mathbf{a}}$ as a torus invariant divisor in the polystelahedral variety $PS_{\mathbf{a}}$. This is because the polypermutohedron with cage $\mathbf{a}$ is the facet of the polystellahedron with cage $\mathbf{a}$ corresponding to the ray spanned by $e_A=(1,\ldots,1) \in\bb{Z}^A\subseteq\mathbb{R}^{A}$.

\subsection{Wonderful compactifications}
Next, we introduced the necessary tools from Li's theory of wonderful compactifications \cite{li2009wonderful}. Throughout this subsection we fix $Y$ to be a nonsingular variety over $\bb{C}$.

\begin{definition}[{\cite[Section~5.1]{li2009wonderful}}]
    Let $A,A_1,\dots,A_k,B$ be smooth subvarieties of $Y$. Then,
    \begin{itemize}
        \item The intersection of $A$ and $B$ is said to be \emph{clean} if the set-theoretic intersection $A\cap B$ is smooth, and we have the following relation between tangent spaces:
        \[
            T_{A\cap B,y}=T_{A,y}\cap T_{B,y},\text{ for all }y\in A\cap B.
        \]
        \item The intersection of $A_1,\dots,A_k$ is said to be \emph{transversal}, if either $k=1$, or for all $y\in Y$
        \[
            \operatorname{codim}\left( \bigcap_{i=1}^k T_{A_i,y},T_y\right)
            =
            \sum_{i=1}^k \operatorname{codim}(A_i,Y).
        \]   
    \end{itemize}
\end{definition}

\begin{definition}[{\cite[Definition~2.1]{li2009wonderful}}]\label{bg: li arrangement def}
    An \emph{arrangement} of subvarieties of $Y$ is a finite collection $\cal{S}=\{S_i\}$ of properly contained smooth subvarieties $S_i\subsetneq Y$ such that for all $i\neq j$, $S_i$ and $S_j$ intersect cleanly, and $S_i\cap S_j$ is either equal to some $S_k$ or empty. In particular,  $\cal{S}$ is closed under intersections. 
\end{definition}

\begin{definition}[{\cite[Definition~2.2]{li2009wonderful}}]\label{bg: li def building set}
    Let $\cal{S}$ be an arrangement of subvarieties of $Y$. 
    \begin{itemize}
        \item     A subset $\cal{G}\subseteq\cal{S}$ is called a \emph{building set of $\cal{S}$} if, for all $S\in\cal{S}$, the minimal elements of $\{G\in\cal{G}\,\vert\, S\subseteq G\}$ intersect transversally, and their intersection is $S$. In particular,  this condition is automatically satisfied if $S\in\cal{G}$. 
        \item  The collection $\cal{G}$ is called a \emph{building set} if the set $\cal{S}$ of all possible intersections of elements of $\cal{G}$ is an arrangement, and furthermore $\cal{G}$ is a building set of $\cal{S}$.
    \end{itemize}
\end{definition}
We arrive at the main definition in our subsection 

\begin{definition}[{\cite[Definition~1.1]{li2009wonderful}\label{bg: def wonderful}}]
    Let $\cal{G}$ be a nonempty building set of subvarieties of $Y$, and $Y^\circ=Y\setminus \bigcup_{G\in\cal{G}}G$. Then, the closure of the image of the locally closed embedding
    \[
        Y^\circ\hookrightarrow\prod_{G\in\cal{G}}\operatorname{Bl}_G Y
    \]
    is called the \emph{wonderful compactification} of $\cal{G}$ on $Y$. It is denoted $Y_\cal{G}$.
\end{definition}    

\begin{definition}[{\cite[Definition~2.7]{li2009wonderful}}]
    Let $Z$ be a nonsingular subvariety of a nonsingular variety $Y$, and $\pi: \Bl_Z Y \to Y$ be the blow-up of $Y$ along $Z$. For any subvariety $V$ of $Y$, we define the \emph{dominant transform of $V$} to be the strict transform of $V$ if $V\nsubseteq Z$, and to be the scheme-theoretic inverse $\pi^{-1}(V)$ if $V\subseteq Z$.
\end{definition}

\begin{proposition}[{\cite[Theorem 1.2 and Proposition 2.13]{li2009wonderful}}] \label{prop: li blow-up building set}
    Let $Y$ be a nonsingular variety and $\mathcal{G} = \{S_1 , S_2 ,\dots , S_n\}$ be a building set of an arrangement of subvarieties of $Y$. Then,
    \begin{enumerate}
        \item[(1)] The wonderful compactification $Y_{\mathcal{G}}$ is a nonsingular variety. Moreover, for every $S_i\in\mathcal{G}$ there exists an irreducible smooth divisor $D_{S_i}$ such that:
        \begin{enumerate}
            \item[(i)] If $Y^\circ=Y\setminus\bigcup S_i$, then $Y_\mathcal{G}\setminus Y^\circ = \bigcup D_{S_i}$.
            \item[(ii)] Any set of divisors $D_{S_i}$ intersect transversely.
        \end{enumerate}
        \item[(2)] $Y_\mathcal{G}$ is isomorphic to an iterated blow-up $\Bl_\mathcal{G} Y$ of $Y$ along the dominant transforms of the elements in $\mathcal{G}$. 
    \end{enumerate}
\end{proposition}

As Li observes, the order in which the blow-ups are performed in part (2) of the previous Proposition is flexible. We have the following:

\begin{proposition}[{\cite[Theorem~1.3]{li2009wonderful}}] \label{bg: li order of blow-up doesnt matter}
    Let $Y$ be a nonsingular variety and $\mathcal{G} = \{S_1 , S_2 ,\dots , S_n\}$ be a building set of an arrangement of subvarieties of $Y$. Suppose the elements of $\mathcal{G}$ are ordered in such a way that the first $k$ terms $S_1,S_2,\dots, S_k$ form a building set for all $1\leq k\leq n$. Then, the wonderful compactification $Y_\mathcal{G}$ is isomorphic to an iterated blow-up $\Bl_\mathcal{G} Y$ of $Y$ along the dominant transforms of the elements of $\mathcal{G}$ in this chosen order. 
\end{proposition}

\subsection{Fulton-MacPherson compactifications and related constructions}
\label{sec:FMpreliminary}
Our compactifications will depend on a weight vector.
\begin{definition}\label{def:weights}
The domain of admissible weights for the weighted compactification of $n$ points in the flag $ 
\CC^{a_n} \subseteq \CC^{a_{n-1}} \subseteq \cdots \subseteq \CC^{a_1},
$  is the set
\[
    \cal{D}^{\operatorname{FM}}_n
    =
    \{
    (w_1, \ldots, w_n) \in \mathbb{Q}^n \; | \; 0 < w_i \leq 1, \forall i\in[n] \;
    \}
\]
We refer to the elements of these sets as \emph{weight vectors}.
\end{definition}

Let $X$ be a nonsingular variety. For all $I\subseteq[n]$ we define the $I$-diagonal 
\[
    \Delta_I=\{(x_1,\dots,x_n)\in X^n\,\vert\, x_i=x_j, \ \forall i,j\in I\}\subseteq X^n.
\]

Following the convention in \cite{fulton1994compactification}, throughout this article, we will refer to the flat families over our moduli spaces as their ``universal" families. These families are flat and parametrize the desired degenerations. The use of quotation symbols follows from the fact the Fulton-MacPherson space represents a functor of so-called screens \cite[Theorem~4]{fulton1994compactification}, but it is unknown if the Fulton-MacPherson spaces represent a functor of pointed degenerations of varieties.     
\begin{definition}[{\cite[Theorems~2~and~3]{routis2014weighted}}]
\label{def:weightedFM}
Let $X$ be a nonsingular variety of dimension $d$.
Given $\mathbf{w} \in \cal{D}^{\operatorname{FM}}_n$, define the building set $\mathcal{K}_{\mathbf{w}}$ on $X$ as
\[
\mathcal{K}_{\mathbf{w}} := 
\bigg\{
\Delta_I \subseteq X^n \; \big| \; I \subseteq [n],\; |I|\geq 2 \text{ and } a_I > 1
\bigg\}.
\]
Then, the weighted Fulton-MacPherson
compactification $X_{\mathbf{w}}[n]$ is the wonderful compactification of the open set
\[
X^n \setminus \bigcup_{\Delta_I \in \mathcal{K}_{\mathbf{w}} }  \Delta_I.
\]

Moreover, there is a smooth ``universal" family 
$\phi_{\mathbf{w}}: X_{\mathbf{w}}[n]^+ \to X_{\mathbf{w}}[n]$ 
equipped with $n$ sections 
$\sigma_i: X_{\mathbf{w}}[n] \to  X_{\mathbf{w}}[n]^+$ whose images lie in the relative smooth locus of $\phi_{\mathbf{w}}$. 
The map $\phi_{\mathbf{w}}$ is a flat morphism between nonsingular varieties whose fibers are the $n$-pointed $\mathbf{w}$-stable degenerations of $X$.
\end{definition}

In particular, if $\mathbf{w} = (1, \ldots, 1)$, the space $X_{\mathbf{w}}^{[n]}$ coincides with the classical Fulton-MacPherson compactification $X^{[n]}$ of \cite{fulton1994compactification}.

\section{Moduli of points in a flag and polypermutohedral varieties}
\label{sec:Polypermuta_moduli}

In this section we prove Theorem~\ref{GFM: thm:weightedFM_forFlags}. The proof is carried out in stages: First, in Subsection~\ref{subsection: GFM}, we construct a compactification of the configuration space of $n$ labeled weighted points in an affine flag. We then use this construction in Subsection~\ref{subsec: taw} to define the moduli space $T^{\mathbf{a}}_{\mathbf{w}}$. 
The toric compactification $T^{\mathbf{a}}_{LM}$ and its isomorphism with the polypermutohedral variety are described in Subsection~\ref{sec:toricTaw}. The nontrivial fibration of $T^{\mathbf{a}}_{\mathrm{LM}}$ is established in Subsection~\ref{sec:FibrationPP}. The proof of the theorem relies on the previous subsections, and is given in Subsection~\ref{sec:proofMainThmPolyPermetohedra}. The section concludes with Subsection~\ref{PP: subsec refinement}, where we prove the existence of a toric map from the standard Losev-Manin space to the compactification $T^{\mathbf{a}}_{LM}$.

Throughout this section we assume all cages $\mathbf{a}=(a_1,\dots,a_n)$ over $[n]$ are such that $a_1\geq\cdots\geq a_n$; see Definition~\ref{bg: def cages}. This constitutes no loss of generality and simplifies the presentation significantly.

\subsection{Weighted Fulton-MacPherson compactification of points in a flag
}\label{subsection: GFM}
Let \[
\CC^{a_n} \subseteq \CC^{a_{n-1}} \subseteq \cdots \subseteq \CC^{a_1}
\] be a flag of linear subspaces, where $\CC^{a_{i+1}}$ is the $a_{i+1}$-dimensional subspace of $\CC^{a_{i}}$ generated by the first $a_{i+1}$ coordinates.
We will constantly use the fact that given 
$\mathbf{a}=(a_1,\dots,a_n)$ and $\mathbf{d}=(a_1,\ldots,a_1)$,  we have a closed immersion 
\[
i_{\mathbf{a}}: \CC^{\mathbf{a}} = \prod_{i=1}^n \CC^{a_i} \hookrightarrow \CC^{\mathbf{d}} = (\CC^{a_1})^n
\]
given by the componentwise inclusion of each $\CC^{a_i}$ into $\CC^{a_1}$ as the subspace spanned by the first $a_i$ coordinates.

\begin{definition}
    Consider the cages $\mathbf{a}=(a_1,\dots,a_n)$ and $\mathbf{d}=(a_1,\dots,a_1)$. Suppose $I\subseteq [n]$ is such that $|I|\geq 2$. The \emph{$I$-diagonal of the $\mathbf{a}$-flag} is
    \[
        \Delta_I^{\mathbf{a}} = \{(x_1,\dots,x_n)\in (\CC^{a_1})^n \;|\; x_i\in\CC^{a_i}\text{ and }x_i=x_j\text{ for all }i,j\in I\}.
    \]
    When $a_1=\cdots=a_n$, the $I$-diagonal $\Delta_I^\mathbf{a}$ is the usual $I$-diagonal $\Delta_I$ in $(\CC^{a_1})^n$.
\end{definition}

\begin{lemma}\label{GFM: building set}
    Let $\mathbf{w}=(w_1,\dots,w_n)$ be a weight vector in $\mathcal{D}^{\FM}_n$; see Definition~\ref{def:weights}. The set 
    \[
        \mathcal{K}_\mathbf{w}^{\mathbf{a}} := \left\{\Delta_I^{\mathbf{a}} \,\vert\,  \mathbf{w}_I> 1\right\}.
    \]
    is a building set in $\CC^{\mathbf{a}}$ in the sense of Li \cite{li2009wonderful}; see Definition~\ref{bg: li def building set}. 
\end{lemma}

We point out that the case $a_1 = \cdots =a_n$ is a particular case of Routis' results in \cite{routis2014weighted}.

\begin{proof}[Proof of Lemma ~\ref{GFM: building set}]
    Let $\langle\mathcal{K}_\mathbf{w}^{\mathbf{a}}\rangle$ be the set of all possible scheme-theoretic intersections of the elements $\mathcal{K}_\mathbf{w}^{\mathbf{a}}$. Since $\Delta_I^{\mathbf{a}}\cap\Delta^{\mathbf{a}}_J=\Delta^{\mathbf{a}}_{I\cup J}$ if and only if $I\cap J\neq\emptyset$, we have the following decomposition:
    \[
        \langle\mathcal{K}_\mathbf{w}^{\mathbf{a}}\rangle=
        \left\{
        \Delta_{I_1}^\mathbf{a}\cap\cdots\cap\Delta_{I_r}^\mathbf{a}\;\big|\; r\geq 2,\, \Delta_{I_j}^\mathbf{a}\in\cal{K}_\mathbf{w}^\mathbf{a},\text{ and }I_j\cap I_k=\emptyset\text{ for all }j\neq k
        \right\}
        \cup
        \mathcal{K}_\mathbf{w}^{\mathbf{a}}. 
    \]
    Notice that all elements in $\langle\mathcal{K}_\mathbf{w}^{\mathbf{a}}\rangle$ are smooth varieties. 
    Moreover, the set $\langle\mathcal{K}_\mathbf{w}^{\mathbf{a}}\rangle$ is an arrangement of subvarieties in the sense of Li (see Definition~\ref{bg: li arrangement def}), this is, any two elements intersect cleanly in $\CC^{\mathbf{a}}$. 
    To see this, recall that two nonsingular subvarieties $A,B\subseteq Y$ of a variety $Y$ are said to intersect cleanly if their intersection is smooth and
    their tangent spaces satisfy $T_{A,y}\cap T_{B,y}=T_{A\cap B,y}$ for all $y\in A\cap B$. In the case at hand every element of $\langle\mathcal{K}_\mathbf{w}^{\mathbf{a}}\rangle$ is a linear subspace of $\CC^{\mathbf{a}}$, in which case this condition is immediate, since each one of the subvarieties can be identified with its tangent space.

    It only remains to prove that $\mathcal{K}_\mathbf{w}^{\mathbf{a}}$ is a building set for the arrangement $\langle\mathcal{K}_\mathbf{w}^{\mathbf{a}}\rangle$. By definition, this amounts to proving that given any $S\in\langle\mathcal{K}_\mathbf{w}^{\mathbf{a}}\rangle\setminus\mathcal{K}_\mathbf{w}^{\mathbf{a}}$, the minimal elements of $\{\Delta_I^{\mathbf{a}}\in\mathcal{K}_\mathbf{w}^{\mathbf{a}}:\Delta_I^{\mathbf{a}}\supseteq S\}$ intersect transversely and their intersection is $S$. This is a direct dimension count. 

    Let us once again identify each one of the diagonals with their tangent spaces. By 
    \cite[Section 5.1.2]{li2009wonderful}
     a collection $\Delta_{I_1}^{\mathbf{a}},\dots,\Delta_{I_r}^{\mathbf{a}}$ of diagonals with $\{I_1,\ldots,I_r\}$ pairwise disjoint 
     intersects transversely if and only if
    \begin{equation}\label{GFM: eqn codim}
        \sum_{j=1}^r\codim\left(\Delta_{I_j}^{\mathbf{a}},\CC^{\mathbf{a}}\right) 
        = \codim\left(\bigcap_{j=1}^r \Delta_{I_j}^{\mathbf{a}},\CC^{\mathbf{a}}\right).
    \end{equation}
    Define $i_j$ to be the largest element of $I_j$. Then, a direct dimension count shows that 
    \[
        \codim(\Delta_{I_j}^{\mathbf{a}},\CC^{\mathbf{a}}) = \sum_{k\in I_j}a_k - a_{i_j}, 
    \]  
    so that the left hand side of Equation~\eqref{GFM: eqn codim} becomes $\sum_{k\in\cup I_j} a_k - a_{i_1} - \cdots - a_{i_r}$. 
    On the other hand, computing the right-hand side by directly counting the codimension of $\Delta_{I_1}^{\mathbf{a}}\cap\cdots\cap\Delta_{I_r}^{\mathbf{a}}$ in $\mathbb{C}^{\mathbf{a}}$ shows that it is also equal to 
    $\sum_{k\in\cup I_j} a_k - a_{i_1} - \cdots - a_{i_r}$.    
\end{proof}

\begin{definition} \label{GFM: definition.weighted.FM.for.flags}
    The \emph{weighted compactification} 
    $\mathbb{C}_{\mathbf{w}}^{[\mathbf{a}]}$ of $\CC^{\mathbf{a}}\setminus\bigcup_{\Delta_I^{\mathbf{a}}\in\mathcal{K}_\mathbf{w}^{\mathbf{a}}}\Delta_I^{\mathbf{a}}$ is the wonderful compactification of the building set $\mathcal{K}_\mathbf{w}^{\mathbf{a}}$ in $\mathbb{C}^{\mathbf{a}}$.
    We refer to $\mathbb{C}_{\mathbf{w}}^{[\mathbf{a}]}$
    as the \emph{weighted Fulton-MacPherson compactification of points in a flag.}  
\end{definition}

\begin{example}\label{GFM: remark evangelos weighted fulton macpherson}
    If $\mathbf{d}=(d,\dots,d)$ is a constant cage, then the compactifications $\bb{C}^{[\mathbf{d}]}_\mathbf{w}$ specialize to Routis' weighted Fulton-MacPherson for configurations of points in $\mathbb{C}^d$; see \cite[Definition~2.6]{routis2014weighted}. 
\end{example}

The wonderful compactification $\mathbb{C}^{[\mathbf{a}]}_{\mathbf{w}}$ is endowed with a structural map $\psi_{\mathbf{a}}: 
\mathbb{C}^{[\mathbf{a}]}_{\mathbf{w}} \rightarrow \mathbb{C}^{\mathbf{a}}$, which will be important for the following results.
For each $I_0 \subseteq [n]$ such that $|I_0|\geq 2$ and $\Delta^{\mathbf{a}}_{I_0} \in \mathcal{K}_\mathbf{w}^{\mathbf{a}}$, 
the morphism $\psi_{\mathbf{a}}: \mathbb{C}^{[\mathbf{a}]}_{\mathbf{w}} \rightarrow \mathbb{C}^{\mathbf{a}}$ is given by the composition 
\begin{align}\label{eq:structuralMorphism}
\mathbb{C}^{[\mathbf{a}]}_{\mathbf{w}} 
\rightarrow
X_{\mathbf{a}} 
:= \prod_{\Delta_I^{\mathbf{a}} \in \mathcal{K}_\mathbf{w}^{\mathbf{a}}}\Bl_{\Delta_I^{\mathbf{a}}} \CC^{\mathbf{a}}
\rightarrow
\Bl_{\Delta_{I_0}^{\mathbf{a}}} \CC^{\mathbf{a}}
\rightarrow
\CC^{\mathbf{a}},
\end{align}
where the map on the left is the inclusion $\mathbb{C}^{[\mathbf{a}]}_{\mathbf{w}} \subseteq X_{\mathbf{a}}$ in the definition of $\mathbb{C}^{[\mathbf{a}]}_{\mathbf{w}}$, the middle map is the projection on the indicated factor and the map on the right is the indicated blow-up. 
The morphism $\psi_{\mathbf{a}}:\mathbb{C}^{[\mathbf{a}]}_{\mathbf{w}} \rightarrow \mathbb{C}^{\mathbf{a}}$ is independent of the choice of $I_0$, since when restricted to the open dense subset $\mathcal{U}^{\mathbf{a}} \subseteq \mathbb{C}^{[\mathbf{a}]}_{\mathbf{w}}$ it agrees with the inclusion $\mathcal{U}^{\mathbf{a}} \hookrightarrow \mathbb{C}^{\mathbf{a}}$. 
In particular, $\psi_{\mathbf{a}}$ is birational and it is the identity when resticted to $\mathcal{U}^{\mathbf{a}} \subseteq \mathbb{C}^{[\mathbf{a}]}_{\mathbf{w}}$ on the source and $\mathcal{U}^{\mathbf{a}} \subseteq \mathbb{C}^{\mathbf{a}}$ on the target. 
Moreover, \cite[Theorem 1.3]{li2009wonderful} shows that $\psi_{\mathbf{a}}: \mathbb{C}^{[\mathbf{a}]}_{\mathbf{w}} \rightarrow \mathbb{C}^{\mathbf{a}}$ is an iterated blow-up along the dominant transforms of the varieties $\Delta^{\mathbf{a}}_{I} \in \mathcal{K}_\mathbf{w}^{\mathbf{a}}$ in any order such that any initial sub-list of the ordered list is a building set (see Proposition~\ref{bg: li order of blow-up doesnt matter}).

%%%%%%%%%%%%%%%%%%%%%%%%%%%%%%%%%%%%%%%%%%%%%%%%%%%%%%%%%%%%%%%%%

The following lemma shows how a smooth subvariety interacts with the wonderful compactification of a larger ambient space, assuming that its intersections with the ambient building set form a building set.

\begin{lemma}\label{lemma.wonderful.subvariety}
Let $X$ be a smooth variety, 
$Y$ be a smooth subvariety of $X$
and $\mathcal{G}$  be a building set on $X$. 
Suppose that the collection of the scheme-theoretic intersections with $Y$ of the varieties in $\mathcal{G}$ forms a building set $\mathcal{H}$ in $Y$.  
Then, there exists a unique closed immersion $j: Y_{\mathcal{H}} \hookrightarrow X_{\mathcal{G}}$ 
which makes the following diagram commute:  
\[
\begin{tikzcd}
    Y_{\mathcal{H}} 
    \ar[r,hook, "j"] \ar[d, "g"'] & 
    X_{\mathcal{G}} 
    \ar[d, "f"] 
    \\
    Y \ar[r,hook, "i"] & 
    X,
\end{tikzcd}
\]
where $f: X_{\mathcal{G}} \rightarrow X$ and $g: Y_{\mathcal{H}} \rightarrow Y$ are the wonderful compactification morphisms  
and $i: Y \hookrightarrow X$ is the inclusion. 
In particular, $Y_{\mathcal{H}}$ is the strict transform of $Y$ under the morphism $f:X_{\mathcal{G}} \rightarrow X$ considered as an iterated blow-up. 
\end{lemma}

\begin{proof}
Fix any subvariety $Z$ in $\mathcal{G}$. 
By \cite[Corollary 7.15]{hartshorne-book}, the closed immersion $i:Y \hookrightarrow X$ induces a closed immersion $i_Z: \Bl_{Z \cap Y} Y
        \hookrightarrow
          \Bl_{Z} X$ such that the following diagram commutes
    \[
    \begin{tikzcd}
        \Bl_{Z \cap Y} Y
        \ar[r,hook, "i_Z"]\ar[d] & 
          \Bl_{Z} X 
        \ar[d]
        \\
        Y\ar[r,hook, "i"] & X.
    \end{tikzcd}
    \]
Then, if we let $Z$ vary and we take the product over all $Z \in \mathcal{G}$, then we get a closed immersion $\widetilde{Y} \hookrightarrow \widetilde{X}$ where
\[
    \widetilde{Y} := \prod_{Z \in \mathcal{G}} \Bl_{Z \cap Y} Y
    \longrightarrow
    \widetilde{X} := \prod_{Z \in \mathcal{G}} \Bl_{Z} X.  
\]
Let $U_X = X \setminus \bigcup_{Z \in \mathcal{G}} Z$
and 
$U_Y = Y \setminus \bigcup_{Z \in \mathcal{G}} (Z \cap Y)$.
From the construction of the wonderful compactifications, the morphisms $f$ and $g$ are isomorphisms over $U_X$ and $U_Y$. 
We notice that $U_Y = Y \cap U_X$ and hence also $U_Y = \widetilde{Y} \cap U_X$. 
Therefore, we get a commutative diagram
\[
\begin{tikzcd}
    \widetilde{Y}
    \ar[r,hook] & 
    \widetilde{X}
    \\
    U_Y\ar[r,hook]\ar[u,hook] & U_X\ar[u,hook]
\end{tikzcd}
\]
where the vertical maps are open immersions and the horizontal maps are closed immersions. 
This diagram induces an identification between the closure of $U_Y$ in $\widetilde{Y}$, that is $Y_{\mathcal{H}}$, and the closure of $U_Y$ in $\widetilde{X}$. 
The diagram also induces a closed immersion from the closure of $U_Y$ in $\widetilde{X}$ into the closure of $U_X$ in $\widetilde{X}$, that is $X_{\mathcal{G}}$. 
Then we obtain a closed immersion $j:Y_{\mathcal{H}} \hookrightarrow X_{\mathcal{G}}$. 

If we consider the diagram in the statement using the closed immersion $j$ that we constructed as the top morphism, the diagram commutes because by construction it commutes over the dense open subset $U_Y$ of $Y_{\mathcal{H}}$.

Now, notice that there can be at most one morphism $Y_{\mathcal{H}} \hookrightarrow X_{\mathcal{G}}$ making the diagram in the statement commute because such a morphism has prescribed values over the dense open subset $U_Y$ of $Y_{\mathcal{G}}$.  
Indeed, by the commutativity of the diagram, when the morphism $Y_{\mathcal{H}} \hookrightarrow X_{\mathcal{G}}$ is restricted to $U_Y$ it agrees with the composition
\[
(f|_{U_X})^{-1}
\circ i \circ 
(g|_{U_Y}). 
\]
Then, the existence and uniqueness of the closed immersion $j$ making the diagram commute follows. 
Finally, the image of $Y_{\mathcal{H}}$ under $j$ is a subvariety of $X_{\mathcal{G}}$ with the same dimension as $Y$ and which maps birationally onto $Y$ under $f$, but this characterizes the strict transform of $Y$ under $f$, and this completes the proof.  
\end{proof}

In the case of the weighted compactification $\mathbb{C}_{\mathbf{w}}^{[\mathbf{a}]}$, we arrive at the following result.

\begin{lemma}\label{GFM:lemma.closed.immersion} Given positive integers $a_1 \geq \cdots \geq a_n$, let $\mathbf{a}=(a_1,\dots,a_n)$ and $\mathbf{d}=(a_1,\dots,a_1)$, and let $\mathbf{w} \in \mathcal{D}^{\FM}_n$ be a set of admissible weights. Then, there exists a unique closed immersion $j_{\mathbf{a}}:\mathbb{C}^{[\mathbf{a}]}_{\mathbf{w}} \hookrightarrow \mathbb{C}^{[\mathbf{d}]}_{\mathbf{w}}$ which makes the following diagram commute:
\[
\begin{tikzcd} \mathbb{C}^{[\mathbf{a}]}_{\mathbf{w}} \ar[r,hook,"j_{\mathbf{a}}"]\ar[d,"\psi_{\mathbf{a}}"'] & \mathbb{C}^{[\mathbf{d}]}_{\mathbf{w}} \ar[d,"\psi_{\mathbf{d}}"] \\ \CC^{\mathbf{a}}\ar[r,hook,"i_{\mathbf{a}}"] & \CC^{\mathbf{d}} \end{tikzcd} 
\]
where $\psi_{\mathbf{a}}$ and $\psi_{\mathbf{d}}$ are the wonderful compactification morphisms (see Equation~\eqref{eq:structuralMorphism}), and $i_\mathbf{a}$ is the inclusion. 
\end{lemma}
\begin{proof} This result follows as a corollary of Lemma \ref{lemma.wonderful.subvariety}. Indeed, let $X = \CC^{\mathbf{d}}$, $Y = \CC^{\mathbf{a}}$, and let $\mathcal{G}=\mathcal{K}_\mathbf{w}^{\mathbf{d}}$  be the set of all diagonals $\Delta_I^{\mathbf{d}}$ associated with $\mathbf{w}$. 
Consider any $I_0 \subseteq [n]$ such that $|I_0|\geq 2$. 
Notice that the $I_0$-diagonal $\Delta_{I_0}^{\bf{a}}$ is in $\mathcal{K}_\mathbf{w}^{\mathbf{a}}$ if and only if 
the $I_0$-diagonal $\Delta_{I_0}^{\bf{d}}$ is in $\mathcal{K}_\mathbf{w}^{\mathbf{d}}$.  
Also, the $I_0$-diagonal $\Delta_{I_0} = \Delta_{I_0}^{\bf{d}} \subseteq \CC^{\mathbf{d}}$ intersects $\CC^{\mathbf{a}}$, since even the deepest diagonal $\Delta_{[n]} \subseteq \Delta_{I_0}$ clearly intersects $\CC^{\mathbf{a}}$.  
Moreover, $\Delta_{I_0}^{\bf{a}}  =  \CC^{\mathbf{a}}\cap \Delta_{_0}^{\bf{d}}$ scheme-theoretically.
Then, the intersections of the diagonals in $\mathcal{G}=\mathcal{K}_\mathbf{w}^{\mathbf{d}}$ with $Y = \CC^{\mathbf{a}}$ form the corresponding building set $\mathcal{H}=\mathcal{K}_\mathbf{w}^{\mathbf{a}}$ of all diagonals $\Delta_I^{\mathbf{a}}$ in $\CC^{\mathbf{a}}$ associated with $\mathbf{w}$. 
By applying Lemma \ref{lemma.wonderful.subvariety} in this setting, we obtain the existence and uniqueness of the desired closed immersion $j_\mathbf{a}: \mathbb{C}^{[\mathbf{a}]}_{\mathbf{w}} \hookrightarrow \mathbb{C}^{[\mathbf{d}]}_{\mathbf{w}}$, making the diagram commute. \end{proof}

%%%%%%%%%%%%%%%%%%%%%%%%%%%%%%%%%%%%%%%%%%%%%%%%%%%%%%%%%%%%%%%%%

To conclude this series of lemmas let us construct the ``universal" family of the configuration spaces $\bb{C}^{[\mathbf{a}]}_\mathbf{w}$. Routis constructs the ``universal" family of the weighted Fulton-MacPherson compactification when the cage is constant; see Example~\ref{GFM: remark evangelos weighted fulton macpherson}. In light of the previous Lemma, this gives a ``universal" family for our compactifications $\bb{C}^{[\mathbf{a}]}_\mathbf{w}$ for all cages $\mathbf{a}$.

\begin{lemma}[{\cite[Theorem~2.3]{routis2014weighted}}]
    Let $\mathbf{d}$ be a constant cage. There exists a ``universal" family $(\bb{C}^{[\mathbf{d}]}_\mathbf{w})^+\to\bb{C}^{[\mathbf{d}]}_\mathbf{w}$ equipped with $n$ sections $\sigma_i:\bb{C}^{[\mathbf{d}]}_\mathbf{w}\to(\bb{C}^{[\mathbf{d}]}_\mathbf{w})^+$ whose images are disjoint and give the marked points. The ``universal" family is a flat morphism between nonsingular varieties, whose fibers are $\mathbf{w}$-stable degenerations of $\bb{C}^\mathbf{d}$.
\end{lemma}

\begin{definition}\label{GFM: definition universal family}
    Consider the cages $\mathbf{a}=(a_1,\dots,a_n)$ and $\mathbf{d}=(a_1,\dots,a_1)$, and let $\mathbf{w}\in\cal{D}^{\FM}_n$. Define $(\bb{C}^{[\mathbf{a}]}_\mathbf{w})^+$ the ``universal" family of the configuration spaces $\bb{C}^{[\mathbf{a}]}_\mathbf{w}$ to be the pullback in the following Cartesian diagram:
    \begin{equation}\label{sec moduli: diagram definition family}
    \begin{tikzcd}
        (\bb{C}_\mathbf{w}^{[\mathbf{a}]})^+\ar[r,hook]\ar[d] & (\mathbb{C}_\mathbf{w}^{[\mathbf{d}]})^+\ar[d] \\
        \bb{C}^{[\mathbf{a}]}_\mathbf{w}\ar[r,hook] & \bb{C}_\mathbf{w}^{[\mathbf{d}]},
    \end{tikzcd}
    \end{equation}
    where the bottom horizontal arrow is the closed immersion from Lemma~\ref{GFM:lemma.closed.immersion}. The map $(\bb{C}_\mathbf{w}^{[\mathbf{a}]})^+\to\bb{C}_\mathbf{w}^{[\mathbf{a}]}$ is equipped with $n$ sections  $\sigma_i:\bb{C}_\mathbf{w}^{[\mathbf{a}]}\to(\bb{C}_\mathbf{w}^{[\mathbf{a}]})^+$, $i\in[n]$, whose images are disjoint, obtained from the universal property of the Cartesian product via the identity map of $\bb{C}_\mathbf{w}^{[\mathbf{a}]}$ and the composition 
    \[
    \begin{tikzcd}
        \bb{C}_\mathbf{w}^{[\mathbf{a}]}\ar[r,hook]& \bb{C}_\mathbf{w}^{[\mathbf{d}]}\ar[r,"\sigma_i"]&(\bb{C}_\mathbf{w}^{[\mathbf{d}]})^+.
    \end{tikzcd}
    \]
\end{definition}

\subsection{Moduli of $n$ weighted points in a flag up to translation and scaling
}\label{subsec: taw}
The goal of this section is to give a direct construction of the moduli space $T_\mathbf{w}^\mathbf{a}$ of points in the flag, up to scaling and translation by the $\mathbb{C}^{a_n}$-coordinates, as a wonderful compactification, and to construct its ``universal" family.

\begin{definition}\label{def:taw}
    Consider a cage $\mathbf{a}$ and a collection of weight data $\mathbf{w}\in\cal{D}^{\FM}_n$ such that $\mathbf{w}_{[n]}>1$. Let $D_{[n]}^\mathbf{a}=\pi^{-1}(\Delta_{[n]}^\mathbf{a})\subseteq \bb{C}_\mathbf{w}^{\mathbf{a}}$ be the exceptional divisor of $\pi:\mathbb{C}_\mathbf{w}^{[\mathbf{a}]}\to\CC^\mathbf{a}$ above the deepest diagonal $\Delta_{[n]}^\mathbf{a}$. 
    We define the variety 
    \[
        T_\mathbf{w}^\mathbf{a} := (\pi\vert_{D_{[n]}}^\mathbf{a})^{-1}(0)
    \]
    as the fiber above $0\in\Delta_{[n]}^\mathbf{a}\subseteq\bb{C}^\mathbf{a}$. In particular, every square in the following diagram is Cartesian:
    \begin{equation}\label{sec moduli: diagram definition T}
    \begin{tikzcd}
        T_\mathbf{w}^\mathbf{a}\ar[r,hook]\ar[d] & D_{[n]}^\mathbf{a}\ar[r,hook]\ar[d,"\pi\vert_{D_{[n]}}"] & \mathbb{C}_\mathbf{w}^{[\mathbf{a}]}\ar[d,"\pi"] \\
        0\ar[r,hook] & \Delta_{[n]}^\mathbf{a}\ar[r, hook] & \CC^\mathbf{a}.
    \end{tikzcd}
    \end{equation}
\end{definition}
The above construction is similar to the one for the case of $n$ points in affine space; see \cite[Definition~3.1.1]{chen2009pointed} and \cite[Definition~4.10]{gallardo2017wonderful}. We highlight that it's independent of the point $x\in\Delta_{[n]}^\mathbf{a}$ used. However, we make this specific choice for convenience. Indeed, every element in the building set $\mathcal{K}_\mathbf{w}^\mathbf{a}$ is invariant under the action of $\Delta_{[n]}^\mathbf{a}$ on $\bb{C}^\mathbf{a}$ by translation. In particular, for all $x\in\Delta_{[n]}^\mathbf{a}$ this action defines an isomorphism $T_\mathbf{w}^\mathbf{a}\cong(\pi|_{D_{[n]}^\mathbf{a}})^{-1}(x)$. 

The requirement $\mathbf{w}_{[n]}>1$ guarantees that the building set of $\bb{C}^{[\mathbf{a}]}_\mathbf{w}$ is not empty. This motivates the following.

\begin{definition}\label{pp: twa weights}
    The domain of admissible weights for the moduli space of $n$ points in the flag $\bb{C}^\mathbf{a}$ is the set
    \[
        \cal{D}^{\operatorname{T}}_n
        =
        \{
        \mathbf{w}=(w_1, \ldots, w_n) \in \mathbb{Q}^n \; | \; 0 < w_i \leq 1,\, \forall i\in[n],\, \mathbf{w}_{[n]}>1
        \}.
    \]
\end{definition}

Consider cages $\mathbf{a}=(a_1,\dots,a_n)$ and 
$\mathbf{d}=(a_1,\dots,a_1)$ over $[n]$, and let 
$\mathbf{w}\in\mathcal{D}^{\operatorname{T}}_n$. For ease 
of notation let us denote the deepest diagonal as 
$\Delta:=\Delta_{[n]}^{\mathbf{a}}$. Suppose 
$\mathcal{K}_\mathbf{w}^{\mathbf{a}}$ (Definition \ref{GFM: building set}) is endowed with a total order satisfying the condition of Proposition~\ref{bg: li order of blow-up doesnt matter}, and such 
that $\Delta$ is its minimal element (e.g. an order with nondecreasing dimensions). Under this order, the space 
$\bb{C}_\mathbf{w}^{[\mathbf{a}]}$ is an iterated blow-up 
$\pi:\bb{C}_\mathbf{w}^{[\mathbf{a}]}\to\CC^{\mathbf{a}}$, its first step being the blow-up 
$\pi_1:\Bl_{\Delta}\CC^{\mathbf{a}}\to\CC^{\mathbf{a}}$. 
The fiber of $0\in\CC^{\mathbf{a}}$ under $\pi_1$ is $\bb{P}((N_{\Delta\vert\CC^{\mathbf{a}}})^\vee)\vert_0$, where $N_{\Delta\vert\CC^{\mathbf{a}}}$ is the normal bundle of $\Delta\subseteq\bb{C}^\mathbf{a}$. The conormal bundle $(N_{\Delta\vert\CC^{\mathbf{a}}})^\vee$ has constant fiber $(\bb{C}^\mathbf{a}/\Delta)^\vee$ over every point of $\Delta \cong \mathbb{C}^{a_n}$, so in particular this central fiber is isomorphic to the projectivization of this vector space. Therefore, it is possible to coordinatize this fiber via an isomorphism
\begin{equation} \label{equation.fiber.coordinates}
    \bb{P}((\bb{C}^\mathbf{a}/\Delta)^\vee)\cong\bb{P}^{a_1+\cdots + a_{(n-1)}-1}.
\end{equation}
The points in $\bb{P}((\bb{C}^\mathbf{a}/\Delta)^\vee)$ are in correspondence with the nonzero vectors in $\bb{C}^\mathbf{a}/\Delta$, up to scaling. 
In turn, the nonzero vectors in $\bb{C}^\mathbf{a}/\Delta$, up to scaling, are in correspondence with equivalence classes of configurations $(q_1,\dots,q_n)\in\mathbb{C}^\mathbf{a} \setminus \Delta$, up to scaling and translation by elements in $\Delta$.  
Each such equivalence class has a unique representative of the form $(p_1,\dots,p_{n-1},0)\in\bb{C}^\mathbf{a} \setminus \{ 0 \}$, up to scaling. 
Then, we fix the identification in (\ref{equation.fiber.coordinates}) as the unique isomorphism that maps the point of $\bb{P}((\bb{C}^\mathbf{a}/\Delta)^\vee)$ associated to $(p_1,\dots,p_{n-1},0)$ to the point
\[
    [x_{11}:\dots:x_{1a_1}: x_{21}:\dots:x_{2a_2}:\dots: x_{(n-1)1}:\dots:x_{(n-1)a_{n-1}}]\in\bb{P}^{a_1+\cdots+a_{n-1} - 1},
\]
where $p_i=(x_{i1},\dots,x_{ia_i})$, for each $i\in [n-1]$.

For all $I\subsetneq [n]$, the image $\delta_I^\mathbf{a}$ of the $I$-diagonal $\Delta_I^\mathbf{a}$ under this isomorphism is
\begin{align}\label{PP: delta I definition}
    \delta_I^\mathbf{a} = \left\{
    \begin{aligned}
    &\bigcap_{i\in I\setminus\{n\}} V(x_{i1},\dots,x_{ia_i}),\text{ if }n\in I;\\
    &\;\,\bigcap_{i,j \in I, i \leq j} V(x_{i1}-x_{j1},\dots,x_{ia_j}-x_{ja_j}, x_{i(a_j+1)},\dots, x_{ia_i}),\text{ if }n\not\in I.
    \end{aligned}
    \right.
\end{align}

\begin{lemma}\label{lemma.intersections.linear.spaces.with.exceptional.divisor}
Let $\mathbb{C}^m$ be an affine space and let $V \subsetneq W$ be linear subspaces of $\mathbb{C}^m$. 
Consider the blow-up $\operatorname{Bl}_V \mathbb{C}^m \rightarrow \mathbb{C}^m$, with exceptional divisor $E$ and let $\widetilde{W}$ be the strict transform of $W$. 
Then any fiber of $E \rightarrow V$ can be identified with $\mathbb{P}((\mathbb{C}^m/V)^{\vee})$, and under this identification the intersection of  $\widetilde{W}$ with that fiber of $E \rightarrow V$ is the image of $\mathbb{P}((W/V)^{\vee})$ under the closed immersion $\mathbb{P}((W/V)^{\vee}) \rightarrow \mathbb{P}((\mathbb{C}^m/V)^{\vee})$ induced by the inclusion $V \subsetneq W$. 

\begin{proof}
The exceptional divisor $E$ is the projective bundle $E=\mathbb{P}((N_{V|\mathbb{C}^m})^{\vee})$ over $V$. Identifying the tangent space of a linear subspace with the linear subspace itself, we see that the bundle $(N_{V|\mathbb{C}^m})^{\vee}$ has constant fiber $(\mathbb{C}^m/V)^{\vee}$.

To describe the intersection of the strict transform $\widetilde{W}$ with each fiber of $E \rightarrow V$, we start by identifying this intersection set-theoretically. 
Such an intersection is a variety (reduced and irreducible) and hence is determined by its underlying set. 
Fix any element of $V$, which we can assume to be $0$ by homogeneity. We then have a Cartesian diagram representing the fiber of $E \rightarrow V$ over $0 \in V$ as follows:
\[
\begin{tikzcd}
     \mathbb{P}((\mathbb{C}^m/V)^{\vee}) \arrow[r, ""] \arrow[d, ""] & E=\mathbb{P}((N_{V|\mathbb{C}^m})^{\vee}) \arrow[d, ""] \\
    0 \arrow[r, ""] & V.
\end{tikzcd}
\]

The intersection of $\widetilde{W}$ with this fiber is obtained set-theoretically by projectivizing the image of the injective linear map $W/V \rightarrow \mathbb{C}^m/V$. This injective linear map induces a surjective linear map $(\mathbb{C}^m/V)^{\vee} \rightarrow (W/V)^{\vee}$, and therefore a closed immersion $\mathbb{P}((W/V)^{\vee}) \rightarrow \mathbb{P}((\mathbb{C}^m/V)^{\vee})$.
It follows that the intersection of the strict transform $\widetilde{W}$ with any fiber of $E \rightarrow V$ is the image of $\mathbb{P}((W/V)^{\vee})$ under this closed immersion, as desired.
\end{proof}
\end{lemma}

%%%%%%%%%%%%%%%%%%%%%%%%%%%%%%%%%%%%%%%%%%%%%%%%%%%%%%%%%%%%%%%%%%

\begin{proposition}\label{proposition:taw_is_wonderful}
    Consider a cage $\mathbf{a}=(a_1,\dots,a_n)$ and $\mathbf{w}\in\cal{D}^{\operatorname{T}}_n$. Then, $T_\mathbf{w}^{\mathbf{a}}$ is the wonderful compactification of $\bb{P}^{a_1+\cdots+a_{(n-1)}-1}$ with respect to the building set 
    \[
        \cal{H}_\mathbf{w}^\mathbf{a}=\{\delta_I^\mathbf{a}\,\vert\, I\subsetneq [n],\, \mathbf{w}_I>1\}.
    \]
\end{proposition} 

\begin{proof} 
From Definition~\ref{def:taw} we get the following commutative diagram consisting of Cartesian squares: 
\begin{equation} \label{commutative.diagram.wonderful.Twa}
\begin{tikzcd}
    T_\mathbf{w}^\mathbf{a} \arrow[r, ""] \arrow[d, ""] & \mathbb{C}_\mathbf{w}^{[\mathbf{a}]} \arrow[d, ""] \\
    \bb{P}((\bb{C}^\mathbf{a}/\Delta)^\vee) \cong \mathbb{P}^N \arrow[r, ""] \arrow[d, ""] & \operatorname{Bl}_{\Delta} \CC^\mathbf{a} \arrow[d, ""] \\
    0 \arrow[r, ""] & \CC^\mathbf{a},
\end{tikzcd}
\end{equation}
where $N:= a_1+\cdots + a_{(n-1)}-1$. Let us recall that $\mathbb{C}_\mathbf{w}^{[\mathbf{a}]}$
is the wonderful compactification of the building set $\mathcal{G}:=\mathcal{K}_\mathbf{w}^{\mathbf{a}}$ in $\CC^\mathbf{a}$.
Since $\Delta$ is minimal in $\mathcal{G}$ by \cite[Proposition 2.8]{li2009wonderful} the dominant transforms to $\operatorname{Bl}_{\Delta} \CC^\mathbf{a}$ of the elements in $\mathcal{G}$ form a building set $\widetilde{\mathcal{G}}$ in $\operatorname{Bl}_{\Delta} \CC^\mathbf{a}$. 
From the inductive construction of wonderful compactifications in \cite[Definition 2.12]{li2009wonderful}, we know that $\mathbb{C}_\mathbf{w}^{[\mathbf{a}]}$ is the wonderful compactification 
of the building set $\widetilde{\mathcal{G}}$ in $\operatorname{Bl}_{\Delta} \CC^\mathbf{a}$. 

Moreover, the morphism $\mathbb{C}_\mathbf{w}^{[\mathbf{a}]} \rightarrow \operatorname{Bl}_{\Delta} \CC^\mathbf{a}$ 
is the iterated blow-up of the dominant transforms of the elements in $\widetilde{\mathcal{G}}$ in any order compatible with inclusion of subvarieties. % (i.e. if $Z_1 \subseteq Z_2$ are two of the subvarieties, the blow-up corresponding to $Z_1$ comes before than the one corresponding to $Z_2$).
Since the exceptional divisor $E^{\mathbf{a}}_{[n]}$ of the blow-up $\operatorname{Bl}_{\Delta} \CC^\mathbf{a} \rightarrow \CC^\mathbf{a}$ is maximal in $\widetilde{\mathcal{G}}$, we can order the iterated blow-up process so that the last blow-up step is that of a dominant transform of $E^{\mathbf{a}}_{[n]}$. Such dominant transform is a Cartier divisor, and then that last blow-up step can be omitted without changing the resulting variety $\mathbb{C}_\mathbf{w}^{[\mathbf{a}]}$. 
Hence we can consider the building set for the wonderful compactification $\mathbb{C}_\mathbf{w}^{[\mathbf{a}]} \rightarrow \operatorname{Bl}_{\Delta} \CC^\mathbf{a}$ to be $\widetilde{\mathcal{G}}' = \widetilde{\mathcal{G}} \setminus \{ E^{\mathbf{a}}_{[n]} \}$.
Moreover, since $\Delta$ is contained in all elements in $\mathcal{G}$, then all the dominant transforms that yield the elements in $\widetilde{\mathcal{G}}'$ from elements in $\mathcal{G} \setminus \{\Delta\}$ are actually strict transforms. 

Since $\CC^\mathbf{a}$ is an affine space and all elements in $\mathcal{G} \setminus \{\Delta\}$ are linear subspaces containing $\Delta$, 
then Lemma~\ref{lemma.intersections.linear.spaces.with.exceptional.divisor} computes their intersections with $\bb{P}((\bb{C}^\mathbf{a}/\Delta)^\vee) \cong \mathbb{P}^N$ seen as the fiber over $0 \in \Delta$ of the exceptional divisor of $\operatorname{Bl}_{\Delta} \CC^\mathbf{a} \rightarrow \CC^\mathbf{a}$.
For any $I \subsetneq [n]$ with $|I| \geq 2$ such that $\Delta_I^{\mathbf{a}}$ is in $\mathcal{G} \setminus \{\Delta\} = \mathcal{K}_\mathbf{w}^{\mathbf{a}} \setminus \{\Delta\}$,
the intersection of the strict transform $\widetilde{\Delta_I^{\mathbf{a}}}$ of $\Delta_I^{\mathbf{a}}$ in $\operatorname{Bl}_{\Delta} \CC^\mathbf{a}$
with $\bb{P}((\bb{C}^\mathbf{a}/\Delta)^\vee) \cong \mathbb{P}^N$ is precisely $\delta_I^{\mathbf{a}}$ by Lemma~\ref{lemma.intersections.linear.spaces.with.exceptional.divisor}.  

Therefore, the set of intersections of the elements of the building set $\widetilde{\mathcal{G}}'$ with $\bb{P}^N=\bb{P}^{a_1+\cdots+a_{(n-1)}-1}$ 
is precisely the building set $\cal{H}_\mathbf{w}^\mathbf{a}=\{\delta_I^\mathbf{a}\,\vert\, I\subsetneq [n],\, \mathbf{w}_I>1\}$ in $\bb{P}^N$.
From the diagram (\ref{commutative.diagram.wonderful.Twa}) with Cartesian squares, we have that $T_\mathbf{w}^\mathbf{a}$ is the strict transform of $\bb{P}^N$ in $\mathbb{C}_\mathbf{w}^{[\mathbf{a}]}$ under the iterated blow-up $\mathbb{C}_\mathbf{w}^{[\mathbf{a}]} \rightarrow \operatorname{Bl}_{\Delta} \mathbb{C}^{\mathbf{a}}$.  
Therefore by Lemma~\ref{lemma.wonderful.subvariety}, $T_\mathbf{w}^\mathbf{a}$ is equal to the wonderful compactification of the building set $\cal{H}_\mathbf{w}^\mathbf{a}=\{\delta_I^\mathbf{a}\,\vert\, I\subsetneq [n],\, \mathbf{w}_I>1\}$ in $\bb{P}^N=\bb{P}^{a_1+\cdots+a_{(n-1)}-1}$, and this completes the proof.
\end{proof}

%%%%%%%%%%%%%%%%%%%%%%%%%%%%%%%%%%%%%%%%%%%%%%%%%%%%%%%%%%%%%%%%%%

\begin{corollary}
   Consider a cage $\mathbf{a}=(a_1,\dots,a_n)$ and $\mathbf{w}\in\cal{D}^{\operatorname{T}}_n$. Then, $T_\mathbf{w}^{\mathbf{a}}$ is the iterated blow-up of $\bb{P}^{a_1+\cdots+a_{(n-1)}-1}$ along the (strict transform of the) elements of $\cal{H}_\mathbf{w}^\mathbf{a}$ in any order compatible with inclusion of subvarieties. In particular,  $T_\mathbf{w}^\mathbf{a}$ is a toric variety if $\mathbf{w}$ is such that $\mathbf{w}_{[n-1]}\leq 1$.
\end{corollary}
\begin{proof}
 By Proposition~\ref{proposition:taw_is_wonderful}, we know that $T_\mathbf{w}^{\mathbf{a}}$ is a wonderful compactification. So, the claim about the iterated blow-up follows from the general theory of Li \cite[Theorem~1.3]{li2009wonderful}. The claim about the case  $\mathbf{w}_{[n-1]}\leq 1$ follows from the fact that such condition implies $n\in I$ if $\delta_I^\mathbf{a} \in \cal{H}_\mathbf{w}^\mathbf{a}$. By \eqref{PP: delta I definition}, such locus $\delta_I^\mathbf{a}$ is a torus-invariant subvariety of $\bb{P}^{a_1+\cdots+a_{(n-1)}-1}$ under the standard toric structure of projective space.
\end{proof}

\begin{lemma}
    Consider the cages $\mathbf{a}=(a_1,\dots,a_n)$ and $\mathbf{d}=(a_1,\dots,a_1)$, 
    and let $\mathbf{w}\in\mathcal{D}^{\operatorname{T}}_n$. 
    Then, there exists a unique closed immersion $f_{\mathbf{a}}:T^{\mathbf{a}}_{\mathbf{w}} \hookrightarrow T^{\mathbf{d}}_{\mathbf{w}}$ 
    which makes the following diagram commute:  
    \[
    \begin{tikzcd}
        T^{\mathbf{a}}_{\mathbf{w}} 
        \ar[r,hook, "f_{\mathbf{a}}"]\ar[d, "\theta_{\mathbf{a}}"'] & 
          T^{\mathbf{d}}_{\mathbf{w}} 
        \ar[d, "\theta_{\mathbf{d}}"]
        \\
        \bb{P}^{a_1+\cdots+a_{(n-1)}-1}\ar[r,hook,"g_{\mathbf{a}}"] & \bb{P}^{a_1(n-1)-1}
    \end{tikzcd}
    \]
where $\theta_{\mathbf{a}}$ and $\theta_{\mathbf{d}}$ are the wonderful compactification morphisms  
and $g_{\mathbf{a}}$ is the inclusion.
\end{lemma}

\begin{proof}
We define $N:=a_1+\cdots+a_{(n-1)}-1$ and recall our convention that $a_1 \geq \cdots \geq a_{n}$. 
Let us recall the description of the inclusion $g_{\mathbf{a}}:\bb{P}^{N} \hookrightarrow \bb{P}^{a_1(n-1)-1}$.  
We have inclusions $\CC^{a_i} \hookrightarrow \CC^{a_1}$ as the span of the first $a_i$ coordinates, for each $i$. 
By taking the product, we get an inclusion of affine spaces 
$\mathbb{C}^{a_1+\cdots+a_{n-1}} \hookrightarrow \mathbb{C}^{a_1(n-1)}$, which induces the closed immersion 
$g_{\mathbf{a}}:\mathbb{P}^N \hookrightarrow \mathbb{P}^{a_1(n-1)-1}$ between their projectivizations. 

Consider any $I_0 \subsetneq [n]$. 
Notice that $\delta_{I_0}^{\bf{a}}$ is in $\mathcal{H}_\mathbf{w}^{\mathbf{a}}$ if and only if 
$\delta_{I_0}^{\bf{d}}$ is in $\mathcal{H}_\mathbf{w}^{\mathbf{d}}$.  
Also, $\delta_{I_0} = \delta_{I_0}^{\bf{d}} \subseteq \bb{P}^{a_1(n-1)-1}$ intersects $\bb{P}^{N}$, since even $\delta_{[n]} \subseteq \delta_{I_0}$ clearly intersects $\bb{P}^{N}$.  
Moreover, $\delta_{I_0}^{\bf{a}}  =  \bb{P}^{N} \cap \delta_{_0}^{\bf{d}}$ scheme-theoretically.

Then, the intersections with $Y = \bb{P}^{N}$ of the varieties $\delta_I^{\mathbf{a}}$ in the building set $\mathcal{G}=\mathcal{H}_\mathbf{w}^{\mathbf{d}}$ in $X=\bb{P}^{a_1(n-1)-1}$ form the building set $\mathcal{H}=\mathcal{H}_\mathbf{w}^{\mathbf{a}}$ in $\bb{P}^{N}$. 
By applying Lemma \ref{lemma.wonderful.subvariety} in this setting, we obtain the existence and uniqueness of the desired closed immersion $f_{\mathbf{a}}:T^{\mathbf{a}}_{\mathbf{w}} \hookrightarrow T^{\mathbf{d}}_{\mathbf{w}}$, making the diagram commute. 
 \end{proof}

The variety $T_\mathbf{w}^\mathbf{a}$ acquires a ``universal" family from its embedding into $T^\mathbf{d}_\mathbf{w}$. 
\begin{definition}
    The ``universal" family $(T_\mathbf{w}^\mathbf{a})^+$ of $T_\mathbf{w}^\mathbf{a}$ is the pullback of the ``universal" family $\phi:(T^\mathbf{w}_{d,n})^+\to T^\mathbf{d}_\mathbf{w}$ along the closed immersion $T_\mathbf{w}^\mathbf{a}\hookrightarrow T^\mathbf{d}_\mathbf{w}$. 
    That is, the following diagram is Cartesian:
    \[
        \begin{tikzcd}
            (T_\mathbf{w}^\mathbf{a})^+\ar[r]\ar[d] & (T^\mathbf{d}_\mathbf{w})^+\ar[d] \\
            T_\mathbf{w}^\mathbf{a}\ar[r] & T^\mathbf{d}_\mathbf{w}.
        \end{tikzcd}
    \]
    In particular, $(T^\mathbf{a}_\mathbf{w})^+  \rightarrow T^\mathbf{a}_\mathbf{w}$ is flat and it is equipped with $n$ sections $\sigma_i :  T^\mathbf{a}_\mathbf{w} \to (T^\mathbf{a}_\mathbf{w})^+$
    whose images give the marked points.
\end{definition}

For the description of the fibers in the ``universal" family $(T^\mathbf{d}_\mathbf{w})^+ \to T^\mathbf{d}_\mathbf{w}$, we refer the reader to \cite[Section~2.4]{gallardo2017wonderful}.
Let us remark that there is a Zariski open $\cal{V}^\mathbf{a}_{\mathbf{w}}$ such that the geometric fibers of $(T_\mathbf{w}^\mathbf{a})^+\to T_\mathbf{w}^\mathbf{a}$ over $\cal{V}^\mathbf{a}_{\mathbf{w}}$ are all isomorphic to a flag of projective spaces 
$\mathbb{P}^{a_n} \subseteq \cdots \subseteq \mathbb{P}^{a_1}$ with a flag of hyperplanes $H_n \subseteq \cdots \subseteq H_1$ and $n$ points away from such hyperplanes -- that is, all points are in a flag of affine linear spaces
$\mathbb{C}^{a_n} \subseteq \cdots \subseteq \mathbb{C}^{a_1}$.
This configuration of $n$ is defined up to global scaling and translation by a vector in $\mathbb{C}^{a_n}$.  The study of the degenerations of these fibers will be part of forthcoming work.

\begin{remark}
Here is an alternative interpretation of the moduli problem that we study in this section. 
   If we consider $\mathbb{P}^{a_1} \supseteq \mathbb{C}^{a_1}$, 
   it is straightforward to verify that the group of automorphisms of $\mathbb{P}^{a_1}$ that fix the vector spaces in the flag $ 
\CC^{a_n} \subseteq \CC^{a_{n-1}} \subseteq \cdots \subseteq \CC^{a_1}
$ as sets and the hyperplane at infinity $\mathbb{P}^{a_1}\setminus\mathbb{C}^{a_1}$ pointwise, is the same as the subgroup of the automorphisms of  $\mathbb{C}^{a_1}$ generated by scalings and translations by the elements of $\CC^{a_n}$. 
Then, $T_\mathbf{w}^\mathbf{a}$ provides a compactification of the space of weighted configurations of distinct points in a flag of affine spaces up to either of these equivalent equivalence relations.   
    \end{remark}

\begin{remark}
In the case of a constant cage $\mathbf{d}=(d,\dots,d)$, 
the moduli spaces $T_{\mathbf{w}}^\mathbf{d}$
are the moduli spaces of pointed, stable rooted trees studied 
for weights equal to one in \cite{chen2009pointed} and for general weights in \cite{gallardo2017wonderful}. 
Let us mention how our notation relates to these works. The space $T_\mathbf{w}^\mathbf{a}$ with $\mathbf{w} = (1, \ldots, 1)$ and 
$\mathbf{a} = (d, \ldots, d)$ is denoted as $T_{d,n}$ in 
\cite{chen2009pointed}. The space 
$T_\mathbf{w}^\mathbf{a}$ with $\mathbf{w} \in \cal{D}^{T}_n$ and 
$\mathbf{a} = (d, \ldots, d)$ is denoted as $T_{d,n}^{\mathbf{w}}$ in \cite{gallardo2017wonderful}.
\end{remark}

\subsection{The Losev-Manin compactification and the polypermutohedral variety}\label{sec:toricTaw}

Among the compactifications constructed in Section 
\ref{subsec: taw}, there is one that is a higher-dimensional analog to the Losev-Manin compactification of $M_{0,n}$. 
In this subsection we prove that this generalization, which we denote as $T^\mathbf{a}_{LM}$, is isomorphic to the polypermutohedral variety $P_{\mathbf{a}^-}$ (Definition~\ref{defn:Polypermutohedral}) with cage $\mathbf{a}^-=(a_1,\dots,a_{n-1})$.

\begin{definition}
    Consider the cage $\mathbf{a}=(a_1,\dots,a_n)$, and let $\mathbf{w}=(w_1,\dots,w_n)\in\cal{D}^{\operatorname{T}}_n$ be such that $w_{[n-1]}\leq 1$ and $w_i+w_n>1$ for all $i\in[n-1]$. Then, we call the space $T^\mathbf{a}_\mathbf{w}$ the \emph{Losev-Manin compactification} for this moduli problem. We denote it as $T^\mathbf{a}_{LM}$.
\end{definition}

\begin{theorem}\label{PP: proof that toric FM agrees with PP}
    Consider the cages $\mathbf{a}=(a_1,\dots,a_n)$ and $\mathbf{a}^-=(a_1,\dots,a_{n-1})$. Then, the Losev-Manin compactification $T^{\mathbf{a}}_{LM}$ is isomorphic to the polypermutohedral variety with cage $\mathbf{a}^-$. This is,
    \[
    \operatorname{T}^{\mathbf{a}}_{LM}    \cong P_{(\mathbf{a}^-)}.
    \]
\end{theorem}

When $\mathbf{a}=\mathbf{1}^n=(1,\dots,1)\in\bb{Z}_{>0}^n$, the polypermutohedron $P_{(a^-)}=P_{\mathbf{1}^{(n-1)}}$ coincides with the standard permutohedron of dimension $n-2$. On the other hand, the toric variety corresponding to the $m$-dimensional permutohedron is isomorphic to the Losev-Manin compactification $\overline{M}_{0,m+3}^{LM}$ of $M_{0,m+3}$; see \cite{losev2000new}. This shows the following.

\begin{corollary}\label{PP: losev manin as polypermutohedron}
    Let $\mathbf{1}^n=(1,\dots,1)\in\bb{Z}_{>0}^n$. Then, there is an isomorphism
    \[
        T_{LM}^{\mathbf{1}^n} \cong P_{\mathbf{1}^{(n-1)}} \cong
        \overline{M}_{0,n+1}^{LM}.
    \]  
\end{corollary}

In order to prove Theorem~\ref{PP: proof that toric FM agrees with PP} let us introduce the necessary results from the theory of combinatorial building sets. We mainly follow \cite{zelevinsky-nested,postnikov-reiner-williams}. 

\begin{definition}
     A \emph{building set} $\mathcal{B}$ on a finite set $E$ is a collection of nonempty subsets of $E$ such that:
    \begin{enumerate}
        \item $\{i\}\in\mathcal{B}$ for all $i\in E$.
        \item If $I,J\in\mathcal{B}$ and $I\cap J\neq\emptyset$, then $I\cup J\in\mathcal{B}$.
    \end{enumerate}
    
    Denote by $\mathcal{B}_{max}\subseteq \mathcal{B}$ the set of maximal (with respect to inclusion) elements of $\mathcal{B}$. The building set $\mathcal{B}$ is said to be \emph{connected} if $\mathcal{B}_{max}=\{E\}$.
\end{definition}

\begin{definition}
    Given a building set $\mathcal{B}$, call a subset $N\subseteq\mathcal{B}-\mathcal{B}_{max}$ \emph{nested} if it satisfies the following condition: for any $k\geq 2$ and any $S_1,\dots, S_k\in N$ such that none of the $S_i$ contains another one, the union $S_1\cup\cdots\cup S_k$ is not in $\mathcal{B}$. 
\end{definition}

Let $\mathcal{B}$ be a connected building set on $E$. The \emph{nested fan associated to $\mathcal{B}$} is the fan $\Sigma(\mathcal{B})$ in $\bb{R}^E/\bb{R}e_{E}$ whose cones are
    \[
        \sigma_N:=\{\cone(\overline{e}_{S_1},\dots,\overline{e}_{S_k})\,\vert\,N=\{S_1,\dots,S_k\}\text{ is a nested set of }\mathcal{B}\}.
    \]
The fan $\Sigma(\cal{B})$ is the normal fan of a smooth polytope $P_\cal{B}$, known as the \emph{nestohedron} of $\mathcal{B}$; see \cite[Theorem~6.1]{zelevinsky-nested}. In particular, this fan is unimodular with respect to the image lattice $\bb{Z}^E\subseteq\bb{R}^E$ in $\bb{R}^E/\bb{R}e_{E}$; see \cite[Theorem~5.1~and~Corollary~5.2]{zelevinsky-nested}.   

The nested fan associated to a connected building set can be constructed via an iterated process of star subdivision of $\Sigma_{E}$ as follows:

\begin{construction}[{\cite[Theorem~4]{feichtner-yuzvinsky}}]\label{bg: blow up construction}
    Start with the fan $\Sigma_{E}$ in $\bb{R}^E/\bb{R}e_{E}$ whose rays have been labeled by the singletons $i$ for $i\in E$. Given a nonsingleton set $I\in\mathcal{B}$, define the cone $\sigma_I:=\cone(\overline{e}_i\,\vert\, i\in I)\in\Sigma_{E}$. Finally, define a total order of $\mathcal{B}=\{I_1,\dots,I_k\}$ compatible with inclusion, so that larger sets come before smaller sets. Then, define the fan $\Sigma(\mathcal{B})$ as that one obtained from $\Sigma_{E}$ after performing iterated star subdivisions along the cones $\sigma_I$ as $I$ ranges through the elements of $\mathcal{B}$ in the order we defined.
\end{construction}

\begin{example}
    The building set $\cal{B}_{LM}:=2^{[n]}\setminus\emptyset$ is known as \emph{boolean} building set. Its corresponding nestohedron is the $(n-1)$-dimensional permutohedron. 
\end{example}

\begin{example}\label{PP: def graphical building sets}
    A large class of examples of building sets and nestohedra is given by the so-called \emph{graphical building sets} and their corresponding \emph{graph associahedra}, first introduced in \cite{carr2006coxeter}. Given a connected graph $G$ with vertex set $E$, its corresponding graphical building set is
    \[
        \cal{B}(G) = \{T\subseteq[n]\,\vert\, G|_T\text{ is a connected subgraph of }G\}.
    \]
    Examples include permutohedra, when $G$ is a complete graph, and stellahedra, when $G$ is a star graph.
\end{example}

\begin{definition}
    Consider a building set $\cal{B}$ over $E$, we denote by $X(\cal{B})$ the toric variety constructed from the nested fan $\Sigma(\cal{B})$.
\end{definition}

Let us now introduce the notion of \emph{pullback} of a combinatorial building set $\mathcal{B}$ over $E$ along a surjective function $\pi:A\to E$ of finite sets. This construction already appears in \cite{crowley2022bergman}, and is equivalent to Eur and Larson's definition of \emph{expansion} for base polytopes of a polymatroid. As particular examples, we shall prove that every polypermutohedron is the pullback of the boolean building set $\cal{B}_{LM}=2^{[n]}\setminus\emptyset$.
\begin{definition}
    Consider a building set $\mathcal{B}$ on $E$, and let $\pi:A\to E$ be a surjective function of finite sets. The \emph{pullback of $\mathcal{B}$ along $\pi$} is the building set over $A$ defined as
    \[
        \pi^*\mathcal{B} = \{\{i\}\,\vert\,i\in A\}\cup\{\pi^{-1}(I)\,\vert\,I\in\mathcal{B}\}.
    \]
\end{definition}

The verification that $\pi^*\cal{B}$ is a building set is direct, and left to the reader. If $\pi$ is surjective, then $(\pi^*\mathcal{B})_{max}=\pi^{-1}(\mathcal{B}_{max})$, and hence $\cal{B}$ is connected if and only if $\pi^*(\cal{B})$ is connected.

Let us compute explicitly the nested sets of the pullback $\pi^*\mathcal{B}$, as well as the cones of its nested fan $\Sigma(\pi^*\cal{B})$. Recasting Construction~\ref{bg: blow up construction} in this context shows that $\Sigma(\pi^*\cal{B})$ is obtained from the normal fan of the simplex $\Sigma_A$ by a series of star subdivisions.

\begin{definition}\label{PP: def pi pair}
    Consider a building set $\mathcal{B}$ on $E$, and let $\pi:A\to E$ be a surjective function of finite sets. A \emph{$\pi$-pair} $(I, N)$ consists of a subset $I\subsetneq A$, not containing any fibers of $\pi$, and a nested set $N$ of $\mathcal{B}$. Here both cases $I=\emptyset$ and $N=\emptyset$ are allowed.
\end{definition}

\begin{lemma}\label{bg: lemma cones}
    Let $\pi^*\mathcal{B}$ be the pullback of the building set $\mathcal{B}$ on $E$ along $\pi:A\to E$. Then, there is a bijection between $\pi$-pairs and nested sets of $\pi^*\cal{B}$. Explicitly, every nested set of $\pi^*\mathcal{B}$ is of the form
    \[
        N_{(I,N)} = \{\{i\}\,\vert\, i\in I\}\cup\{\pi^{-1}(S)\,\vert\,S\in N\},
    \]
    for a unique $\pi$-pair $(I,N)$. Conversely, every $\pi$-pair defines a unique nested set in this manner.
\end{lemma}

\begin{proof}
    Consider a nested set $M$ of $\pi^*\mathcal{B}$. If $|M|=1$, then the result is direct, so let us consider the case $|M|\geq 2$. Then, by definition, $M\subseteq 2^A\setminus\emptyset$ consists of a collection of singletons and sets of the form $\pi^{-1}(S)$ for subsets $S\in\cal{B}$. Let $I$ denote the collection of singletons in $M$ that are not full fibers of $\pi$, in other words, all singletons excluding those that are a full fiber. If $I$ were to contain any fiber $\pi^{-1}(i)$ for some $i\in E$, this fiber would have size at least two, the collection of all singletons making up this fiber would be disjoint, and their union $\pi^{-1}(i)$ would be an element of $\pi^*\mathcal{B}$ with at least two elements, which contradicts $M$ being a nested set. Therefore, $I$ cannot contain any full fibers of size at least two. 
    
    On the other hand, consider any $k\geq 2$ elements $\pi^{-1}(S_1),\dots,\pi^{-1}(S_k)\in M$, so that in particular $S_1,\dots,S_k\in\cal{B}$. We remark that in particular, it is possible for some $\pi^{-1}(S_i)$ to be a singleton. Then, $\pi^{-1}(S_i)\subseteq\pi^{-1}(S_j)$ if and only if $S_i\subseteq S_j$. Similarly, $\pi^{-1}(S_1)\cup\cdots\cup\pi^{-1}(S_k)\in\pi^*\mathcal{B}$ if and only if $S_1\cup\cdots\cup S_k\in\mathcal{B}$. Therefore, the collection of all subsets $S\in\cal{B}$ such that $\pi^{-1}(S)\in M$ forms a nested set $N$ of $\cal{B}$. 
    
    Conversely, given a $\pi$-pair $(I,N)$, it is a direct verification to check that $N_{(I,N)}:=\{\{i\}\,\vert\, i\in I\}\cup\{ \pi^{-1}(S)\,\vert\, S\in N\}$ defines a nested set of $\pi^*\mathcal{B}$.
\end{proof}

\begin{corollary}
    Given a $\pi$-pair $(I,N)$, the cone corresponding to the nested set $N_{(I,N)}$ in $\Sigma(\pi^*\cal{B})$ is
    \[
        \sigma_{(I,N)} := \cone(\bar{e}_i\,\vert\, i\in I)+\cone(\bar{e}_{\pi^{-1}(S)}\,\vert\, S\in N)\subseteq\bb{R}^A/\bb{R}e_A.
    \]
    Moreover, every cone in the nested fan of $\pi^*\mathcal{B}$ is of this form.
\end{corollary}

\begin{corollary}\label{pp: lemma cone decomposition pullback}
    Given a $\pi$-pair $(I,N)$, the cone $\sigma_{(I,N)}\in\Sigma(\pi^*\cal{B})$ can be decomposed as
    \[
        \sigma_{(I,N)} = \sigma_{(I,\emptyset)}+\sigma_{(\emptyset, N)}.
    \]
\end{corollary} 

With these results at hand, we can prove that polypermutohedra are pullbacks of the boolean building set. 

\begin{proposition}\label{PP: PP is pullback of boolean bs}
    Let $\cal{B}_{LM}=2^{[n]}\setminus\emptyset$ be the boolean building set, and let $\pi:A\to [n]$ be a caging with cage $\mathbf{a}$. Then, the nested fan $\Sigma(\pi^*\cal{B}_{LM})$ is the inner normal fan of the polypermutohedron with cage $\mathbf{a}$.
\end{proposition}

\begin{proof}
    It is well-known that nested sets of the boolean building set $\cal{B}_{LM}$ over $[n]$ are in bijection with flags of nonempty subsets $S_1\subsetneq S_2\subsetneq\cdots\subsetneq S_k\subsetneq[n]$. Indeed, a nested set of $\cal{B}_{LM}$ cannot contain two incomparable sets, as their union is also an element of $\cal{B}_{LM}$. Therefore, in light of Corollary~\ref{pp: lemma cone decomposition pullback}, the cones of $\pi^*\cal{B}_{LM}$ are
    \[
        \sigma_{(I,N)} = \sigma_{(I,\emptyset)}+\sigma_{(\emptyset, N)},
    \]
    where $N=\{S_1\subsetneq S_2\subsetneq\cdots\subsetneq S_k\subsetneq[n]\}$. This is precisely the Bergman fan of the boolean polymatroid with cage $\pi$, as introduced in \cite[Definition~1.2]{crowley2022bergman}. It is proved in \cite[Appendix A]{crowley2022bergman} that this is the inner normal fan of the polypermutohedron.
\end{proof}

\begin{definition}
    A polymatroid $P=(A,f)$ is said to be \emph{boolean} if its base polytope is a polypermutohedron. In other words, if its base polytope is the nestohedron corresponding to $\pi^*\cal{B}_{LM}$ for some caging $\pi:A\to[n]$.
\end{definition}

\begin{proof}[Proof of Theorem~\ref{PP: proof that toric FM agrees with PP}]
    Recall that the compactification $T^{\mathbf{a}}_\mathbf{w}$ is the iterated blow-up of the projective space $\bb{P}^{a_1+\cdots+a_{(n-1)}-1}$ along the dominant transforms of the varieties in the geometric building set
    \[
        \cal{H}_\mathbf{w}^\mathbf{a} = \left\{ \delta_I^\mathbf{a}\,\vert\,\mathbf{w}_I>1\right\}.
    \]  
    When $\mathbf{w}$ is such that $\mathbf{w}_{[n-1]}\leq 1$ and $w_i+w_n>1$ for all $i\in[n-1]$, every subvariety in $\cal{H}_\mathbf{w}^\mathbf{a}$ is of the form
    \[
        \delta_I^\mathbf{a} = \bigcap_{i\in I\setminus\{n\}} V(x_{i1},\dots,x_{ia_i}),
    \]  
    for some $I\subsetneq [n]$ such that $n\in I$ and $|I|\geq 2$; see Equation~\eqref{PP: delta I definition}. These subvarieties are invariant under the standard torus action on $\bb{P}^{a_1+\cdots+a_{(n-1)}-1}$. 
    
    On the other hand, let $\cal{B}_{LM}=2^{[n-1]}\setminus\emptyset$ be the boolean building set over $[n-1]$, and $\rho:A^-\to[n-1]$ be a caging with cage $\mathbf{a}^-$. Then, for every $I$ as above, the cone $\sigma_{\rho^{-1}(I\setminus \{n\})}\in\Sigma(\rho^*\cal{B}_{LM})$ corresponds to the subvariety $\delta_I^\mathbf{a}$ under the toric orbit-cone correspondence. Moreover, every cone in $\Sigma(\rho^*\cal{B}_{LM})$ is of this form by definition. Therefore, both constructions can be obtained as iterated blow-ups of $\bb{P}^{a_1+\cdots+a_{(n-1)}-1}$ along the same subvarieties, in an order that can be assumed to be the same, so they are the same.
\end{proof}

\subsection{Pullback of building sets and nontrivial toric fibrations}
\label{sec:FibrationPP}

Throughout this section we fix a caging $\pi:A\to[n]$ with cage $\mathbf{a}=(a_1,\dots,a_n)$. The relation between the nested fan of a building set $\cal{B}$ and that one of its pullback $\pi^*\cal{B}$ gives rise to interesting combinatorics. A witness of this is the fact that the toric variety constructed from the pullback fan is a nontrivial, locally trivial fibration over a product of projective spaces with fibers corresponding to the base building set. Define the subsets $A_i:=\pi^{-1}(i)\subseteq A$ for all $i=1,\dots,n$, so that $|A_i|=a_i$.

\begin{theorem}\label{fib: theorem fibration}
    Let $\Sigma(\mathcal{B})$ and $\Sigma(\pi^*\cal{B})$ denote the nested fans of $\mathcal{B}$ and $\pi^*\mathcal{B}$, respectively. Then, there is a nontrivial splitting of the fan $\Sigma(\pi^*\cal{B})$,
    \[
        0\to \Sigma(\mathcal{B})\to \Sigma(\pi^*\cal{B})\to \Sigma_{A_1}\times\cdots\times\Sigma_{A_n}\to 0,
    \]
    where $\Sigma_{A_i}$ is the normal fan to the $(a_i-1)$-dimensional simplex in $\bb{R}^{A_i}/\bb{R}\cdot e_{A_i}$. In particular, there is a nontrivial, locally trivial toric fibration of projective toric varieties
    \[
    X(\pi^*\mathcal{B})\to \bb{P}^{A_1}\times\cdots\times\bb{P}^{A_n}
    \]
    with fibers isomorphic to $X(\mathcal{B})$.
\end{theorem}

Let us recall the definition of a splitting of a fan. 

\begin{definition}[{\cite[Definition~3.3.18]{Cox-Little-Schenck-Toric-Varieties-Book}}]\label{fib: cls definition splitting}
    Consider three lattices $N_0,N$ and $N'$ and a short exact sequence
    \[
    \begin{tikzcd}[column sep = small]
        0\ar[r]&
        N_0\ar[r,"\iota"]&
        N\ar[r,"\varphi"]&
        N'\ar[r]&
        0.
    \end{tikzcd}
    \]
    Suppose that $\Sigma_0,\Sigma,\Sigma'$ are fans in $(N_0)_\bb{R}, N_\bb{R}$ and $N'_\bb{R}$ and the maps $\iota$ and $\varphi$ are compatible with them. We say that $\Sigma$ is \emph{split by $\Sigma_0$ and $\Sigma'$}, or that the sequence defines a \emph{splitting of $\Sigma$}, if there exists a subfan $\widehat{\Sigma}$ of $\Sigma$ such that
    \begin{enumerate}
        \item Given cones $\widehat{\sigma}\in\widehat{\Sigma}$ and $\sigma_0\in\Sigma_0$, the sum $\widehat{\sigma}+\sigma_0$ is in $\Sigma$, and every cone of $\Sigma$ arises in this way.
        \item $\varphi_\bb{R}$ maps every cone in $\widehat{\Sigma}$ bijectively to a cone $\sigma'\in\Sigma'$ such that $\varphi(\widehat{\sigma}\cap N)=\sigma'\cap N'$. Furthermore, the map $\widehat{\sigma}\mapsto\sigma'$ defines a bijection between $\widehat{\Sigma}$ and $\Sigma'$.
    \end{enumerate}
    
\end{definition}

\begin{proposition}[{\cite[Theorem~3.3.19]{Cox-Little-Schenck-Toric-Varieties-Book}}]\label{fib: CLS fibration result}
    If $\Sigma$ is split by $\Sigma_0$ and $\Sigma'$ as in Proposition~\ref{fib: cls definition splitting}, then $X(\Sigma)$ is a locally trivial fibration (fiber bundle) over $X(\Sigma')$ with fiber $X(\Sigma_0)$. In particular, all fibers of $X(\Sigma)\to X(\Sigma')$ are isomorphic to $X(\Sigma_0)$.
\end{proposition}

\begin{proof}[Proof of Theorem~\ref{fib: theorem fibration}]
    By Lemma~\ref{pp: lemma cone decomposition pullback}, every cone of $\Sigma(\pi^*\cal{B})$ corresponds to a $\pi$-pair $(I,N)$ and has the form
    \[
        \sigma_{(I, N)} = \sigma_{(I,\emptyset)}+\sigma_{(\emptyset, N)},
    \]
    where
    \[
        \sigma_{(I,\emptyset)}=\cone(\overline{e}_i\,\vert\, i\in I)\quad\text{and}\quad\sigma_{(\emptyset, N)}=\cone(\overline{e}_{\pi^{-1}(S)}\,\vert\, S\in N)\subseteq\bb{R}^A/\bb{R}e_A,
    \]
    and, by the definition of a $\pi$-pair, $I\subsetneq A$ does not contain any fiber of $\pi$ and $N$ is a nested set of $\mathcal{B}$.

    Consider the short exact sequence of vector spaces
    \[
    \begin{tikzcd}[column sep = small]
        0\ar[r]&
        \bb{R}^{[n]}/\bb{R}e_{[n]}\ar[r,"\iota"]&
        \bb{R}^{A}/\bb{R}e_{A}\ar[r,"\varphi"]&
        \bb{R}^{A_1}/\bb{R}e_{A_1}\times\cdots\times\bb{R}^{A_n}/\bb{R}e_{A_n}\ar[r]&
        0.
    \end{tikzcd}
    \]
     Here, the maps in the sequence are defined as $\iota(\overline{e}_i)=\overline{e}_{A_i}$ for all $i\in [n]$, and if $i\in A_i\subseteq A$, then $\varphi(\overline{e}_i)=(0,\dots,\overline{e}_i,\dots,0)$. One can check directly that these maps are well-defined and that the sequence is exact.

    Let us show that $\iota$  defines an inclusion of fans. Consider a cone $\sigma_N\in\Sigma(\cal{B})$ corresponding to a nested set $N$ of $\mathcal{B}$. Then, $\sigma_N=\cone(\overline{e}_S\,\vert\,S\in N)$, and its image under $\iota$ is precisely the cone $\sigma_{(\emptyset, N)}$ in $\bb{R}^A/\bb{R}e_A$. In particular, $\iota$ is an inclusion of fans.
    
    Next, let us see that $\varphi$ defines a surjective map of fans. Consider a cone $\sigma_{(I,N)}=\sigma_{(I,\emptyset)}+\sigma_{(\emptyset, N)}\in\Sigma(\pi^*\cal{B})$. By the exactness of the sequence, $\sigma_{(\emptyset, N)}$ is in the kernel of $\varphi$. On the other hand, the cone $\sigma_{(I,\emptyset)}$ maps to the cone $\sigma_{I\cap A_1}\times\cdots\times\sigma_{I\cap A_n}$ in $\Sigma_{A_1}\times\cdots\times\Sigma_{A_n}$. To show surjectivity note that every cone in the product fan is of the form $\sigma_{J_1}\times\cdots\times\sigma_{J_n}$ for some subsets $\emptyset\subseteq J_i\subsetneq A_i$ (the empty set corresponding to the zero cone). Then, $(I:=J_1\cup\cdots\cup J_n,\emptyset)$ is a $\pi$-pair because $I$ it does not contain any fiber of $\pi$. Moreover, we have that $\varphi(\sigma_{(I,\emptyset)})=\sigma_{J_1}\times\cdots\times\sigma_{J_n}$.

    Let us proceed to show that $\Sigma(\pi^*\cal{B})$ is split by $\Sigma(\cal{B})$ and $\Sigma_{A_1}\times\cdots\times\Sigma_{A_n}$. Condition (1) in Definition~\ref{fib: cls definition splitting} is an immediate consequence of Corollary~\ref{pp: lemma cone decomposition pullback} saying that every cone $\sigma_{(I, N)}$ can be decomposed as the sum $\sigma_{(I,\emptyset)}+\sigma_{(\emptyset, N)}$. To verify condition (2) we define the subfan $\widehat{\Sigma}$ of $\Sigma(\pi^*\cal{B})$ to be the fan consisting of all cones of the form $\sigma_{(I,\emptyset)}$ as $I$ ranges over all the subsets of $A$ not containing any fiber of $\pi$. Then, it follows directly from our discussion on the surjectivity of $\varphi$ that $\varphi$ defines a bijection between the cones of $\widehat{\Sigma}$ and those of $\Sigma_{A_1}\times\cdots\times\Sigma_{A_n}$. Moreover, since $\varphi$ restricts to a map of lattices, it follows that $\varphi(\sigma_{(I,\emptyset)}\cap\bb{Z}^A/\bb{Z}e_A)=\varphi(\sigma_{(I,\emptyset}))\cap\bb{Z}^{A_1}/\bb{Z}e_{A_1}\times\cdots\times\bb{Z}^{A_n}/\bb{Z}e_{A_n}$. 

    To conclude, let us prove that the fibration is not trivial. Consider any (vector space) splitting $\nu:\bb{R}^{A_1}/\bb{R}e_{A_1}\times\cdots\times\bb{R}^{A_n}/\bb{R}e_{A_n}\to\bb{R}^A/\bb{R}e_A$. Then, for the fibration to be trivial we must have that $\nu(\sigma_I)=\sigma_{(I,\emptyset)}$ for all $I=J_1\cup\cdots\cup J_n$ with $J_i\subsetneq A_i$. However, note that the one-dimensional cones of the form $\sigma_{(I,\emptyset)}$ span all of $\bb{R}^A/\bb{R}e_A$. This would imply that the image of $\nu$ is all of $\bb{R}^A/\bb{R}e_A$, which is impossible since the image of $\nu$ is the complement of $\RR^{[n]}/\bb{R}e_{[n]}\hookrightarrow\bb{R}^A/\bb{R}e_A$.
\end{proof}

From Theorem~\ref{fib: theorem fibration} together with Proposition~\ref{fib: CLS fibration result} we recover the following surprising result, initially proved for constant cages $\mathbf{a}=(d,\dots,d)$, $d\geq 1$, 
 and using geometric methods in \cite{ggg-higher-lm}.

\begin{corollary}\label{fib: corollary fibration PP}
    There is a nontrivial, locally trivial toric fibration of projective toric varieties
    \[
        T^{\mathbf{a}}_{LM}\to \bb{P}^{a_1-1}\times\cdots\times\bb{P}^{a_{(n-1)}-1}
    \]
    with fibers isomorphic to $\overline{M}^{LM}_{0,n+1}$, the Losev-Manin compactification of $M_{0,n+1}$.
\end{corollary}
\begin{proof}
    Recall from Corollary~\ref{PP: losev manin as polypermutohedron} that $\overline{M}^{LM}_{0,n+1}$ is the toric variety corresponding to the $(n-2)$-dimensional permutohedron. This permutohedron is the nestohedron defined by the building set $\cal{B}_{LM}=2^{[n-1]}\setminus\emptyset$, so that $\overline{M}^{LM}_{0,n+1}\cong X(\cal{B}_{LM})$. 
    
    On the other hand, by Theorem~\ref{PP: proof that toric FM agrees with PP}, we know that $T^{\mathbf{a}}_{LM}\cong P_{(\mathbf{a}^{-})}$, where $\mathbf{a}^-=(a_1,\dots,a_{n-1})$. Therefore, $T^{\mathbf{a}}_{LM}$ is isomorphic to the toric variety $X(\rho^*\cal{B}_{LM})$, where $\rho:A^{-}\to[n-1]$ is any caging with cage $\mathbf{a}^-$. The result then follows from Theorem~\ref{fib: theorem fibration}. 
\end{proof}

\subsection{Proof of Theorem~\ref{GFM: thm:weightedFM_forFlags}}
\label{sec:proofMainThmPolyPermetohedra}
The different proofs of our theorem are spread throughout the previous subsections. The following is meant to aid the reader in navigating those results.

\begin{proof}
    After constructing the compactifications $\mathbb{C}^{\mathbf{a}}_{\mathbf{w}}$ in Definition~\ref{GFM: definition.weighted.FM.for.flags}, the moduli spaces $T^{\mathbf{a}}_{\mathbf{w}}$ are defined in Subsection~\ref{subsec: taw}; see also Theorem~\ref{def:taw}. The claim that $T^{\mathbf{a}}_{\mathbf{w}}$ is smooth, and its boundary divisors are normal crossings is a direct consequence of Li's general theory of wonderful compactifications; see Proposition~\ref{prop: li blow-up building set}. 

    The proof of part (i) of the Theorem is given in Theorem~\ref{PP: proof that toric FM agrees with PP}. Part (iii) is Corollary~\ref{fib: corollary fibration PP}. Let us proceed to prove part (ii).

    Define $N=a_1+\cdots+a_{(n-1)}-1$. Let $\mathbf{w}_{LM}$ be any weight vector such that $w_{[n-1]}<1$ and $w_i+w_n$ for all $i\in[n-1]$, and let $\mathbf{w}=(1,\dots,1)$. Then, $T_{LM}^\mathbf{a}$ and $T^\mathbf{a}$ are the wonderful compactifications of $\bb{P}^{N}$ with respect to the geometric building sets $\cal{H}_{\mathbf{w}_{LM}}^\mathbf{a}$ and $\cal{H}_{\mathbf{w}}^\mathbf{a}$, respectively. In particular, since $\mathbf{w}_{LM}\leq\mathbf{w}$, we have that $\cal{H}_{\mathbf{w}_{LM}}^\mathbf{a}\subseteq \cal{H}_{\mathbf{w}}^\mathbf{a}$. 
    For simplicity let us denote these building sets as $\cal{H}_{LM}$ and $\cal{H}$. Suppose that
    \[
        \cal{H} = \{U_1,\dots,U_p,V_1,\dots,V_q\}
    \]  
    is such that $\cal{H}_{LM}=\{U_1,\dots,U_p\}$ and $\cal{H} \setminus \cal{H}_{LM}=\{V_1,\dots,V_q\}$. Each of the subvarieties $U_k$ equals a $\delta_{I_k}^\mathbf{a}$ for some $I_k\subsetneq [n]$ such that $|I_k|\geq 2$ and $n\in I_k$. Similarly, each $V_k$ is of the form $\delta_{J_k}^\mathbf{a}$ for some $J_k\subsetneq [n]$ such that $|J_k|\geq 2$ and $n\not\in J_k$; see Equation~\eqref{PP: delta I definition}. We refer to the sets $I_k,J_k$ as the index sets of $U_k$ and $V_k$, respectively.

    Let us order $\cal{H}_{LM}$ so that the subvarieties with larger index sets appear before those with smaller index sets. Similarly, order $\cal{H}$ so that the elements of $\cal{H}_{LM}\subseteq \cal{H}$ come first---appearing in the order we just fixed---and are followed by the subvarieties $V_k$ ordered so that those subvarieties with larger index sets appear before those with smaller index sets. 

    We claim that these orders for $\cal{H}_{LM}$ and $\cal{H}$ both satisfy the conditions of Proposition~\ref{bg: li order of blow-up doesnt matter}, so that $T^\mathbf{a}_{LM}$ and $T^\mathbf{a}$ are both the iterated blow-ups of $\bb{P}^{N}$ along the building sets in the orders we just prescribed for their building sets. This would imply the result, as the $p$th step of the iterated blow-up process of $T^\mathbf{a}$ is precisely that one in the construction of $T_{LM}^\mathbf{a}$.

    The fact that $\cal{H}_{LM}$ in this order satisfies the conditions of Proposition~\ref{bg: li order of blow-up doesnt matter} is a direct consequence of the fact that, for any $1\leq r\leq p$, the subset $\{U_1,\dots,U_r\}\subseteq\cal{H}_{LM}$ is closed under intersections of its own elements. In particular, every such subset is a building set so the claim follows.

    Consider a subset $\cal{H}_r:=\{U_1,\dots,U_p,V_1,\dots,V_r\}\subseteq\cal{H}$ for some $1\leq r\leq q$, and let $\cal{S}_r$ be the induced arrangement of subvarieties of $\cal{H}_r$ in $\bb{P}^{N}$. Every element $S\in\cal{S}_r$ can be written as an intersection of the form $S= U_I\cap V_{J_1}\cap\cdots\cap V_{J_k}$ for some $0\leq k\leq r$, $I\subsetneq [n]$ such that $n\in I$ and $|I|\geq 2$, and subsets $J_i\subsetneq [n]$ such that $|J_i|\geq 2$ and $n\not\in J_i$. This follows from the fact that $\cal{H}_{LM}$ is closed under intersection. Moreover, as in the proof of Lemma~\ref{GFM: building set}, we can assume that all the subsets $I,J_1,\dots,J_k$ are disjoint. 
    
    In order to show $\cal{H}_r$ is a building set we prove that the intersection $U_I\cap V_{J_1}\cap\cdots\cap V_{J_k}$ is transverse. As in Lemma~\ref{GFM: building set}, this is a direct dimension count. Since all the subvarieties $U_I, V_{J_1},\cdots, V_{J_k}$ are linear, we can identify them with their tangent spaces. Then, transversality will follow from the equality
    \[
        \codim(U_I\cap V_{J_1}\cap\cdots\cap V_{J_k}, \bb{P}^{N}) = \codim(U_I,\bb{P}^N)+\sum_{i=1}^k\codim(V_{J_i},\bb{P}^N).
    \]
   This verification is left to the reader, as the dimension count is very similar to that one in the proof of Lemma~\ref{GFM: building set}, using the fact that all the subsets $I,J_1,\dots,J_k$ are disjoint.

   Since the previous argument works for all elements $S$ of the arrangement $\cal{S}_r$, we obtain that $\cal{H}_r$ is a building set. Since the argument holds for all values of $r$, the conditions of Proposition~\ref{bg: li order of blow-up doesnt matter} are met, so $T^\mathbf{a}$ is the iterated blow-up of $\bb{P}^N$ along the subvarieties in $\cal{H}$ in the order we described previously. In particular, we have a sequence of blow-up maps $T^\mathbf{a}\to T^\mathbf{a}_{LM}\to\bb{P}^N$, as desired.
\end{proof}
 
\subsection{Pullback of building sets and refinements}\label{PP: subsec refinement}

The data of a caging $\pi:A\to[n]$ is equivalent to that of a partition of the set $A$ with $n$ elements. Given two cagings $\pi:A\to[n]$ and $\pi':A\to[n']$, we say that $\pi$ \emph{refines} $\pi'$, written $\pi\succeq\pi'$, if for every $i\in[n]$ one has $\pi^{-1}(i)\subseteq\pi'^{-1}(i')$ for some $i'\in[n']$. Eur and Larson in \cite[Proposition~2.2]{eur2024intersection} leverage this observation to prove that given a fixed set $A$ and two cagings $\pi$ and $\pi'$ such that $\pi\succeq\pi'$, the nested fans of the two corresponding polystellahedra are related by an explicit series of star subdivisions. Here we prove that this phenomenon is also present for pullbacks of building sets. As noted by Eur and Larson, this is a direct consequence of \cite[Theorem~4.2]{feichtner-muller} and \cite[Proposition~2]{feichtner-yuzvinsky}.

\begin{lemma}\label{PP: polypermuto refinement}
    Let $\pi:A\to[n]$ and $\pi':A\to[n']$ be two cagings such that $\pi\succeq\pi'$. Consider the building sets $\cal{B}_{LM}=2^{[n]}\setminus\emptyset$ and $\cal{B}'_{LM}=2^{[n']}\setminus\emptyset$ on $[n]$ and $[n']$, respectively. Then, the fan $\Sigma(\pi^*\mathcal{B}_{LM})$ is obtained from $\Sigma(\pi'^*\mathcal{B}_{LM}')$ by a sequence of star subdivisions. Moreover, these subdivisions are unimodular. In particular, the toric variety $X(\pi^*\mathcal{B}_{LM})$ is obtained from $X(\pi'^*\mathcal{B}_{LM}')$ by a sequence of smooth blow-ups along torus-invariant subvarieties.
\end{lemma}
\begin{proof}
    The nonempty poset $\cal{L}=2^{A}$ is finite and there exists a greatest lower bound for any non-empty subset.  
    Hence, $\cal{L}$ is an atomic meet-semilattice (see \cite[Section~2]{feichtner-muller} for relevant definitions). 
    Since $\pi\succeq\pi'$, for every element $S'\in\cal{B}'_{LM}$ the subset $S=\pi(\pi'^{-1}(S'))\subseteq [n]$ is such that $\pi^{-1}(S) = \pi'^{-1}(S')$. Therefore, $\pi'^*(\cal{B}'_{LM})\subseteq\pi^*(\cal{B}_{LM})$. Then, the result follows from \cite[Theorem~4.2]{feichtner-muller}. The fact that the resulting fan of these subdivisions is unimodular follows from \cite[Proposition~2]{feichtner-yuzvinsky}.
\end{proof}

Using this result in the context of polypermutohedra we obtain the following.

\begin{proposition}\label{PP: corollary map from LM to TLMa}
    Consider a caging $\pi:A\to[n]$ with cage $\mathbf{a}=(a_1,\dots,a_n)$, and let $m=\sum_{i=1}^{n-1}a_i$. Then, there is a toric map
    \[
        \overline{M}^{LM}_{0,m+2}\to T_{LM}^\mathbf{a}
    \]
    which is the composition of blow-ups along smooth torus-invariant subvarieties of $T_{LM}^\mathbf{a}$.
\end{proposition}
\begin{proof}
    Define $A^{-}=\pi^{-1}([n-1])$ and $\rho=\pi|_{A^-}$. By Lemma~\ref{PP: polypermuto refinement}, there is a reduction map $X(P_{(\mathbf{1}^{m})})\to X(P_{(\mathbf{a}^-)})$. By Theorem~\ref{PP: proof that toric FM agrees with PP} and Corollary~\ref{PP: losev manin as polypermutohedron}, we have that $X(P_{\mathbf{1}^{m}})\cong T_{LM}^{\mathbf{1}^{(m+1)}}\cong\overline{M}_{0,m+2}^{LM}$. On the other hand, $X(P_{(\mathbf{a}^-)})\cong T_{LM}^\mathbf{a}$. The result follows.
\end{proof}

%%%%%%%%%%%%%%%%%%%%%%%%%%%%%%%%%%%%%%%%%%%%%%%%%%%%%
%%%%%%%%%%%%%%%%%%%%%%%%%%%%%%%%%%%%%%%%%%%%%%%%%%%%%

\section{Generalized Fulton-MacPherson compactification and polystellahedral varieties}\label{sec:FM_polystella}

This section is dedicated to the proof of Theorem~\ref{thm:mainPolyStella}. We begin with the construction of our moduli spaces in Section~\ref{sec:GFM_PS}. To establish the isomorphism between the moduli perspective and the combinatorial side, we introduce a novel presentation of the polystellahedral variety as a blow-up of a product of projective spaces, which is achieved in Section \ref{sec:PolyGeometry}. 
The presentation of the polypermutohedral variety as a quotient of the polystellahedral variety appears in Subsection~\ref{subsection.quotients}. The proof of the theorem synthesizes the key results from these subsections, and is given in Subsection~\ref{subsection.proof.main.thm.polystella}.

\subsection{Generalized Fulton-MacPherson compactification}\label{sec:GFM_PS}
Let $X$ be a complex nonsingular algebraic variety, and let $D$ be a nonsingular closed proper subvariety of $X$. In \cite{kim2009generalization} Kim and Sato introduced the following compactification $X_D^{[n]}$ of the configuration space of $n$ labeled points in $X \setminus D$, characterized by allowing the points to coincide, but not to meet $D$.

\begin{definition}[{\cite[Section~1.1]{kim2009generalization}}] \label{def:GFM}
    The compactification $X_D^{[n]}$ of $n$ points in $X$ away from $D$ is the closure of $X^n \setminus \bigcup_{S} D_{S}$ diagonally embedded in 
    \[
    X^n \times \prod_{\substack{S \subseteq [n]\\S\neq\emptyset}}\mathrm{Bl}_{D_{S}} X^n,
    \]
    where $D_{S} \subseteq X^n$ is the locus parametrizing the collection of points $x$ whose $i$th component $x_i$ is in $D$ if $i \in S \subseteq[n]$.
\end{definition}

As noted in \cite[Section~2.2]{kim2009generalization}, the collection 
\begin{align}\label{ps: ks building set}
 \cal{G}_{X,D}:=\{D_S\,\vert\,S \subseteq[n],S\neq\emptyset\}  
\end{align}
in $X^n$ is a building set in the sense of Li \cite{li2009wonderful}. It follows that $X_D^{[n]}$ is the wonderful compactification of $X^n$ with respect to $\cal{G}_{X,D}$. 
\begin{theorem}[{\cite[Theorem~1]{kim2009generalization}}]
\label{thm:GeneralizedFM}
Let $X$ and $D$ be as above, then the following holds:
\begin{enumerate}
    \item The variety $X_D^{[n]}$ is nonsingular.
    \item There is a ``universal" family $X_D^{[n]+} \to X_D^{[n]}$. It is a flat family of stable degenerations of $X$ with $n$ smooth labeled points away from $D$.
    \item The boundary $X_D^{[n]} \setminus (X^n \setminus \bigcup_{S \in \cal{G}_{X,D}} D_{S})$ is a union of divisors $\widetilde{D}_{S}$ corresponding to $D_{S}$, $|S| \geq 1$. Any set of these divisors intersects transversally.
\end{enumerate}
\end{theorem}
Let $\pi:A\to [n]$ be a caging with cage $\mathbf{a}=(a_1,\dots,a_n)$, where $a_i = |\pi^{-1}(i)|$. The main goal of this section is to generalize $X_{D}^{[n]}$
to the case of configurations of $n$ points in a flag $\mathbb{P}^{a_n}\subseteq\cdots\subseteq\mathbb{P}^{a_1}$, lying away from a subflag $H_n\subseteq\cdots\subseteq H_1$ of hyperplanes $H_i\subseteq\mathbb{P}^{a_i}$.

Consider the cages $\mathbf{a}$ and $\mathbf{d} = (d, \cdots, d)$, where $d=a_1$ and $a_1\geq\cdots\geq a_n$. Define $\bb{P}^\mathbf{a}:=\prod_{i=1}^n\bb{P}^{a_i}$. Throughout this section we fix a flag
\[
    \mathbb{P}^{a_n}\subseteq \mathbb{P}^{a_{n-1}} \subseteq \cdots \subseteq
    \mathbb{P}^{a_1}, 
\]
together with hyperplanes $H_i\subseteq\bb{P}^{a_i}$ for all $i\in [n]$ satisfying the compatibility condition $H_{i+1}=\mathbb{P}^{a_{i+1}}\cap H_{i}$ for all $i=1,\dots,n-1$. Observe that the subvarieties $\mathbb{C}^{a_i} \cong \mathbb{P}^{a}_i \setminus H_i$ form a flag
\[
    \mathbb{C}^{a_n}\subseteq\mathbb{C}^{a_{n-1}}\subseteq\cdots\subseteq\mathbb{C}^{a_1}.
\]
For simplicity we will consider $H_1=V(X_0)$, and we often write $H$ instead of $H_1$. Let $\mathbf{d}=(a_1,\dots,a_1)$, then the flag induces an embedding $\iota_{\mathbf{a}}:\mathbb{P}^{\mathbf{a}} \hookrightarrow \mathbb{P}^{\mathbf{d}}=\prod_{i=1}^n \mathbb{P}^d$
such that
\[
\iota_{\mathbf{a}}(\bb{P}^{a_1}\times\cdots\times H_i\times\cdots\times\bb{P}^{a_n}) = (\bb{P}^{a_1}\times\cdots\times H\times\cdots\times\bb{P}^{a_1})\cap\iota_\mathbf{a}(\bb{P}^\mathbf{a}).
\]

\begin{lemma}\label{polystella: lemma pullbackBuildingsetGFM}
    Let $\mathcal{G}_{\mathbb{P}^d,H}$ be the building set on $\mathbb{P}^\mathbf{d}$ defining $(\mathbb{P}^d)^{[n]}_{H}$ as in Equation~\eqref{ps: ks building set}, and let     
    \[
        \mathcal{G}_{\mathbf{a}} = \{ H_S\cap\mathbb{P}^{\mathbf{a}}\,\vert\ \,
        H_S \in \mathcal{G}_{\mathbb{P}^d,H}    
        \}
    \]
    be its restriction to $\bb{P}^\mathbf{a}$. 
    Then $\mathcal{G}_{\mathbf{a}}$ is a building set on $\mathbb{P}^{\mathbf{a}}$. 
\end{lemma}

\begin{proof}
    The building set $\mathcal{G}_{\mathbb{P}^d,H}$ is closed under intersection, so $\mathcal{G}_{\mathbf{a}}$, is as well. Therefore, it is enough to show that the intersection of its elements is clean.
       
    Since $H_{i+1}=\mathbb{P}^{a_{i+1}}\cap H_{i}$, we obtain that $H_{i+1} = \mathbb{P}^{a_{i+1}}\cap H$. Therefore, the subvariety $H_S\cap\mathbb{P}^{\mathbf{a}}$ is the product subvariety of $\mathbb{P}^{\mathbf{a}}$ obtained by taking $\mathbb{P}^{a_i}$ in the components labeled by $i\in [n]\setminus S$, and $H_i$ in those labeled by $i\in S$. The fact that the elements of $\mathcal{G}_{\mathbf{a}}$ intersect cleanly is a direct consequence of the fact that $H_S\cap\mathbb{P}^{\mathbf{a}}$ is a multilinear subspace. Indeed, by identifying $H_S\cap\mathbb{P}^{\mathbf{a}}$ with its embedded tangent space one obtains that 
    \[
        T_{H_{S_1}\cap\mathbb{P}^{\mathbf{a}},y}\cap T_{ H_{S_2}\cap\mathbb{P}^{\mathbf{a}},y} = 
        H_{S_1}\cap\mathbb{P}^{\mathbf{a}}\cap H_{S_2}\cap\mathbb{P}^{\mathbf{a}} 
        =H_{S_1}\cap H_{S_2}\cap\mathbb{P}^{\mathbf{a}}
        = T_{H_{S_1}\cap H_{S_2}\cap\mathbb{P}^{\mathbf{a}}, y}
    \]
    for all $y\in  H_{S_1}\cap H_{S_2}\cap\mathbb{P}^{\mathbf{a}}$ and $S_1,S_2\subseteq [n]$.
\end{proof} 

\begin{definition}
The generalized Fulton-MacPherson compactification of configurations of points in the flag 
$\mathbb{P}^{\mathbf{a}}$ lying outside of the subflag $H_n\subseteq\cdots\subseteq H_1$ is the wonderful compactification $\mathbb{P}^{[\mathbf{a}]}_{H}$ of $\mathbb{P}^\mathbf{a}$ with respect to the building set $\mathcal{G}_{\mathbf{a}}$ defined in Lemma~\ref{polystella: lemma pullbackBuildingsetGFM}.
\end{definition}

\begin{proposition}\label{polystella: main proposition}
   Given $\mathbf{a}=(a_1,\dots,a_n)$ and $\mathbf{d} = (a_1, \ldots, a_1)$. Then, there exists a unique closed immersion which makes the following diagram commute 
    \[
    \begin{tikzcd}
        \mathbb{P}^{[\mathbf{a}]}_{H} 
        \ar[r,hook]\ar[d] & 
          \mathbb{P}^{[\mathbf{d}]}_{H} 
        \ar[d]
        \\
        \mathbb{P}^{\mathbf{a}}\ar[r,hook, "\iota_\mathbf{a}"] & \mathbb{P}^{\mathbf{d}}.
    \end{tikzcd}
    \] 
\end{proposition}
\begin{proof}
Let $X = \PP^{\mathbf{d}}$, $Y = \PP^{\mathbf{a}}$, and let $\mathcal{G}=\mathcal{G}_{\mathbf{d}}$ and $\mathcal{H}=\mathcal{G}_{\mathbf{a}}$. 
Lemma~\ref{polystella: lemma pullbackBuildingsetGFM} shows that the intersections with $Y$ of the elements of the building set $\mathcal{G}$ on $X$ are precisely the elements of the building set $\mathcal{H}$ on $Y$, and therefore the conclusion follows by Lemma \ref{lemma.wonderful.subvariety}.  
\end{proof}
Finally, we remark that there is a ``universal" family $(\mathbb{P}^{[\mathbf{a}]}_{H})^+\to\mathbb{P}^{[\mathbf{a}]}_{H}$ parametrizing stable degenerations of $\mathbb{P}^{\mathbf{a}}$ with $n$ smooth labeled points away from $H$.  This family is defined to be the pullback of the universal family $(\mathbb{P}^{[\mathbf{d}]}_{H})^+\to\mathbb{P}^{[\mathbf{d}]}_{H}$ constructed by Kim and Sato, as described in Theorem~\ref{thm:GeneralizedFM}.

\subsection{The polystellahedral variety as a blow-up of a product of projective spaces}
\label{sec:PolyGeometry}

In this section we construct a collection of fans that interpolate between the fan of a given polystellahedral variety and the fan of a product of projective spaces. 
This interpolation is by explicit star subdivisions, and the cones in the intermediate fans are defined using the notion of \emph{compatible triples}, which we introduce as a refinement of the compatible pairs from \cite[Definition 1.2]{eur2024intersection}. 
This provides an explicit description of any given polystellahedral variety as an iterated blow-up of a product of projective spaces in Theorem~\ref{PS: construction as blow-up of Pa}.

\begin{notation}
    Let us denote the canonical basis of $\mathbb{R}^{A}$ by $\{ e_{i} : i \in A \}$ and for each $I \subseteq A$ let $e_I=\sum_{i \in I} e_i$.
    We let $\{ e^{*}_{i} : i \in A \}$ denote the associated dual basis of the dual vector space $\left(\mathbb{R}^{A}\right)^*$. 
We denote the classes of $e_i$ and $e_I$ in the quotient vector space $\mathbb{R}^{A}/\langle e_A \rangle$ by $\overline{e}_i$ and $\overline{e}_I$, for any $e_i$ and $e_I$. 
\end{notation}

%%%%%%%%%%%%%%%%%%%%%%%%%%%%%%%%%%%%%%%%%%%%%%%%%%%%%%%%%%%%%%

\begin{definition}[{\cite[Definition~1.2]{eur2024intersection}\label{def.compatible.pairs}}] A \emph{compatible pair} with respect to a surjective function between finite sets $\pi: A \rightarrow E$ is a pair $I \leq \mathcal{F}$ consisting of a subset $I \subseteq A$ and a chain 
$$
\mathcal{F}=\left\{F_1 \subsetneq F_2 \subsetneq F_3 \subsetneq \cdots \subsetneq F_{k} \subsetneq F_{k+1}=E\right\}
$$ 
of proper subsets of $E$ such that 
if $\pi^{-1}(S) \subseteq I$ for a subset $S \subseteq E$, then $S \subseteq F_1$. 
\end{definition}

\begin{remark}
    In the language of building sets, Larson and Eur \cite[Proof of Proposition~2.2]{eur2024intersection} show that the polystellahedron with caging $\pi:A\to E$ is the nestohedron coming from the building set
    \[
        \{\{i\}\,\vert\, i\in A\}\cup \{\{0\}\cup S\,\vert\, \emptyset\subseteq S\subseteq E\}.
    \]
\end{remark}

%%%%%%%%%%%%%%%%%%%%%%%%%%%%%%%%%%%%%%%%%%%%%%%%%%%%%%%%%%%%%%

\begin{definition} \label{def.compatible.triple} 
A \emph{compatible triple} with respect to a surjective function between finite sets $\pi: A \rightarrow E$ is a triple $(I,\mathcal{F},J)$ consisting of subsets $I \subseteq A$ and $J \subseteq E$, and a chain 
\[
\mathcal{F}=\left\{F_1 \subsetneq F_2 \subsetneq F_3 \subsetneq \cdots \subsetneq F_{k} \subsetneq F_{k+1} \subseteq E \setminus J \right\}
\]
of proper subsets of a set $F_{k+1} \subseteq E \setminus J$, satisfying the following compatibility condition. 
For each $S \subseteq E$ such that $\pi^{-1}(S) \subseteq I$:  
If $\mathcal{F} \neq \emptyset$ then $S \subseteq F_1$ and 
if $\mathcal{F} = \emptyset$ then  $S \cap J = \emptyset$.
\end{definition}

\begin{remark}
Definition~\ref{def.compatible.triple} builds upon the preexisting terminology of compatible pairs. 
Compatible triples provide the parameters needed to describe all the cones in the explicit blow-up construction of all polystellahedral varieties that we provide in this section. 
The conditions in Definition~\ref{def.compatible.triple} are interpreted in analogy to the case of compatible pairs, as we now explain. 
The set $F_{k+1}$ is implicitly a part of $\mathcal{F}$.  
Like in the case of compatible pairs, while the sets $F_l$ are allowed to be empty, 
the chain $\mathcal{F}$ is considered to be empty if there are no sets $F_l$ for $1 \leq l \leq k$.   
For any compatible triple $(I,\mathcal{F},J)$, we have that $S \cap J = \emptyset$ for any $S \subseteq E$ such that $\pi^{-1}(S) \subseteq I $. 
This holds by definition if $\mathcal{F} = \emptyset$ and it follows from the other assumptions in Definition~\ref{def.compatible.triple} when $\mathcal{F} \neq \emptyset$. 
Indeed, if $\mathcal{F} \neq \emptyset$ and $\pi^{-1}(S) \subseteq I $, we have 
$S \cap J   
\subseteq  F_1 \cap J 
\subseteq F_{k+1} \cap J = \emptyset$, so $S \cap J = \emptyset$.
 
\end{remark}

\begin{definition}       \label{definition.associated.cone.pairs}
Given a compatible pair $I\leq \mathcal{F}$ as in Definition~\ref{def.compatible.pairs}, \cite[Definition 1.2]{eur2024intersection} defines its \emph{associated cone} $\sigma_{I \leq \mathcal{F}}$ as the cone in the space $\mathbb{R}^A$ given by
\begin{align*}
    \sigma_{I \leq \mathcal{F}} 
    &= \operatorname{cone}(e_i \,\vert\, i \in I) 
        + \operatorname{cone}\left(-e_{A \setminus \pi^{-1}(F_1)}, \ldots, -e_{A \setminus \pi^{-1}(F_k)}\right).
\end{align*}
\end{definition}

\begin{definition}       \label{definition.associated.cone.triple}
Given a compatible triple $(I,\mathcal{F},J)$ as in Definition~\ref{def.compatible.triple}, we define its \emph{associated cone} $C_{I,\mathcal{F},J}$ as the cone in the space $\mathbb{R}^A$ given by
\begin{align*}
    C_{I,\mathcal{F},J} 
    &= \operatorname{cone}(e_i \,\vert\, i \in I) 
        + \operatorname{cone}\left(-e_{A \setminus \pi^{-1}(F_1)}, \ldots, -e_{A \setminus \pi^{-1}(F_k)}\right)
    + \operatorname{cone}(-e_{\pi^{-1}(j)} \,\vert\, j \in J).
\end{align*}
\end{definition}

\begin{remark}
The cones $\sigma_{I \leq \mathcal{F}} \subseteq \mathbb{R}^{A}$ associated to the compatible pairs in \cite[Definition 1.2]{eur2024intersection} are precisely the cones in Definition~\ref{definition.associated.cone.triple}
associated to the compatible triples satisfying $F_{k+1}=E$.   
Indeed, for any compatible pair $I \leq \mathcal{F}$ we have that $F_{k+1} = E$; there is a unique possible extension to a compatible triple $(I,\mathcal{F},J)$, namely with $J = \emptyset$; and $\sigma_{I \leq \mathcal{F}} = C_{I,\emptyset,\mathcal{F}}$. 
Conversely, given any compatible triple $(I,\mathcal{F},J)$ with $F_{k+1}= E$, then $I \leq \mathcal{F}$ is a compatible pair, $J = \emptyset$, and $C_{I,\mathcal{F},J} = C_{I,\emptyset,\mathcal{F}} = \sigma_{I \leq \mathcal{F}}$.
\end{remark}

\begin{lemma} \label{lemma.cones.from.compatible.triples.are.smooth}
Let $(I,\mathcal{F},J)$ be a compatible triple with respect to $\pi:A \rightarrow E$.  
Then, the vectors 
\[
\{e_i : i \in I \}   \cup \
\left\{-e_{A \setminus \pi^{-1}(F_1)}, \ldots, -e_{A \setminus \pi^{-1}(F_k)}\right\}   \cup \ 
\{-e_{\pi^{-1}(j)} : j \in J \} 
\] 
are the minimal set of integral generators of the cone $C_{I,\mathcal{F},J}$. 
Moreover, the cone $C_{I,\mathcal{F},J}$ is smooth. 
\end{lemma}

\begin{proof}
    Recall the vector space basis $\{ e^{*}_{i} : i \in A \}$ of $\left(\mathbb{R}^{A}\right)^*$ dual to the basis $\{ e_{i} : i \in A \}$ of $\mathbb{R}^{A}$.
It is enough to exhibit lattice vectors in $\left(\mathbb{R}^{A}\right)^*$, $u_i$ for each $i \in I$; $v_l$ for each $1 \leq l \leq k$; and $w_j$ for each $j \in J$, such that for each $i_0 \in i$, $1 \leq l_0 \leq k$, and $j_0 \in J$, we have  
\begin{equation}  \label{equation.desired.vectors}
\hspace*{-0.9cm} % Adjust this value as needed
\begin{minipage}{0.48\linewidth}
\[
u_{i_0}(v) = 
\begin{cases} 
1 & \text{if } v=e_{i_0}, \\
0 & \text{if } v=e_{i},\ i \in I \setminus \{i_0\}, \\
0 & \text{if } v=-e_{A \setminus \pi^{-1}(F_l)},\ 1 \leq l \leq k, \\
0 & \text{if } v=-e_{\pi^{-1}(j)}, \ j \in J.
\end{cases}
\]
\end{minipage}
\begin{minipage}{0.48\linewidth}
\[
\begin{aligned}
v_{l_0}(v) &= 
\begin{cases} 
1 & \text{if } v=-e_{A \setminus \pi^{-1}(F_{l_0})}, \\
0 & \text{if } v=-e_{A \setminus \pi^{-1}(F_l)}, \ 1 \leq l \leq k, \ l \neq l_0, \\
0 & \text{if } v=-e_{\pi^{-1}(j)}, \ j \in J.
\end{cases}  \\
w_{j_0}(v) &= 
\begin{cases} 
1 & \text{if } v=-e_{\pi^{-1}(j_0)}, \ j \in J,  \\
0 & \text{if } v=-e_{\pi^{-1}(j)}, \ j \in J, \ j \neq j_0.
\end{cases}
\end{aligned}
\]
\end{minipage}
\end{equation}
Let us first assume $\mathcal{F} \neq \emptyset$.
For each $i \in I$ with $\pi(i) \notin F_1$, choose an element $i' \in A \setminus I$ with $\pi(i)=\pi(i')$.    
For each $1 \leq l \leq k+1$, choose an element $i_{l} \in A$ such that $\pi(e_{i_l}) \in F_l \setminus F_{l-1}$, where $F_0:=\empty$.  
For each $j \in J$, choose an element $i_j \in A$ such that $\pi(i_j) =j$ (notice we could also assume $i_j \notin I$). 
Now, for each $i \in I$, each $1 \leq l \leq k$, and each $j \in J$, let us define 
\begin{equation} \label{equation.dual.vectors}
\begin{aligned}
u_i &= 
\begin{cases} 
e^{*}_{i} & \text{if } \pi(i) \in F_1, \\
e^{*}_{i} - e^{*}_{i'} & \text{if } \pi(i) \notin F_1, 
\end{cases} \qquad &
v_l &= e^{*}_{i_l} - e^{*}_{i_{l+1}}, \qquad &
w_j &= - e_{i_j}^*.
\end{aligned}
\end{equation}
It is straightforward to verify that the vectors in Equation~\eqref{equation.dual.vectors} satisfy the desired conditions in Equation~\eqref{equation.desired.vectors}.
Now, let us consider the remaining case, namely that $\mathcal{F} = \emptyset$.  
For each $i \in A$ such that $i \in I \cap \pi^{-1}(J)$ we choose $i' \in A$ such that $\pi(i')=\pi(i)$ and $i' \notin I$, and for each $j \in J$, we choose an element $i_j \in A$ such that $\pi(i_j) =j$ (notice we could also assume $i_j \notin I$).  
For each $i \in I$ and each $j \in J$, let us define 
\begin{equation} \label{equation.dual.vectors.empty.case}
\begin{aligned}
u_i &= 
\begin{cases} 
e^{*}_{i} & \text{if } \pi(i) \notin J, \\
e^{*}_{i} - e^{*}_{i'} & \text{if } \pi(i) \in J, 
\end{cases} \qquad & 
w_j &= - e_{i_j}^*.
\end{aligned}
\end{equation}
It is straightforward to verify that the vectors in Equation~\eqref{equation.dual.vectors.empty.case} satisfy the desired conditions in Equation~\eqref{equation.desired.vectors} relevant to the case $\mathcal{F} = \emptyset$, and this completes the proof.
\end{proof}

%%%%%%%%%%%%%%%%%%%%%%%%%%%%%%%%%%%%%%%%%%%%%%%%%%%%%%%%%%%%%%%%%%%%%%%%%%

\begin{remark}  \label{faces.of.compatible.triples} 
    In view of Lemma~\ref{lemma.cones.from.compatible.triples.are.smooth} for each compatible triple $(I,\mathcal{F},J)$ the faces of the cone $C_{I,\mathcal{F},J}$ are precisely the cones generated by subsets of the minimal set of generators, and hence each such face is itself the cone $C_{I',\mathcal{F}',J'}$ associated to some compatible triple $(I',\mathcal{F}',J')$. 
\end{remark}

\begin{remark}     \label{form.maximal.cones}
    The cones that are maximal with respect to containment among the collection of all cones associated to compatible triples have the following two forms depending on whether the flag $\mathcal{F}$ is empty.   
    If the flag $\mathcal{F}$ is empty, there exists a section $\theta: J \rightarrow A$ of $\pi$ over $J$ (i.e., $\pi \circ \theta = \operatorname{Id}_{J}$), such that $I = A \setminus \theta(J)$. Hence, in this case $C_{I,\mathcal{F},J}$ is determined by the function $\theta: J \rightarrow A$ and has the form 
    \[
    C_{I,\mathcal{F},J} = C_{\theta} = \operatorname{cone}(e_i \,\vert\, i \in A \setminus \theta(J)) + \operatorname{cone}(-e_{\pi^{-1}(j)} \,\vert\, j \in J).  
    \]
    If the flag $\mathcal{F}$ is not empty, there exists a section $\theta: (E \setminus F_1) \rightarrow A$ of $\pi$ over $E \setminus F_1$ (i.e., $\pi \circ \theta = \operatorname{Id}_{E \setminus F_1}$), such that $I = A \setminus \theta(E \setminus F_1)$.
    Moreover, in this case $F_{k+1} = E \setminus J$ and $\mathcal{F}$ is a flag of proper subsets of $E \setminus J$, with $F_1$ possibly empty, which is full between $F_1$ and $F_{k+1}$, that is, $|F_{l}| = |F_{1}| + l -1$ for each $1 \leq l \leq k+1$. Hence, in this case $C_{I,\mathcal{F},J}$ has the form 
\[
    C_{I,\mathcal{F},J} 
    = \operatorname{cone}(e_i \,\vert\, i \in A \setminus \theta(E \setminus F_1)) 
        + \operatorname{cone}\left(-e_{A \setminus \pi^{-1}(F_1)}, \ldots, -e_{A \setminus \pi^{-1}(F_k)}\right)
    + \operatorname{cone}(-e_{\pi^{-1}(j)} \,\vert\, j \in J).
 \]   
\end{remark}

\begin{remark}
All cones that are maximal with respect to containment among the collection of all compatible triples are themselves maximal cones, that is, they are full-dimensional in $\mathbb{R}^A$.     
We can see this by verifying the cardinality of the minimal generating set is always $|A|$. 
We use the notation and description of the maximal cones in Remark~\ref{form.maximal.cones} and we use Lemma~\ref{lemma.cones.from.compatible.triples.are.smooth} to obtain the minimal generating sets.  
When $\mathcal{F}$ is empty we have $\operatorname{dim} C_{I,\mathcal{F},J} = |A| - |\theta(J)| + |J| = |A|$. 
Similarly, when $\mathcal{F}$ is not empty we have 
\[
\operatorname{dim} C_{I,\mathcal{F},J} = |A \setminus \theta(E \setminus F_1)|+|k|+|J| = (|A| - |E| + |F_1|) + (|F_{k+1}|- |F_1|) +(|E| - |F_{k+1}|) = |A|.
\] 
\end{remark}

%%%%%%%%%%%%%%%%%%%%%%%%%%%%%%%%%%%%%%%%%%%%%%%%%%%%%%%%%%%%%%%%%%%%%%%%%%

\begin{definition} \label{definition.associated.fans.triple}
    In the setting of Definition~\ref{definition.associated.cone.triple}, for each $0 \leq s \leq |E|$, let
        \[
    \Delta_{\pi,s} = \left\{ \ C_{I,\mathcal{F},J} \ : \  \textnormal{($I, \mathcal{F},J)$ is a compatible triple with respect to $\pi$ and $|F_{k+1}|=s$} \  \right\}.  
    \]
\end{definition}

\begin{remark}  \label{remark.triples.with.empty.fan}
For a compatible triple $(I,\mathcal{F},J)$ with $\mathcal{F} $ empty, the condition that $C_{I,\mathcal{F},J}$ is in $\Delta_{\pi,s}$ is equivalent to requiring that $|J| \leq |E| - s$ and that $I$ does not contain $\pi^{-1}(j)$ for any $j \in J$.  
\end{remark}

\begin{remark}
    Each collection $\Delta_{\pi,s}$ is closed under taking faces, as we can see from Remark~\ref{faces.of.compatible.triples}. We will prove in Lemma~\ref{technical.lemma.blow-up.triples} that they are actually smooth projective fans.  
\end{remark}

\begin{example}  \label{example.fan.product.projective.spaces}
The collection of cones $\Delta_{\pi,0}$ is the standard fan in $\mathbb{R}^{A} = \mathbb{R}^{a_1} \times \mathbb{R}^{a_2} \times \cdots \times \mathbb{R}^{a_r}$
defining the product of projective spaces $\mathbb{P}^{a_1} \times \mathbb{P}^{a_2} \times \cdots \times \mathbb{P}^{a_r}$ as a toric variety. 
To see this, consider a cone $C_{I,\mathcal{F},J} \in \Delta_{\pi,0}$ where $(I,\mathcal{F},J)$ is a compatible triple. 
Since $|F_{k+1}|=s=0$, then the flag $\mathcal{F}$ is empty. 
Moreover, $I \subseteq A$ and $J \subseteq E$ can be any subsets such that $I$ contains no fibers of $\pi$ over elements of $J$, by the compatibility condition. 
Hence, the cone $C_{I,\mathcal{F},J}$ has the form 
\[
\begin{aligned}
    C_{I,\mathcal{F},J}
&=\operatorname{cone}(e_i \,\vert\, i \in I)+\operatorname{cone}(-e_{\pi^{-1}(j)} \,\vert\, j \in J) \\
&= \sum_{j \in E} \operatorname{cone}(e_i \,\vert\, i \in I \cap \pi^{-1}(j))+ \sum_{j \in J} \operatorname{cone}(-e_{\pi^{-1}(j)}) = \sum_{j \in E}  \sigma_j   
\end{aligned}
\]
where for each $j \in E$, 
\[
\sigma_j =
\begin{cases} 
\operatorname{cone}(e_i \,\vert\, i \in I \cap \pi^{-1}(j)) & \text{if } j \notin J, \\
\operatorname{cone}(e_i \,\vert\, i \in I \cap \pi^{-1}(j)) + \operatorname{cone}(-e_{\pi^{-1}(j)}) & \text{if } j \in J.  
\end{cases} 
\]
By the compatibility condition, when $j \in J$ we have that $I \cap \pi^{-1}(j) \subsetneq \pi^{-1}(j)$. 
Therefore, for each $j \in E$ the possible cones $\sigma_j$ are precisely all the cones generated by proper subsets of the set of vectors $\{ e_i : i \in \pi^{-1}(j) \} \cup \{ -e_{\pi^{-1}(j)} \}$.
Since for each $j \in E$ the cones in $\mathbb{R}^{a_j}$ generated by all proper subsets of set of vectors $\{ e_i : i \in \pi^{-1}(j) \} \cup \{ -e_{\pi^{-1}(j)} \}$ form the standard fan of $\mathbb{P}^{a_j}$, we see that $\Delta_{\pi,0}$ is the fan of $\mathbb{P}^{a_1} \times \mathbb{P}^{a_2} \times \cdots \times \mathbb{P}^{a_r}$ as we claimed. 
\end{example}

\begin{example}       \label{example.fan.polystellahedron}
The collection of cones $\Delta_{\pi,|E|}$ is the fan of the polystellahedral variety $PS_\pi$.
Indeed, the compatible triples $(I,\mathcal{F},J)$ such that $C_{I,\mathcal{F},J} \in \Delta_{\pi,|E|}$ are those with $F_{k+1} = E$.
Therefore, given one such compatible triple 
the set $J$ is empty; 
$$
\mathcal{F}=\left\{ F_1 \subsetneq F_2 \subsetneq F_3 \subsetneq \cdots \subsetneq F_{k} \subsetneq F_{k+1}=E\right\}
$$ is a possibly empty flag of proper subsets of $F_{k+1}=E$; and $I \subseteq A$ is any subset such that if $\pi^{-1}(j) \subseteq I$ for some $j \in E$, then $j \in F_1$.  
Hence, $I \leq \mathcal{F}$ is a compatible pair as in Definition~\ref{def.compatible.pairs}.    
Moreover, all compatible pairs arise in this way from a unique compatible triple with $F_{k+1} = E$,  
since for each compatible pair $I \leq \mathcal{F}$, we have that $(I,\mathcal{F},\emptyset)$ is a compatible triple with $F_{k+1} = E$. 
It is clear that $C_{I,\mathcal{F},\emptyset} = \sigma_{I \leq \mathcal{F}}$ for any compatible triple $(I,\mathcal{F},J)$ with $J = \emptyset$.  
Therefore, $\Delta_{\pi,|E|}$ is precisely the fan of the polystellahedral variety $PS_\pi$ presented in \cite[Section 2.1 and Proposition 2.3]{eur2024intersection}.
\end{example}

\begin{definition} \label{definition.cones.CJ}
    Given $\pi : A \rightarrow E$, for each $J \subseteq E$ we define the cone $C_J$ as the cone associated to the compatible triple $(\emptyset,\emptyset,J)$, that is, 
    \[
    C_J := C_{\emptyset,\emptyset,J} = \operatorname{cone}( - e_{\pi^{-1}(j)} \,\vert\, j \in J).
    \] 
\end{definition}

%%%%%%%%%%%%%%%%%%%%%%%%%%%%%%%%%%%%%%%%%%%%%%%%%%%%%%%%%%%%%%%%%%%%%%%%%%%%%%%%%

\begin{remark}
    By Remark~\ref{remark.triples.with.empty.fan}, for each $J \subseteq E$  the cone $C_{J} = C_{\emptyset,\emptyset,J}$ is in $\Delta_{\pi,s}$  for each $s \leq |E|-|J|$. 
\end{remark}

\begin{lemma}  \label{technical.lemma.cone.containments}
        Let $0 \leq s \leq |E|-1$, let $(I,\mathcal{F},J)$ be a compatible triple such that $C_{I,\mathcal{F},J}$ is in $\Delta_{\pi,s}$,  
        and let $J_0 \subseteq E$ be any subset. 
        If $C_{J_0} \subseteq C_{I,\mathcal{F},J}$, then $J_0 \subseteq J$. In particular, 
        \begin{enumerate}
            \item \label{technical.lemma.cone.containments.part.1} If $|J| < |E|-s$, then $C_{I,\mathcal{F},J}$ contains no cone $C_{J_0}$ such that $J_0 = |E|-s$. 
            \item \label{technical.lemma.cone.containments.part.2} If $|J| = |E|-s$, then $C_{I,\mathcal{F},J}$ contains exactly one cone $C_{J_0}$ such that $J_0 = |E|-s$, namely $C_J$. 
        \end{enumerate} 
\end{lemma}

\begin{proof}
Let us assume that $C_{J_0} \subseteq C_{I,\mathcal{F},J}$ and let us show that $J_0 \subseteq J$. 
We have that $C_{J_0}=\operatorname{cone}(- e_{\pi^{-1}(j)} \,\vert\, j \in J_0) \subseteq C_{I,\mathcal{F},J}$, and then 
$- e_{\pi^{-1}(j)} \in C_{I,\mathcal{F},J}$ for each $j \in J_0$. 
We will now show that $- e_{\pi^{-1}(j)} \notin C_{I,\mathcal{F},J}$ for each $j \in E \setminus J$, from which $J_0 \subseteq J$ clearly follows. 
Let us then fix $j \in E \setminus J$. 
 
If the flag $\mathcal{F}$ is empty, we define $I' = I \cup \pi^{-1}(j)$.
We have that $(I',\mathcal{F},J)$ is a compatible triple and $C_{J_0} \subseteq C_{I,\mathcal{F},J} \subseteq C_{I',\mathcal{F},J}$. 
We see that $e_{\pi^{-1}(j)} \in \operatorname{cone}( e_i  \,\vert\, i \in I')  \subseteq C_{I',\mathcal{F},J} $. 
But then $- e_{\pi^{-1}(j)} \notin C_{I',\mathcal{F},J}$ because the cone $C_{I',\mathcal{F},J}$ is pointed, since it is smooth.  
Hence, $- e_{\pi^{-1}(j)} \notin C_{I,\mathcal{F},J}$, as desired. 

We can now suppose that the flag $\mathcal{F}=\left\{ F_1 \subsetneq F_2 \subsetneq \cdots \subsetneq F_{k} \subsetneq F_{k+1}=E \setminus J \right\}$ is not empty. 
We consider the cases $j \in F_1$ and $j \notin F_1$ separately.

If $j \in F_1$ we proceed as before. That is,  we define $I' = I \cup \pi^{-1}(j)$ and then $(I',\mathcal{F},J)$ is a compatible triple such that $C_{J_0} \subseteq C_{I,\mathcal{F},J} \subseteq C_{I',\mathcal{F},J}$. 
Then, $- e_{\pi^{-1}(j)} \notin C_{I',\mathcal{F},J}$ because $C_{I',\mathcal{F},J}$ is pointed, and then $- e_{\pi^{-1}(j)} \notin C_{I,\mathcal{F},J}$, as desired.

If $j \notin F_1$, then $j \in E \setminus (F_1 \cup J)$. 
 By the compatibility condition, for each $j' \in E \setminus F_1$ the fiber $\pi^{-1}(j')$ is not contained in $I$.
Then, for each $j' \in E \setminus F_1$ we can choose $\theta(j') \in A$ such that $\theta(j') \in \pi^{-1}(j')$ and $\theta(j') \notin I$. 
We consider the basis $\{e^*_i \, | \, i \in A \}$ of $(\mathbb{R}^A)^*$ which is dual to the basis  $\{e_i \, | \, i \in A \}$ of $\mathbb{R}^A$.
We notice that $J$ is not empty since $|J| = |E| - s \neq 0$ and that $J \subseteq E \setminus F_1$.  
Then we can define the nonzero functional $\psi:= - \sum_{j' \in J} e^*_{\theta(j')} \in (\mathbb{R}^A)^*$.  
 For each $j' \in J$ we have $\pi^{-1}(j') \subseteq \pi^{-1}(J) \subseteq A \setminus \pi^{-1}(F_l)$ for each $1 \leq l \leq k$,  
since $J \subseteq E \setminus F_{k+1} \subseteq E \setminus F_l$ for each $1 \leq l \leq k$. 
By construction, the values of $\psi$ on the minimal integral generating set of the cone $C_{I,\mathcal{F},J}$ satisfy 
$\psi(e_i) = 0$ for all $i\in I$, $\psi(-e_{A \setminus \pi^{-1}(F_l)}) > 0$ for all $1\leq l \leq k$, and $\psi(-e_{\pi^{-1}(j')}) > 0$ for all $j' \in J$.
Then $\psi$ induces a supporting hyperplane of the cone $C_{I,\mathcal{F},J}$ and it defines the face $C_{I,\mathcal{F},J} \cap \{ u \in \mathbb{R}^A \, | \, \psi(u)=0 \} = \operatorname{cone}(e_i \,\vert\, i \in I)$ of $C_{I,\mathcal{F},J}$. 
Since $j \notin J$, we see that $\psi(-e_{\pi^{-1}(j)}) = 0$. 
Then, if we assume by contradiction that $-e_{\pi^{-1}(j)} \in C_{I,\mathcal{F},J}$,  
we would deduce that $-e_{\pi^{-1}(j)} \in \operatorname{cone}(e_i \,\vert\, i \in I)$. 
But this is a contradiction because the functional $\gamma = e^*_{\theta(j)}$ takes the value zero on all elements in $\{ e_i \, \vert \, i \in I\}$ and the negative value $-1$ on $-e_{\pi^{-1}(j)}$.
This completes the proof that $C_{J_0} \subseteq C_{I,\mathcal{F},J}$ implies $J_0 \subseteq J$, and from this (\ref{technical.lemma.cone.containments.part.1}) 
follows since $J_0 \subseteq J$ is not possible in that case
and (\ref{technical.lemma.cone.containments.part.2}) follows since only $J_0=J$ is possible in that case. 
%This completes the proof.   
\end{proof}

%%%%%%%%%%%%%%%%%%%%%%%%%%%%%%%%%%%%%%%%%%%%%%%%%%%%%%%%%%%%%%%%%%%%%%%%%%%%%%%%%

\begin{lemma}     \label{technical.lemma.blow-up.triples} 
    For each $0 \leq s \leq |E|$, the collection of cones $\Delta_{\pi,s}$ forms a smooth projective fan. 
    Moreover, for each $0 \leq s \leq |E|-1$, the fan $\Delta_{\pi,s+1}$ is the iterated star subdivision of the fan $\Delta_{\pi,s}$ along the cones $C_J=\operatorname{cone}(e_j \,\vert\, j \in \pi^{-1}(J))$, for all $J \subseteq E$ with $|J|=|E|-s$, in any order. 
\end{lemma}

\begin{proof} 
The collection of cones $\Delta_{\pi,0}$ forms a smooth projective fan by Example~\ref{example.fan.product.projective.spaces}. 
Let $0 \leq s \leq |E|-1$ and assume that $\Delta_{\pi,s}$ is a smooth projective fan. 
Let $\Sigma_{s+1}$ be the smooth projective fan obtained by the iterated star subdivision of the fan $\Delta_{\pi,s}$ along the cones $C_J$, for all $J \subseteq E$ with $|J|=|E|-s$, in any order.
Since there is no cone in $\Delta_{\pi,s}$ containing two distinct cones $C_J$ by Lemma~\ref{technical.lemma.cone.containments}, then the resulting fan $\Sigma_{s+1}$ is independent of the order of the star subdivisions along the $C_J$. 
Geometrically these star subdivisions correspond to blow-ups along disjoint loci.
The claims in the statement all follow if we show that $\Sigma_{s+1}=\Delta_{\pi,s+1}$.
Since the collections of cones $\Sigma_{s+1}$ and $\Delta_{\pi,s+1}$ are both closed under taking faces, it suffices to show that any cone that is maximal within either of the collections $\Sigma_{s+1}$ or $\Delta_{\pi,s+1}$ also belongs to the other collection. 

Let $\tau$ be a maximal cone in $\Sigma_{s+1}$. 
Hence, there exists a maximal cone $\tau_0$ in $\Delta_{\pi,s}$ such that $\tau$ is one of the cones obtained from $\tau_0$ by the star subdivisions along all cones $C_J$, for all $J \subseteq E$ with $|J|=|E|-s$. 
By Lemma~\ref{technical.lemma.cone.containments} we know that $\tau_0$ is subdivided by at most one of these star subdivisions. 
Since $\tau_0 \in \Delta_{\pi,s}$, there exists a compatible triple $(I_0, \mathcal{F}_0,J_0)$ such that $\tau_0=C_{I_0,\mathcal{F}_0,J_0}$ and $|F_{k=1}| = s$, 
where $\mathcal{F}_0$ has the form $\mathcal{F}_0=\left\{ F_1 \subsetneq F_2 \subsetneq \cdots \subsetneq F_{k} \subsetneq F_{k+1} \subseteq E \setminus J_0 \right\}$ (here $\mathcal{F}_0$ may be empty, which by our conventions means $k=0$ and $\mathcal{F}_0=\left\{ F_{k+1} \subseteq E \setminus J_0 \right\}$).  
Whether or not the flag $\mathcal{F}_0$ of proper subsets of $F_{k+1}$ is empty, we know that $J_0 \subseteq E \setminus F_{k+1}$, so $|J_0| \leq |F_{k+1}| = |E| - s$.  
We consider the cases $|J_0| < |E|-s$ and $|J_0| = |E|-s$.

If $|J_0| < |E|-s$ then the cone $\tau_0=C_{I_0,\mathcal{F}_0,J_0}$ does not contain any of the cones $C_J$ with $|J_0| = |E|-s$, and  therefore $\tau = \tau_0$. 
We must then show that $\tau = \tau_0$ is in $\Delta_{\pi,s+1}$.
Since $J_0 \subsetneq E \setminus F_{k+1}$ we can choose and element $j_0 \in E \setminus F_{k+1}$ such that $j_0 \notin J_0$. 
We define a new flag $\mathcal{F}'_0$ by $\mathcal{F}'_0=\left\{ F'_1 \subsetneq F'_2 \subsetneq \cdots \subsetneq F'_{k} \subsetneq F'_{k+1}=E \setminus J_0 \right\}$ where 
$F'_l := F'_l$ for each $1 \leq l \leq k$ and $F'_{k+1}:=F_{k+1} \cup \{j_0 \}$. 
It follows that $(I_0, \mathcal{F}'_0,J_0)$ is a compatible triple and since $|F'_{k+1}| = |F_{k+1}| +  1 =  s +1$, then $C_{I_0, \mathcal{F}'_0,J_0} \in \Delta_{\pi,s+1}$.
By construction $C_{I_0, \mathcal{F}'_0,J_0} = C_{I_0, \mathcal{F}_0,J_0}$,
and hence $\tau = \tau_0 = C_{I_0, \mathcal{F}'_0,J_0} \in \Delta_{\pi,s+1}$, as desired.

If $|J_0| = |E|-s$, then $J_0 = E \setminus F_{k+1}$, and by Lemma~\ref{technical.lemma.cone.containments}(\ref{technical.lemma.cone.containments.part.2}) the cone $\tau_0=C_{I_0,\mathcal{F}_0,J_0}$ contains exactly one of the cones $C_J$ with $|J| = |E|-s$, namely for $J=J_0$. 
We deduce that $\tau$ is one of the $|J_0|$ maximal cones that arise when we perform the star subdivision of the cone $C_{I_0,\mathcal{F}_0,J_0}$ along $C_{J_0}$. 
The minimal integral generators of $C_{I_0,\mathcal{F}_0,J_0}$ are 
$\{e_i: i \in I_0 \} \cup \{-e_{A \setminus \pi^{-1}(F_1)},\ldots,-e_{A \setminus \pi^{-1}(F_k)} \} \cup \{-e_{\pi^{-1}(j)} : j \in J_0 \}$ and $C_{F_1}$ is minimally generated by $\{-e_{\pi^{-1}(j)} : j \in J_0 \}$.  
Then there exists $j_0 \in J_0$ such that $\tau$ is generated by the set of vectors obtained by taking the minimal generators of $C_{I_0, \mathcal{F}_0,J_0}$, and removing the vector $-e_{\pi^{-1}(j_0)}$ and including the vector 
\[
\sum_{j \in J_0} -e_{\pi^{-1}(j)} 
= \sum_{j \in E \setminus F_{k+1}} -e_{\pi^{-1}(j)} 
= -\sum_{i \in \pi^{-1}(E \setminus F_{k+1})} e_i 
= -\sum_{i \in A \setminus \pi^{-1}(F_{k+1})} e_i 
= -e_{A \setminus \pi^{-1}(F_{k+1})}. 
\]
Hence $\tau$ is generated by $\{e_i: i \in I_0 \} \cup \{-e_{A \setminus \pi^{-1}(F_1)},\ldots,-e_{A \setminus \pi^{-1}(F_{k+1})} \} \cup \{-e_{\pi^{-1}(j)} : j \in J_0 \setminus \{j_0\} \}$. 

Let us define $I_1:=I_0$, $J_1:=J_0 \setminus \{j_0\}$, and 
$\mathcal{F}_1:=\left\{ F'_1 \subsetneq F'_2 \subsetneq \cdots \subsetneq F'_{k+1} \subsetneq F'_{k+2} \subseteq E \setminus J_1 \right\}$ 
where $F'_{l}:=F_{l}$ for each $1 \leq l \leq k+1$ and $F'_{k+2} = F_{k+1} \cup \{ j_0 \}$. 
Let us show that $(I_1,\mathcal{F}_1,J_1)$ is a compatible triple. 
Since $(I_1,\mathcal{F}_1,J_1)$ has the form required in Definition~\ref{def.compatible.triple}, we must simply verify the compatibility condition, and we do this considering the cases $k=0$ and $k>0$ separately. Note that the flag $\mathcal{F}_1$ is not empty, so we only need to show that any fiber of $\pi$ contained in $I_1=I_0$ is contained in $F'_1=F_1$.   
If $k = 0$, since $(I_0,\mathcal{F}_0,J_0)$ is a compatible triple then any fiber of $\pi$ contained in $I_1=I_0$ is contained in $E \setminus J_0 = F_{k+1}=F_1=F'_1$, as needed.  
If $k > 0$, since $I_1=I_0$ and $F'_1=F_1$, the compatibility condition for $(I_1,\mathcal{F}_1,J_1)$ is the same as that for $(I_0,\mathcal{F}_0,J_0)$, which holds since the latter is a compatible triple. 
Therefore, $(I_1,\mathcal{F}_1,J_1)$ is a compatible triple. 
The minimal integral generators set of the cone $C_{I_1,\mathcal{F}_1,J_1}$ are exactly the generators of $\tau$ that we described before, 
and therefore $C_{I_1,\mathcal{F}_1,J_1} = \tau$.
Since $|F'_{k+2}| = |F_{k+1}|+1=s+1$, then $\tau = C_{I_1, \mathcal{F}_1,J_1} \in \Delta_{\pi,s+1}$, as desired.

Now, let $\tau$ be a maximal cone in $\Delta_{\pi,s+1}$. 
Then there exists a compatible triple $(I,\mathcal{F},J)$ of the form 
$I \subseteq A$, $\mathcal{F}=\left\{ F_1 \subsetneq F_2 \subsetneq F_3 \subsetneq \cdots \subsetneq F_{k} \subsetneq F_{k+1} \subseteq E \setminus J \right\}$ for some $k \geq 0$, such that $|F_{k+1}|=s+1$ and 
\[
\tau = 
C_{I,\mathcal{F},J} 
= 
\operatorname{cone}(e_i \,\vert\, i \in I) +
\operatorname{cone}\left(-e_{A \setminus \pi^{-1}(F_1)}, \ldots, -e_{A \setminus \pi^{-1}(F_k)}\right) + 
\operatorname{cone}(-e_{\pi^{-1}(j)} \,\vert\, j \in J).
\]      
We show that $\tau$ is in $\Sigma_{s+1}$ considering the cases $k=0$ and $k>0$.  

In the case $k=0$, we have $|F_{k+1}| = |F_{1}| = s+1 > 0$, then we can choose $j_0 \in F_{k+1} = F_1$.  
We define $I_1:=I$, $J_1 := J$, and $\mathcal{F}_1 := \{F'_{1} \subseteq E \setminus J_0   \}$
where $F'_{1}:=F_1 \setminus \{j_0\}$. 
Let us show that $(I_1,\mathcal{F}_1,J_1)$ is a compatible triple. 
We have that $(I_1,\mathcal{F}_1,J_1)$ has the form in Definition~\ref{def.compatible.triple} so it is enough to verify the compatibility condition in the case that the flag is empty. 
Since $I_1:=I$, $J_1 := J$ the compatibility condition for $(I_1,\mathcal{F}_1,J_1)$ states the same as that for $(I,\mathcal{F},J)$, which holds as the latter is a compatible triple. 
Then, $(I_1,\mathcal{F}_1,J_1)$ is a compatible triple as well.  
Since $|F'_{k+1}| = |F_{k+1}| -1 = (s+1) -1 =s$, then $C_{I_1,\mathcal{F}_1,J_1}$ is a cone in $\Delta_{\pi,s}$.
From the definition of the cone associated to a compatible triple, we have  $C_{I_1,\mathcal{F}_1,J_1} = C_{I,\mathcal{F},J}$. 
We have that $|J_1| = |J| \leq |E| - |F_{k+1}| =  |E| - (s+1) < |E| -s$. 
Then, by Lemma~\ref{technical.lemma.cone.containments}(\ref{technical.lemma.cone.containments.part.1}) the cone $C_{I_1,\mathcal{F}_1,J_1}$ does not contain any of the cones $C_{J_0}$ with $J_0 \subseteq E$ such that $|J_0|=|E| -s$. 
Therefore the cone $C_{I_1,\mathcal{F}_1,J_1}$ in $\Delta_{\pi,s}$ is not modified by the star subdivisions that produce $\Sigma_{s+1}$ out of $\Delta_{\pi,s}$. 
Hence, $C_{I,\mathcal{F},J} = C_{I_1,\mathcal{F}_1,J_1}$ is in $\Sigma_{s+1}$, as desired. 

Let us assume now that $k>0$. 
Since $\tau = C_{I,\mathcal{F},J}$ is maximal in $\Delta_{\pi,s+1}$ we deduce that $J = E \setminus F_{k+1}$; that there is a function $\theta: E\setminus F_1 \rightarrow A$ such that $\pi \circ \theta$ is the identity over $E\setminus F_1$ and $I = A \setminus \theta(E\setminus F_1)$; and that $\mathcal{F}$ is such that $|F_{l}|  = |F_{1}|+(l-1)$ for each $1 \leq l \leq k+1$. 
Indeed, if one of these conditions did not hold, one could use it to enlarge $I$, $\mathcal{F}$ or $J$ and obtain a new compatible triple that yields a cone in $\Delta_{\pi,s+1}$ strictly containing $C_{I,\mathcal{F},J}$, but this is not possible by the maximality assumption. 
We have that $|F_{k+1} \setminus F_{k}|=1$, and we denote by $j_0 \in E$ the unique element such that $F_{k+1} \setminus F_{k} = \{j_0\}$. Notice that $j_0 \notin F_1$ and $j_0 \notin J$. 

Let us define $I_1:=I$, $J_1 := J \cup \{j_0\}$, and 
$\mathcal{F}_1:=\left\{ F'_1 \subsetneq F'_2 \subsetneq \cdots  \subsetneq F'_{k} \subseteq E \setminus J_1 \right\}$ 
where $F'_{l} :=F_{l}$ for each $1 \leq l \leq k$.  
Let us show that $(I_1,\mathcal{F}_1,J_1)$ is a compatible triple, considering the cases $k=1$ and $k > 1$ separately. 
If $k=1$, then $\mathcal{F}_1:=\left\{ F'_1 \subseteq E \setminus J_1 \right\}$ and we must show that for each $j \in J_1$ the fiber $\pi^{-1}(j)$ is not contained in $I_1=I$. 
But this holds by the compatibility condition of the compatible triple $(I,\mathcal{F},J)$ since $J_1 = J \cup \{j_0\} \subseteq E \setminus F_1$. 
If $k > 1$, we must show that for each $j \in E \setminus F'_1$ the fiber $\pi^{-1}(j)$ is not contained in $I_1=I$. 
But this holds by the compatibility condition of the compatible triple $(I,\mathcal{F},J)$ since $F'_1=F1$ and $I_1=I$. 
Therefore $(I_1,\mathcal{F}_1,J_1)$ is a compatible triple. Since $|F'_{k}| = |F_k| = |F_{k+1}| -1 = s$, then the cone $C_{I_1,\mathcal{F}_1,J_1}$ is in $\Delta_{\pi,s}$. 

We have that $|J_1| = |J| + 1 = |E| - |F_{k+1}| + 1 = |E| - (s+1) +1 = |E| - s$. 
Then, by Lemma~\ref{technical.lemma.cone.containments}(\ref{technical.lemma.cone.containments.part.2}) the cone $C_{I_1,\mathcal{F}_1,J_1}$ contains exactly one of the cones $C_{J_0}$ with $J_0 \subseteq E$ such that $|J_0|=|E| -s$, namely for $J_0 = J_1 = J \cup \{ j_0 \}$. 
Therefore, the effect of the star subdivisions that on $\Delta_{\pi,s}$ that yield $\Sigma_{s+1}$ on the cone $C_{I_1,\mathcal{F}_1,J_1}$ in $\Delta_{\pi,s}$ is to subdivide it by the star subdivision along $C_{J_1}$.    
The minimal integral generators of $C_{I_1,\mathcal{F}_1,J_1}$ are 
\begin{align*}
    &\{e_i: i \in I_1 \} \cup \{-e_{A \setminus \pi^{-1}(F'_1)},\ldots,-e_{A \setminus \pi^{-1}(F'_{k-1})} \} \cup \{-e_{\pi^{-1}(j)} : j \in J_1 \} \\
    = &\{e_i: i \in I \} \cup \{-e_{A \setminus \pi^{-1}(F_1)},\ldots,-e_{A \setminus \pi^{-1}(F_{k-1})} \} \cup \{-e_{\pi^{-1}(j)} : j \in J \} \cup \{ -e_{\pi^{-1}(j_0)} \}
\end{align*}
and $C_{J_1}$ is minimally generated by $\{-e_{\pi^{-1}(j)} : j \in J_1 \} =\{-e_{\pi^{-1}(j)} : j \in J \} \cup \{ -e_{\pi^{-1}(j_0)} \}$. 
Therefore, after performing the star subdivisions on $\Delta_{\pi,s}$ that yield $\Sigma_{s+1}$, we obtain a cone $\tau_1$ in $\Sigma_{s+1}$ whose set of generators is obtained by taking the minimal set of integral generators of $C_{I_1,\mathcal{F}_1,J_1}$, and removing the vector $-e_{\pi^{-1}(j_0)}$ and including the vector 
\begin{align*}
        \sum_{j \in J_1} -e_{\pi^{-1}(j)} 
&= \sum_{j \in J \cup \{j_0\}} -e_{\pi^{-1}(j)}
= \sum_{j \in (E \setminus F_{k+1}) \cup \{j_0\}} -e_{\pi^{-1}(j)}
= \sum_{j \in E \setminus (F_{k+1} \setminus \{j_0\})} -e_{\pi^{-1}(j)}  \\
&= \sum_{j \in E \setminus F_{k}} -e_{\pi^{-1}(j)}
= -\sum_{i \in \pi^{-1}(E \setminus F_{k})} e_i 
= -\sum_{i \in A \setminus \pi^{-1}(F_{k})} e_i 
= -e_{A \setminus \pi^{-1}(F_{k})}. 
\end{align*}
Then, the cone $\tau_1$ in $\Sigma_{s+1}$ is generated by
\[
\{e_i: i \in I \} \cup \{-e_{A \setminus \pi^{-1}(F_1)},\ldots,-e_{A \setminus \pi^{-1}(F_{k})} \} \cup \{-e_{\pi^{-1}(j)} : j \in J \}. 
\]
Hence, by the definition of the cone associated to a compatible triple,  $C_{I,\mathcal{F},J} = \tau_1$. 
Therefore, $C_{I,\mathcal{F},J}$ is in $\Sigma_{s+1}$, and this completes the proof. 
\end{proof}

\begin{remark} \label{remark.toric.strict.ransform}
    Let $X$ be a toric variety with fan $\Delta$, and let $Y$ and $Z$ be torus invariant subvarieties corresponding to the cones $\sigma_Y$ and $\sigma_Z$ in $\Delta$.  
The fan $\widetilde{\Delta}$ of the toric variety $\widetilde{X} = \operatorname{Bl}_Z X$ is given by the star subdivision of $\Delta$ along the cone $\sigma_Z$. 
If we assume that $\sigma_Y$ does not contain $\sigma_Z$, then the cone $\sigma_Y$ is in $\widetilde{\Delta}$ and its corresponding torus invariant subvariety in $\widetilde{X}$ is the strict transform $\widetilde{Y}$ of $Y$ under the blow-up map $\widetilde{X} \rightarrow X$. 
Indeed, first notice that $\sigma_Y$ is in $\widetilde{\Delta}$ by the definition of the star subdivision operation and that by assumption $Y$ is not contained in $Z$. 
Since $\Delta_Y$ maps to itself under the map of fans $\widetilde{\Delta} \rightarrow \Delta$, 
then torus invariant subvariety in $\widetilde{\Delta}$ corresponding to $\Delta_Y$ has the same dimension as $Y$ and maps surjectively onto $Y$ under $\widetilde{X} \rightarrow X$. But, these properties characterize the strict transform $\widetilde{Y}$ of $Y$ among the subvarieties of $\widetilde{X}$.
\end{remark}

%%%%%%%%%%%%%%%%%%%%%%%%%%%%%%%%%%%%%%%%%%%%%%%%%%%%%%%%%%%%%%%%%%%%%%%%%%%%%%%%%

Let $\pi: A \rightarrow E$ be a surjective function between nonempty finite sets. 
We now present a construction of the polystellahedral variety $PS_{\pi}$ as an explicit blow-up of the product of projective spaces $\prod_{j \in E} \mathbb{P}^{|\pi^{-1}(j)|}$.
We may assume without losing generality that $E = \{1,2,\ldots,n\}$, for some $n \in \mathbb{Z}^{+}$, and we define $a_j := |\pi^{-1}(j)|$ for each $j \in E$. 
In this way $PS_{\pi} = PS_{[a_1,\ldots,a_n]}$ is the polystellahedral variety with cage $[a_1,\ldots,a_n]$ and $\prod_{j \in E} \mathbb{P}^{|\pi^{-1}(j)|} = \mathbb{P}^{a_1} \times \cdots \times \mathbb{P}^{a_n}$.

%%%%%%%%%%%%%%%%%%%%%%%%%%%%%%%%%%%%%%%%%%%%%%%%%%%%%%%%%%%%%%%%%%%%%%%%%%%%%%%%%

\begin{theorem}\label{PS: construction as blow-up of Pa}
The polystellahedral variety $PS_{\pi}=PS_{[a_1,\ldots,a_n]}$ is equal to the iterated blow-up of the product of projective spaces $\mathbb{P}^{a_1} \times \cdots \times \mathbb{P}^{a_n}$ along the subvarieties $H_J := \cap_{j \in J} H_j$ and their strict transforms for all nonempty $J \subseteq \{1,2,\ldots,n\}$, in any order such that $|J|$ is nondecreasing, where for each $1 \leq j \leq n$, $H_j$ is the pullback to $\mathbb{P}^{a_1} \times \cdots \times \mathbb{P}^{a_n}$ of a hyperplane in $\mathbb{P}^{a_j}$ under the projection~map.
\end{theorem}

%%%%%%%%%%%%%%%%%%%%%%%%%%%%%%%%%%%%%%%%%%%%%%%%%%%%%%%%%%%%%%%%%%%%%%%%%%%%%%%%%

\begin{proof}
Let $\widetilde{\Sigma}_{0}$ be the standard fan of $\mathbb{P}^{a_1} \times \cdots \times \mathbb{P}^{a_n}$ in $\mathbb{R}^{A} := \mathbb{R}^{a_1} \times \cdots \times \mathbb{R}^{a_n}$. 
We can assume that $H_j$ is the pull back to $\mathbb{P}^{a_1} \times \cdots \times \mathbb{P}^{a_n}$ of the torus invariant divisor in $\mathbb{P}^{a_j}$ associated to the ray spanned by $(-1,\ldots,-1) \in \mathbb{R}^{a_i}$, for each $1 \leq j \leq n$. 
We define $X_0 = \mathbb{P}^{a_1} \times \cdots \times \mathbb{P}^{a_n}$ and for each $1 \leq s \leq n$ we inductively define the variety $X_s$ to be the toric variety obtained by blowing up $X_{s-1}$ along the strict transforms of the varieties $H_J$ for all $J \subseteq \{1,2,\ldots,n\}$ with $|I| = n-s+1$ in some particular order that we chose and fix (we will see below that the variety $X_s$ is independent of this order). 
For each $1 \leq s \leq n$ let $\widetilde{\Sigma}_{s}$ be the fan of the toric variety $X_s$. 

Since $\Delta_{\pi,0}$ is the fan of $\mathbb{P}^{a_1} \times \cdots \times \mathbb{P}^{a_n}$ by Example~\ref{example.fan.product.projective.spaces} and $\Delta_{\pi,n}$ is the fan of $PS_{\pi} = PS_{[a_1,\ldots,a_n]}$ by Example~\ref{example.fan.polystellahedron}, then the theorem follows if we show inductively that for each $0 \leq s \leq n$ the fan $\widetilde{\Sigma}_{s}$ agrees with the fan $\Delta_{\pi,s}$ (see Definition~\ref{definition.associated.fans.triple}).  
Notice that this would also give us that iterated blow-ups to obtain $X_s$ from $X_{s-1}$ can be performed in any order. For $s=0$, we have that $\widetilde{\Sigma}_{0}$ coincides with $\Delta_{\pi,0}$ by Example~\ref{example.fan.product.projective.spaces}.
Let us assume that $\widetilde{\Sigma}_{s}$ coincides with $\Delta_{\pi,s}$ for some $0 \leq s \leq n-1$, and let us show that $\widetilde{\Sigma}_{s+1}$ coincides with $\Delta_{\pi,s+1}$. 

For each $J \subseteq \{1,2,\ldots,n\}$ with $|J|= n-s$, let $\sigma_J$ be the cone in $\widetilde{\Delta}_s$ corresponding to the strict transform of the torus invariant subvariety $H_J$. 
Notice that for each $J \subseteq \{1,2,\ldots,n\}$ with $|J|= n-s$ the cone $\sigma^0_J$ in $\widetilde{\Sigma}_0$ corresponding to $H_J$ does not contain any of the cones along which we performed star subdivisions in the previous steps. 
Then, then $\sigma_J = \sigma^0_J$ by Remark~\ref{remark.toric.strict.ransform}. 
Hence, the fan $\widetilde{\Sigma}_{s+1}$ is obtained as the iterated star subdivision of $\widetilde{\Sigma}_{s}$ along the all cones $\sigma_J =\sigma^0_J$ for $J \subseteq \{1,2,\ldots,n\}$ with $|J|= n-s$, in the order that we fixed before. 
On the other hand, by Lemma~\ref{technical.lemma.blow-up.triples} the fan $\Delta_{\pi,s+1}$ is the iterated star subdivision of the fan $\Delta_{\pi,s}$ along the cones $C_J$ in Definition~\ref{definition.cones.CJ}, for all $J \subseteq E = \{1,2,\ldots,n\}$ with $|J|=|E|-s = n-s$, in any order.
Therefore, to conclude the proof it is enough to show that $\sigma_J = C_J$ for all $J \subseteq E$ with $|J| = |E| - s = n-s$. 

Let $J \subseteq E = \{1,2,\ldots,n\}$ with $|J| = |E| - s = n-s$. 
Let $A_j := \{1,2,\ldots, a_i \}$ for each $1 \leq j \leq n$ and let $\{ e_{ji} \, | \, 1 \leq  j  \leq n, i \in A_j  \}$ denote the canonical basis of $\mathbb{R}^A = \mathbb{R}^{a_1} \times \cdots \times \mathbb{R}^{a_n}$. 
In the present notation, the cone $C_J$ is the cone generated by 
$\{  - \sum_{i \in A_j} e_{ji} \, | \,  j \in J  \}$. 
For each $j \in J$, the cone in $\widetilde{\Sigma}_0$ corresponding to $H_{j}$ is the ray generated by $\{ - \sum_{i \in A_j} e_{ji} \}$ in $\mathbb{R}^A$. 
Since $H_J := \cap_{j \in J} H_j$, then the cone $\sigma_J$ is the sum of the cones in $\mathbb{R}^A$ corresponding to the varieties $H_j$ for $j \in J$. 
Then, $H_J$ is the cone in $\mathbb{R}^A$ generated by $\{  - \sum_{i \in A_j} e_{ji} \, | \,  j \in J  \}$. 
Therefore $\sigma_J = C_J$ and the proof is complete. 
\end{proof}

\subsection{The polypermutohedral variety as a quotient of the polystellahedral variety} \label{subsection.quotients}
In this subsection we prove that the polypermutohedral variety is a geometric quotient of an open subset of the polystellahedral variety in a canonical way.

\begin{lemma} \label{toric.divisor.is.quotient.of.toric.variety}
    Let $ X $ be a smooth toric variety with torus $ T $, and let $ Z $ be a $ T $-invariant subvariety of $ X $. Let $ U $ denote the complement of the union of $ T $-invariant subvarieties of $ X $ that do not intersect $ Z $, and let $ E $ be the exceptional divisor in the blow-up of $ X $ along $ Z $. Then, there exists a canonical geometric quotient:
    \[
    (U \setminus Z) \sslash \mathbb{G}_m \cong E.
    \]
\end{lemma}

\begin{proof}
We show that there exists a canonical morphism of toric varieties $\pi : U \to E $ such that the restriction $ \pi : (U \setminus Z) \to E $ is a geometric quotient for the action of the pointwise stabilizer $\operatorname{Stab}_T(E) \cong \mathbb{G}_m \subseteq T $ on $ U $. 
The desired conclusion for a particular fan also yields the analogous conclusion for a subfan, then we may assume that the maximal cones in the fan of $X$ are top dimensional. 

Let $U_E$ denote the complement of the union of $T$-invariant subvarieties of $\operatorname{Bl}_Z X$ that do not intersect $E$.
Then, the blow-up morphism gives the identification $U_E \setminus E = U \setminus Z$. 
Therefore, we can assume for the rest of the proof that 
$Z=D$ is a $T$-invariant prime divisor on $X$ and $E=D$. 

Let $\Sigma$ be the fan of $X$ in $N_{\mathbb{R}} = N \otimes \mathbb{R}$ and $\rho \in \Sigma$ be the ray corresponding to $D$. 
Let us recall the structure of $D$ as a toric variety; see \cite[Proposition 3.2.7.]{Cox-Little-Schenck-Toric-Varieties-Book} for details.  
Let $N_\rho$ be the sublattice of $N$ spanned by the points in $\rho \cap N$, and let $N(\rho)=N / N_\rho$.
For each cone $\sigma \in \Sigma$ containing $\rho$, let $\overline{\sigma}$ be the image cone in $N(\rho)_{\mathbb{R}}$ under the quotient map
$N_{\mathbb{R}} \longrightarrow N(\rho)_{\mathbb{R}}$. 
We consider the collections of cones   
\[
\operatorname{Star}_X(\rho)=\{ \sigma \subseteq N_{\mathbb{R}} \mid \rho \preceq \sigma \in \Sigma\} 
\textnormal{\quad and \quad }
\operatorname{Star}(\rho)=\{\overline{\sigma} \subseteq N(\rho)_{\mathbb{R}} \mid \rho \preceq \sigma \in \Sigma\}
\]
Then, by \cite[Proposition 3.2.7.]{Cox-Little-Schenck-Toric-Varieties-Book}, 
$\operatorname{Star}(\rho)$ is a fan in $N(\rho)_{\mathbb{R}}$ and $D$ is the toric variety $X_{\operatorname{Star}(\rho)}$. 
Moreover, the assignment $\sigma \longmapsto \overline{\sigma}$ is a bijection between the cones in $\operatorname{Star}_X(\rho)$ and those in $\operatorname{Star}(\rho)$.
We have that 
\[
\Sigma' := \{\sigma \in \Sigma \mid \exists \tau \in \Sigma \textnormal{ such that } \rho \preceq \tau \textnormal{ and } \sigma \preceq \tau \}
\]
is a subfan of $\Sigma$ and $U$ is the toric variety $X_{\Sigma'}$.
Similarly, we have that 
\[
\Sigma'' := \{\sigma \in \Sigma \mid 
\rho \not\preceq \sigma
\textnormal{ but }
\exists \tau \in \Sigma \textnormal{ such that } \rho \preceq \tau \textnormal{ and } \sigma \preceq \tau \}
\]
is a subfan of $\Sigma$ and $U \setminus D$ is the toric variety $X_{\Sigma''}$.

The quotient map $N_{\mathbb{R}} \longrightarrow N(\rho)_{\mathbb{R}}$ induces morphisms of toric varieties $U = X_{\Sigma'} \rightarrow D$ and $(U \setminus D) = X_{\Sigma''} \rightarrow D$. 
We claim that these morphisms are respectively a good quotient and a geometric quotient for the action of $\mathbb{G}_m$. 
These claims can be verified locally on the base. 

The assignment $\sigma \mapsto \sigma + \rho$ is a bijection between the cones in the fan $\Sigma''$ and the cones in $\operatorname{Star}_X(\rho)$.
Hence, the assignment $\sigma \mapsto \overline{(\sigma + \rho)}$ is a bijection between the cones in the fan $\Sigma''$ and the cones in $\operatorname{Star}(\rho)$.
We also notice that any $\sigma$ in $\Sigma''$ and $\sigma + \rho$ both map onto $\overline{(\sigma + \rho)}$ under the projection $N_{\mathbb{R}} \longrightarrow N(\rho)_{\mathbb{R}}$. 

Since the toric variety $X$ is smooth, all toric varieties under consideration are 
smooth, and hence we can choose coordinates and reduce to the following affine setting. 
Let $\tau$ be a fixed maximal cone in $\operatorname{Star}(\rho)$, which we know has the form $\overline{(\sigma + \rho)}$ for a unique cone $\sigma$ in $\Sigma''$.   
We can choose affine coordinates such that $U_{\sigma+\rho} = \operatorname{Spec} \mathbb{C}[x_0,x_1,\ldots,x_n]$, $D \cap U_{\sigma+\rho} = Z(x_0) \subseteq U_{\sigma + \rho}$, and $U_{\sigma} = \operatorname{Spec} \mathbb{C}[x_0,x_1,\ldots,x_n]_{x_0}$. 
We are reduced to considering the action of $\mathbb{G}_m$ on $\mathbb{C}[x_0,x_1,\ldots,x_n]$ via $t \cdot f(x_0,x_1,\ldots,x_n) = f(tx_0,x_1,\ldots,x_n)$, for any $t \in \mathbb{G}_m$ and $f \in \mathbb{C}[x_0,x_1,\ldots,x_n]$.
We must verify that the ring homomorphisms
$\mathbb{C}[x_0,x_1,\ldots,x_n] \rightarrow \mathbb{C}[x_1,\ldots,x_n]$
and 
$\mathbb{C}[x_0,x_1,\ldots,x_n]_{x_0} \rightarrow \mathbb{C}[x_1,\ldots,x_n]$
respectively induce a good quotient $U_{\sigma+\rho} \rightarrow U_{\overline{\sigma+\rho}}$
and a geometric quotient 
$U_{\sigma} \rightarrow U_{\overline{\sigma+\rho}}$
for this action, but this clearly holds since 
$\mathbb{C}[x_0,x_1,\ldots,x_n]^{\mathbb{G}_m} = \mathbb{C}[x_1,\ldots,x_n]$
and 
$\mathbb{C}[x_0,x_1,\ldots,x_n]_{x_0}^{\mathbb{G}_m} = \mathbb{C}[x_1,\ldots,x_n]$
clearly satisfy the required properties.

\end{proof}

\begin{corollary}  \label{corollary.toric.quotient}
    Let $ X $ be a smooth toric variety with torus $ T $, and let $ D $ be a $ T $-invariant prime divisor on $ X $. Let $ U $ denote the complement of the union of $ T $-invariant subvarieties of $ X $ that do not intersect $ D $. 
    Then, there exist a good quotient and a geometric quotient: 
   \[U / {\mathbb{G}_{m}} \cong D  \textnormal{\quad and \quad}  (U \setminus D) \! \sslash \! {\mathbb{G}_{m}} \cong D.\]
\end{corollary}

\begin{proof}
We established these quotients within the proof of  Lemma~\ref{toric.divisor.is.quotient.of.toric.variety} in the case $Z=D$.
\end{proof}

%%%%%%%%%%%%%%%%%%%%%%%%%%%%%%%%%%%%%%%%%%%%%%%%%%%%%%%%%%%%%%%%%%%%%%%%%%%

We achieve the goal of this subsection by relating the polystellahedral and polypermutohedral varieties via a geometric quotient in the following proposition.

\begin{proposition} \label{polypermuta.is.quotient.of.polystella}
    Let $\pi: A \rightarrow E$ be a surjective function between finite sets and let $\mathbf{a} = (a_1, \ldots, a_n) \in \mathbb{Z}_{>0}^n$ be the associated cage. 
    Then, there exists a nonempty open subset $\left({PS}_{\mathbf{a}}\right)^{\circ}$ of the polystellahedral variety ${PS}_{\mathbf{a}}$ with a canonical $\mathbb{G}_m$ action, which admits as a geometric quotient the polypermutohedral variety $P_{\mathbf{a}}$, that is,  
        \[
            P_{\mathbf{a}}  \cong
             \left( {PS}_{\mathbf{a}}  \right)^{\circ}\quotient \mathbb{G}_m.
        \]
\end{proposition}

\begin{proof}
    By \cite[Section 6]{eur2024intersection}, there is an embedding of $P_{\mathbf{a}}$ as a torus invariant divisor in $PS_{\mathbf{a}}$. 
    The result follows by applying Corollary~\ref{corollary.toric.quotient} to this embedding $P_{\mathbf{a}} \hookrightarrow PS_{\mathbf{a}}$. 
    Explicitly, since the polypermutohedron is the facet of the polystellahedron corresponding to the ray spanned by $e_A=(1,\ldots,1) \in \mathbb{R}^{A}$, then $\left({PS}_{\mathbf{a}}\right)^{\circ}$ is defined by the subfan of the polystellahedral fan consisting of the cones which do not contain $e_A$ but together with $e_A$ span a cone. 
\end{proof}

%%%%%%%%%%%%%%%%%%%%%%%%%%%%%%%%%%%%%%%%%%%%%%%%%%%%%%%%%%%%%%%%%%%%%%%

\subsection{Proof of Theorem \ref{thm:mainPolyStella}} \label{subsection.proof.main.thm.polystella}
To conclude this section, we provide a proof of Theorem~\ref{thm:mainPolyStella}. The different parts of these results are spread throughout the previous subsections, so this proof is meant to aid the reader in navigating the previous results.

\begin{proof}[Proof of Theorem~\ref{thm:mainPolyStella}]
The variety  $\mathbb{P}^{[\mathbf{a}]}_{H}$ is defined to be a wonderful compactification with respect to the building set from Lemma \ref{polystella: lemma pullbackBuildingsetGFM}. Since $H$ is part of the toric boundary of $\mathbb{P}^d$, the building set $\mathcal{G}_{\mathbf{a}}$ is contained in the toric boundary of $\mathbb{P}^{\mathbf{a}}$. Therefore, Claim (i) follows from Proposition \ref{prop: li blow-up building set}. 
The varieties $\mathbb{P}^{[\mathbf{a}]}_{H}$ and ${PS}_{\mathbf{a}}=PS_{[a_1, \ldots, a_n]}$ are both iterated blow-ups of 
$\mathbb{P}^{a_1} \times \cdots \times \mathbb{P}^{a_n}$.  
 Therefore,  Claim (ii) follows by showing that the centers of these blow-ups are the same.  However, this is clear from comparing the blow-up centers given by the building set $\mathcal{G}_{\mathbf{a}}$ with the centers described in
 Theorem~\ref{PS: construction as blow-up of Pa}. 
 Lastly, since we have established that $T^{\mathbf{a}}_{LM}$ is isomorphic to the polypermutohedral variety $P_{\mathbf{a}}$ in Theorem~\ref{GFM: thm:weightedFM_forFlags}(i) and $\mathbb{P}^{[\mathbf{a}]}_{H}$ is isomorphic to the polystellahedral variety ${PS}_{\mathbf{a}}$ in Theorem~\ref{thm:mainPolyStella}(ii), then Claim (iii) is given by Proposition~\ref{polypermuta.is.quotient.of.polystella}.
\end{proof}

%\bibliographystyle{alpha}
%\bibliography{bibliography}

\newcommand{\etalchar}[1]{$^{#1}$}

\end{document}